\documentclass[11pt, a4paper]{article}

\usepackage{graphicx} % Required for inserting images
\usepackage{url}
\usepackage{a4wide}
\usepackage{natbib}
\usepackage{enumitem}
\usepackage{hyperref}
\usepackage{setspace}
\usepackage{multirow, booktabs}
\usepackage{amsmath, amsthm, amssymb, amsfonts}
\usepackage{bm, bbm}
\usepackage[ruled, vlined]{algorithm2e}
\usepackage{color}
\usepackage[dvipsnames]{xcolor}
\usepackage{csquotes}
\usepackage{cleveref}
\usepackage{makecell}
\usepackage{xr}

\allowdisplaybreaks
\DeclareMathAlphabet\mathbfcal{OMS}{cmsy}{b}{n}

\newcommand{\mbf}{\mathbf}
\newcommand{\mc}{\mathcal}
\newcommand{\bmx}{\begin{bmatrix}}
\newcommand{\emx}{\end{bmatrix}}

\newcommand{\eps}{\epsilon}
\newcommand{\vep}{\varepsilon}

\renewcommand{\l}{\left}
\renewcommand{\r}{\right}

\def\wh{\widehat}
\def\wt{\widetilde}

\newcommand{\cp}{\theta}

\newcommand{\E}[0]{\mathsf{E}}
\newcommand{\Var}[0]{\mathsf{Var}}
\newcommand{\Cov}[0]{\mathsf{Cov}}
\newcommand{\tr}[0]{\mathsf{tr}}
\newcommand{\p}{\mathsf{P}}
\newcommand{\R}{\mathbb{R}}

\newcommand{\iid}{\text{\upshape iid}}
\newcommand{\nn}{\nonumber}

\newcommand{\Mc}{\text{\tiny \upshape M}} % subscript for McScan
\newcommand{\Qc}{\text{\tiny \upshape Q}} % subscript for QcScan
\newcommand{\Oc}{\text{\tiny \upshape O}} % subscript for OcScan
\newcommand{\Rf}{\text{\tiny \upshape R}} % subscript for Refined Scan

\theoremstyle{definition}
\newtheorem{thm}{Theorem}
\theoremstyle{definition}
\newtheorem{cor}[thm]{Corollary}
\theoremstyle{definition}
\newtheorem{lem}[thm]{Lemma}
\theoremstyle{definition}
\newtheorem{prop}[thm]{Proposition}
\theoremstyle{definition}

\theoremstyle{remark}
\newtheorem{rem}{Remark}[section]
\theoremstyle{definition}

\theoremstyle{definition}

\theoremstyle{definition}

\title{Covariance scanning for adaptively optimal change point detection in high-dimensional linear models}
\author{Haeran Cho \and Housen Li}
\date{\today}

\begin{document}

\maketitle

\begin{abstract}
This paper investigates the detection and estimation of a single change in high-dimensional linear models. We derive the minimax lower bound for the detection boundary which uncovers a phase transition governed by the sparsity of the covariance-weighted differential parameter. This form of {\it inherent sparsity} captures a delicate interplay between the covariance structure of the regressors and the change in regression coefficients on the detectability of a change point. Complementing the lower bound, we introduce two covariance scanning-based methods, McScan and QcSan, which achieve minimax optimal detection performance (up to possible logarithmic factors) in the sparse and the dense regimes, respectively. In particular, QcScan is the first method to achieve consistency in the dense regime and further, we devise a combined procedure which is adaptively minimax optimal across sparse and dense regimes without the knowledge of the sparsity. Computationally, covariance scanning-based methods avoid costly computation of Lasso-type estimators and attain worst-case computation complexity that is linear in the dimension and sample size. Additionally, we consider the post-detection estimation of the differential parameter and the refinement of the change point estimator. Simulation studies support the theoretical findings and demonstrate the computational and statistical efficiency of the proposed covariance scanning methods.
\end{abstract}

%Goal: A complete scenario of minimax optimality in change point detection in high-dimensional linear regressions.
%\begin{itemize}
%\item 
%$\Sigma = I$, Gaussian, Single change point
%\item 
%Using concentration inequality of Gaussian chaos
%\item 
%Detection rates adaptive to sparsity
%\item 
%Localisation rates adaptive to sparsity
%\item 
%What if $\Psi$ is large?
%\end{itemize}

\section{Introduction}

The change point problem in linear models has received great attention over the past several decades, owing to its widespread applicability across diverse scientific and practical disciplines, including genomics, climatology, finance and economics. Numerous methods have been developed that achieve consistency in (multiple) change point detection, both in multivariate settings (e.g.\ \citealp{csorgo1997limit, bai1998estimating}) and high-dimensional scenarios (e.g.\ \citealp{lee2016lasso, kaul2019detection, wang2021statistically, rinaldo2020localizing, cho2022high, xu2024change}). 

When faced with high dimensionality, many existing methods %, such as MOSEG \citep{cho2022high}, 
rely on $\ell_1$-regularised estimators of local regression parameters for change point detection, thereby requiring the sparsity of regression parameters. 
\cite{gao2022sparse} propose a projection-based method that instead assumes sparsity on the \emph{differences} of regression parameters. However, it is applicable only when the dimension is strictly smaller than the sample size. 
More recently, \cite{cho2024detection} introduce McScan, a multiscale covariance scanning approach that achieves consistent change point detection with both statistical and computational efficiency. McScan does not require exact sparsity of the regression parameters or their differences, and accommodates ultra high-dimensional settings. It operates by examining the $\ell_\infty$-norm of the contrasts in local sample covariances between the response and the regressors, and thus is better suited to sparse regimes implicitly. 
We also remark that minimax lower bounds for the change point problem have been investigated, but only when the sparsity of the regression parameters (or their differences) is bounded by the square root of the dimension (see \citealp{rinaldo2020localizing} and \citealp{cho2024detection}). 

Despite these advances, several questions remain to be answered for the change point problem in high-dimensional linear models. In particular, it is hitherto unknown whether it is possible to consistently detect the change point when the change is \emph{non-sparse}. More fundamentally, it is yet to be understood which notion of sparsity --- or more generally, which structural assumption --- determines the \emph{intrinsic difficulty} of the change point problem. 
We provide first answers to these questions under a Gaussian linear model with at most one change (AMOC): The observations $(Y_t, \mbf x_t) \in \R^{1+p}, \, t \in [n]:= \{1, \ldots, n\}$, are distributed as 
\begin{align}
\label{eq:amoc}
Y_t = \begin{cases}
\mbf x_t^\top \bm\beta_0 + \vep_t & \text{for \ } 1 \le t \le \cp, \\
\mbf x_t^\top \bm\beta_1 + \vep_t & \text{for \ } \cp + 1 \le t \le n, 
\end{cases}
\end{align}
where $\mbf x_t \sim_{\iid} \mc N_p(\mbf 0, \bm\Sigma)$ with non-negative definite $\bm\Sigma \in \R^{p \times p}$ and independently, $\vep_t \sim_{\iid} \mc N(0, \sigma^2)$ for all~$t$. We denote by $q \in \{0, 1\}$ the number of change points. In the case of a single change point ($q = 1$), we define the differential parameter as $\bm\delta = \bm\beta_1 - \bm\beta_0 \in \R^p$, the maximum variability of the response as $\Psi = \max(\vert \bm\Sigma^{1/2} \bm\beta_0 \vert_2, \vert \bm\Sigma^{1/2} \bm\beta_1 \vert_2, \sigma)$, 
and the minimal segment length as $\Delta = \min(\cp, n - \cp)$. In the absence of any change point ($q = 0$), where $\bm\beta_0 = \bm\beta_1$, we write $\cp = n$ as a convention.

Under the model~\eqref{eq:amoc}, we are interested in the problem of consistently detecting the presence of any change point and, if $q = 1$, estimating its location by a consistent estimator $\wh\cp$ satisfying $\vert \wh\cp - \cp \vert / \Delta \to 0$ with probability tending to one as $n \to \infty$.
By focusing on the AMOC situation, we separate the two \enquote{orthogonal} aspects of the problem: The design of an adaptively optimal procedure for change point detection, which is the main focus of this paper, and the development of algorithms for multiple change point detection. While our theoretical analysis is framed under independence and Gaussianity, we also discuss extensions to accommodate temporal dependence or non-Gaussianity. 

We summarise our main contributions to the growing literature on change point analysis in high-dimensional linear models, which also serves to outline the structure of the paper.
\begin{enumerate}[label = (\roman*)]
\item \emph{Phase transition and a new notion of sparsity.} 
We derive the minimax lower bounds for the change point detection as summarised in Table~\ref{tab:minimax}. With the (nearly) matching upper bounds presented therein, it reveals a phase transition governed by $\mathfrak{s} := \vert\bm\Sigma^{1/2} \bm\delta\vert_0$, which we term the \emph{inherent sparsity}, rather than the standard notion of sparsity $\mathfrak{s}_\delta := \vert\bm\delta\vert_0$, where $\vert \cdot \vert_0$ denotes the $\ell_0$-semi-norm.  
% Notably, the results also indicate that in the dense regime, where $\mathfrak{s} \ge \sqrt{p\log\log(n)}$, consistent estimation is possible only when $p \ll n^2$. 
See \Cref{sec:minimax}.

\begin{table}[t!]
\caption{Minimax lower bound on $\Psi^{-2} \vert \bm\Sigma^{1/2} \bm\delta \vert_2^2 \Delta$, for the detection of a single change point under the model~\eqref{eq:amoc}, corresponding upper bounds and the methods that attain the upper bounds. Here $\sigma_*$ denotes the smallest \textit{non-zero} eigenvalue of the covariance matrix $\bm\Sigma$ of regressors, $\sigma_X^2$  the largest diagonal entry of $\bm\Sigma$, and $\mathfrak{s} = \vert \bm\Sigma^{1/2} \bm\delta \vert_0$ the inherent sparsity. When $\bm\Sigma$ is known, the dependence on $\sigma_*$, $\sigma_X^2$ and $\Vert \bm\Sigma \Vert$ can be omitted. The lower bound is derived under~\eqref{e:spn} supposing that $\operatorname{rank}(\bm\Sigma) \asymp p$. For the upper bound in the regime of $\mathfrak{s} \ge \sqrt{p{\log\log(n)}}$, we suppose that $p \ge n$; an analogous result holds for $p < n$, see Remark~\ref{r:qcscan}~\ref{i:r:np}.}
\label{tab:minimax}
\centering
\resizebox{\linewidth}{!}{
\begin{tabular}{c cc c}
\toprule
% & \multicolumn{2}{c}{Detection requisite on $\Psi^{-2} \vert \bm\Sigma^{1/2} \bm\delta \vert_2^2 \Delta$} & \\
% \cmidrule(lr){2-3} 
Regime & Lower bound {(\Cref{th:gaussid})}  & Upper bound {(Theorems~\ref{thm:mcscan} and \ref{thm:qcscan})} & Method \\ 
\cmidrule(lr){1-1} \cmidrule(lr){2-2} \cmidrule(lr){3-3} \cmidrule(lr){4-4} 
{$\mathfrak{s} < \sqrt{p\log\log(n)}$} & 
{${\mathfrak{s}}\log(\tfrac{e p\log\log(n)}{\mathfrak{s}^2})$} 
& 
{$\tfrac{\sigma_X^2}{\sigma_*} \mathfrak{s}\log(p\log(n))$} 
& 
{McScan (Section~\ref{sec:sparse})}
\\
\cmidrule(lr){1-1} \cmidrule(lr){2-2} \cmidrule(lr){3-3} \cmidrule(lr){4-4}
{$\mathfrak{s} \ge \sqrt{p\log\log n}$}
& 
{$ \sqrt{p\log\log(n)}$} 
& 
{$\tfrac{ \Vert \bm\Sigma \Vert }{\sigma_*} \sqrt{p \log\log(n)}$}
& {QcScan (Section~\ref{sec:dense})} 
\\
\bottomrule
\end{tabular}}
\end{table}

\item \emph{The first method for consistent change point detection when $\bm\Sigma$ is rank deficient} and in the dense regime.
We first show that McScan \citep{cho2024detection}, which scans the $\ell_\infty$-aggregation of the contrasts in local sample covariance between the response and the regressors, is minimax optimal in the sparse regime up to the factor of $\sigma_*^{-1}\sigma_X^2$ and possible logarithmic factors in $p$ and $n$ (\Cref{sec:sparse}), where $\sigma_*$ denotes the smallest \textit{non-zero} eigenvalue of $\bm\Sigma$ and $\sigma_X^2$ its largest diagonal entry.
Additionally, we introduce the \emph{quadratic covariance scanning} (QcScan) that evaluates the $\ell_2$-aggregation of the local sample covariance contrasts, and show QcScan to be the first method that achieves minimax optimality, up to a factor of $\sigma_*^{-1} \Vert \bm\Sigma \Vert$, in the dense regime.
In short, both McScan and QcScan accommodate rank-deficient $\bm\Sigma$, which is in contrast to existing change point methods that require $\Cov(\mbf x_t)$ to have bounded spectrum.
As a by-product, we also obtain an estimator for $\vert\bm\Sigma\bm\delta\vert_2^2$, which serves as a surrogate of the strength of the signal $\bm\delta^\top\bm\Sigma\bm\delta$ for the change point problem in~\eqref{eq:amoc}. See \Cref{sec:dense}. 

% \item \emph{The first consistent change point detection method in the dense regime.}
% We first show that McScan \citep{cho2024detection}, which scans the $\ell_\infty$-aggregation of contrasts in local sample covariances between regressors and responses, is minimax optimal in the sparse regime, up to a factor of $\sigma_*^{-1}\sigma_X^2$ and possible logarithmic factors in $p$ and $n$. Here $\sigma_*$ denotes the smallest \emph{non-zero} eigenvalue of $\bm\Sigma$, and $\sigma_X^2$ its largest diagonal entry. 
% %In particular, McScan achieves optimality in the sparse regime even when $\bm\Sigma$ is rank deficient.
% For the dense regime, we introduce \emph{quadratic covariance scanning} (QcScan), which evaluates the $\ell_2$-aggregation of local sample covariance contrasts. QcScan is the first method that achieves minimax optimality, up to a factor of $\sigma_*^{-1}\|\bm\Sigma\|$, for  the change point detection in the dense regime. 
% {
% Both McScan and QcScan accommodate rank deficient $\bm\Sigma$, in contrast to existing change point methods, which typically require $\bm\Sigma$ to have bounded spectrum.}
% %
% As a by-product, we also obtain an estimator of $\lvert\bm\Sigma\bm\delta\rvert_2^2$, which serves as a surrogate for the signal strength $\bm\delta^\top\bm\Sigma\bm\delta$ in the change point problem~\eqref{eq:amoc}. See \Cref{sec:upper}.

\item \emph{Adaptation to inherent sparsity.} % To address both sparse and dense regimes, 
We introduce OcScan which, by combining McScan and QcScan, adapts to the unknown level of sparsity and thus achieves optimal statistical performance uniformly across all $\mathfrak{s}$, providing a complete treatment of the change point detection problem in high-dimensional linear models. See \Cref{sec:adaptive}. 

\item \emph{Refinement of change point localisation.} We propose a procedure for further refinement of the change point estimator by incorporating the direct estimation of the differential parameter alongside covariance scanning. This refinement can improve the estimation rate by {up to a factor of $\sqrt p$}, provided that $\bm\delta$ itself is sufficiently sparse. See \Cref{sec:refine}. 

\item \emph{Computational efficiency.} In contrast to a plethora of existing methods that require repeated estimation of local regression parameters, the proposed methods based on covariance scanning, namely McScan, QcScan and OcScan, exhibit ideal worst-case computation complexity of $O(np)$, see \Cref{ss:comp}. 
This demonstrates that for the detection and estimation of a change point, computationally costly operations may be avoided without losing statistical efficiency. 
\end{enumerate}

We further corroborate our theoretical findings with simulation studies in \Cref{sec:sim}, and provide all proofs and technical details in the Appendix.

\paragraph{Notations.}
For a matrix $\mbf A \in \R^{m \times n}$, we denote by $\Vert \mbf A \Vert$ and $\Vert \mbf A \Vert_F$ its operator and Frobenius norms, respectively.
For $0 \le a < b \le m$, we write $\mbf A_{a, b}$ to denote a sub-matrix of $\mbf A$ containing its $i$th rows for $a + 1 \le i \le b$; when $\mbf A$ is a vector, $\mbf A_{a, b}$ denotes the corresponding sub-vector.
By $\lvert\cdot\rvert_\nu$ we denote the $\ell_\nu$-(semi-)norm of a vector for $\nu \in[0, \infty]$.
For sequences $\{a_m\}$ and $\{b_m\}$ of positive numbers, we write $a_m \lesssim b_m$ or equivalently $a_m = O(b_m)$, if $a_m \le C b_m$ for some finite constant $C > 0$. If $a_m \lesssim b_m$ and $b_m \lesssim a_m$, we write $a_m \asymp b_m$. We write $a \vee b = \max(a, b)$ for $a, b \in \R$.
By $\mbf 0$ and $\mbf I$, we denote the vector of zeros and the identity matrix, respectively, whose dimensions depend on the context.

\section{Minimax lower bound}
\label{sec:minimax}

% Important papers: \citet{Bar02,VeVi10,VeGa18}

\subsection{Isotropic case}
\label{sec:isotropic}

Let us define a model space with a single change under~\eqref{eq:amoc}, as
\begin{multline*}
\mc P_{\bm\Sigma, \sigma}^{\mathfrak{s}, n, p}(\tau) = \biggl\{ \p_{\cp, \bm\beta_0, \bm\beta_1} \;:\; \{(Y_t, \mbf x_t)\}_{t \in [n]} \sim \p_{\cp, \bm\beta_0, \bm\beta_1} \text{ \ such that \ } \mbf x_t \sim_{\iid} \mc N_p(\mbf 0, \bm\Sigma) \text{ and}\\
\text{independently \ }  \vep_t \sim_{\iid}\mc N (0,\sigma^2), \; 
Y_t = \mbf x_t^\top \bigl( \bm\beta_0 \cdot \mathbb{I}_{\{t \le \cp\}} + \bm\beta_1 \cdot \mathbb{I}_{\{t > \cp\}} \bigr) + \vep_t,
\\
\bm\delta = \bm\beta_1 - \bm\beta_0, \,
\Delta = \min(\cp, n - \cp), \,
 \, \Psi = \max(\vert \bm\Sigma^{1/2} \bm\beta_0 \vert_2, \vert \bm\Sigma^{1/2} \bm\beta_1 \vert_2, \sigma), 
\\
\vert \bm\Sigma^{1/2} \bm\delta\vert_0 = \mathfrak{s}
\text{ \ and \ } \vert \bm\Sigma^{1/2}\bm\delta\vert_2^2 \Delta \ge \tau \Psi^2\mathfrak{s} \log\Bigl(1 + \tfrac{p {\log\log(n)}}{\mathfrak{s}^2}\vee\tfrac{\sqrt{p {\log\log(n)}}}{\mathfrak{s}}\Bigr)\biggr\},
\end{multline*} 
with some constant $\tau > 0$. We also define the model space without any change point as 
\begin{multline*}
\mc P_{\bm\Sigma, \sigma}^{0, n, p} = \biggl\{\p_{\bm\beta} \;:\; \{(Y_t, \mbf x_t)\}_{t \in [n]} \sim \p_{\bm\beta} \text{ \ such that \ } Y_t = \mbf x_t^\top \bm\beta + \vep_t  \text{ \ with \ }\\ \mbf x_t \sim_{\iid} \mc N_p(\mbf 0, \bm\Sigma)
 \text{ \  and independently, \ } \vep_t \sim_{\iid}\mc N (0,\sigma^2) 
 \biggr\}.
\end{multline*}

\begin{thm}
\label{th:gaussid}
Let $\sigma > 0$ and $n$ be sufficiently large. We consider the case where $\bm\Sigma = \mbf I$.
%\begin{enumerate}[label = (\roman*)]
%\item\label{lem:lb:id:detect}
%\emph{Detection:}
If 
\begin{equation}\label{e:spn}
\mathfrak{s}\log\l(1 + \frac{p\log\log(n)}{\mathfrak{s}^2} \vee \frac{\sqrt{p\log\log(n)}}{\mathfrak{s}} \r) \le n^{\nu}
\end{equation} 
for some $\nu \in (0, 1/2)$, then it holds, for every $\tau \in (0,  1/6]$,
$$
\inf_{\psi}\biggl\{  \sup_{\p \in \mc P_{\mbf I, \sigma}^{0, n, p}} \E_{\p}\l[\psi\l((Y_t, \mbf x_t)_{t\in[n]}\r)\r]  +  \sup_{\p \in \mc P_{\mbf I, \sigma}^{\mathfrak{s}, n, p}(\tau)} \E_{\p}\l[1-\psi\l((Y_t, \mbf x_t)_{t\in[n]}\r)\r]\biggr\} \ge \frac{1}{4e},
$$   
where the infimum is taken over all measurable functions $\psi:\R^{n \times (1+p)}\to[0,1]$. 
%\item\label{lem:lb:id:locate}
%\emph{Estimation:}
%Let $\tau$ be either a large enough constant (in particular, $\tau > \underline{c}_1$), or $\tau \to \infty$ (at an arbitrary rate) as $n\to \infty$. If \(\mathfrak{s}\log\bigl(1 + \tfrac{p\log\log(n)}{\mathfrak{s}^2}\vee\tfrac{\sqrt{p\log\log(n)}}{\mathfrak{s}}\bigr) \le \bigl(\tfrac{n}{\tau}\bigr)^\nu\) for some $\nu \in (0, 1/2)$, then there exists a constant $\underline{c}_2 > 0$ such that
%\[
%\inf_{\wh\theta} \sup_{\p \in \mc P_{\mbf I, \sigma}^{\mathfrak{s}, n, p}(\tau)}\p\l(\bigl\vert\wh\theta\l((Y_t, \mbf x_t)_{t\in[n]}\r) -\theta\bigr\vert \ge \frac{\underline{c}_2 \Psi^2\mathfrak{s}\log\Bigl(1 + \tfrac{p\log\log(n)}{\mathfrak{s}^2}\vee \tfrac{\sqrt{p\log\log(n)}}{\mathfrak{s}}\Bigr)}{\vert\bm\delta\vert_2^2} \r)\ge\frac{1}{4},
%\]
%where the infimum is taken over all measurable functions $\wh\theta:\R^{n \times (p + 1)}\to [n]$. 
%\end{enumerate}
\end{thm}

Theorem~\ref{th:gaussid} provides the minimax lower bound on the detection boundary for the change point problem in~\eqref{eq:amoc}, which reveals a phase transition with respect to the sparsity.  In an isotropic Gaussian setting with $\bm\Sigma = \mbf I$, the sparsity is determined by $\mathfrak{s} = \vert \bm\Sigma^{1/2} \bm\delta \vert_0 = \vert \bm\delta \vert_0 = \mathfrak{s}_\delta$, and the minimax lower bound undergoes a phase transition at $\mathfrak{s} = \sqrt{p{\log\log(n)}}$, from that 
\[
\mathfrak{s} \log\l( 1 + \tfrac{p{\log\log(n)}}{\mathfrak{s}^2} \vee \tfrac{\sqrt {p{\log\log(n)}}}{\mathfrak{s}} \r) \asymp 
\begin{cases}
\mathfrak{s} \log(ep{\log\log(n)}/\mathfrak{s}^2) &\text{ if }{\mathfrak{s} < \sqrt {p{\log\log(n)}}}, \\ \sqrt{p{\log\log(n)}} & \text{ if } {\mathfrak{s} \ge \sqrt {p{\log\log(n)}}}.
\end{cases}
\] 
This phenomenon has hitherto been unknown in the literature, with the existing investigation into minimax lower bounds focusing on the sparse regime of $\mathfrak{s} < \sqrt{p}$ \citep{rinaldo2020localizing, cho2024detection}. 
%Notably, one striking observation in the dense regime is that the minimax rate of estimation becomes a trivial one when $p \gtrsim n^2$. This stands in sharp contrast to the sparse regime, where consistent change point detection remains feasible even in ultra-high dimensions.

\begin{figure}[h]
\centering
\includegraphics[width=\linewidth]{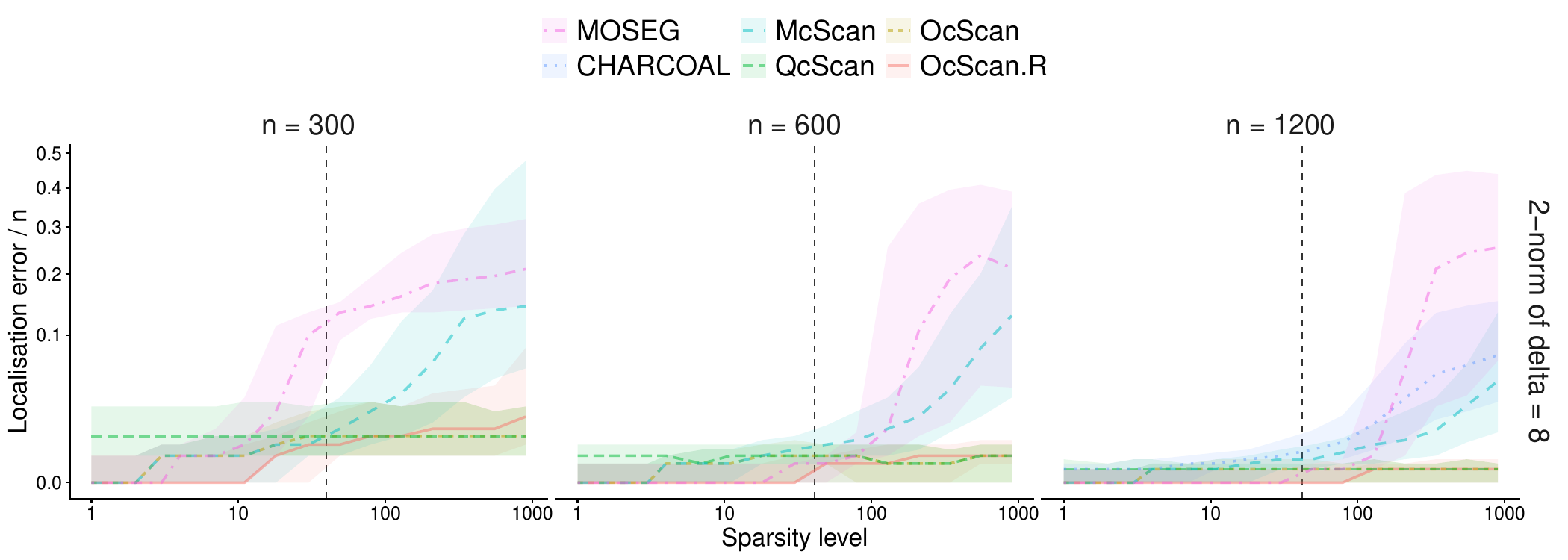}
\caption{Estimation performance of MOSEG \citep{cho2022high}, CHARCOAL (\citealp{gao2022sparse}; only applicable when $n > p$), and the proposed methods:  McScan (Section~\ref{sec:sparse}; a variant of \citealp{cho2024detection}), QcScan (Section~\ref{sec:dense}), OcScan (Section~\ref{sec:adaptive}) and OcScan.R (Section~\ref{sec:refine}), in the isotropic setting ($\bm\Sigma = \mbf I$) for varying $\mathfrak{s}$ ($x$-axis) and sample size $n$ (left to right) while $p = 900$ and $\vert \bm\delta \vert_2 = 8$. 
Results are based on 1000 repetitions, with median error curves shown alongside shaded regions representing the interquartile range. % (upper and lower quartiles). 
The vertical dashed lines correspond to where $\mathfrak{s} = \sqrt{p{\log\log(n)}}$. The $x$-axis is displayed on a log scale, and the $y$-axis on a squared root scale. This is excerpted from Figure~\ref{f:iderr} in \Cref{sec:sim}. \label{f:iderrex}}
\end{figure}

Figure~\ref{f:iderrex} numerically illustrates the phase transition as the sparsity ($x$-axis, in a log scale) and the sample size vary (left to right) % while $\vert \bm\delta \vert_2 = 8$, $p = 900$ and other model parameters remains constant and 
in the isotropic setting.
The figure reports the scaled estimation errors ($y$-axis, in a square root scale) attained by various change point detection methods, including McScan (Section~\ref{sec:sparse}; a variant of \citealp{cho2024detection}) and QcScan (Section~\ref{sec:dense}) which are later shown to be minimax (near-)optimal in the sparse and the dense regimes, respectively. 
While QcScan's estimation error remains roughly constant regardless $\mathfrak{s}$, McScan attains a smaller error for the sparse situations which increases with $\mathfrak{s}$ until it crosses that of QcScan near $\mathfrak{s} = \sqrt{p{\log\log(n)}}$; similar behaviour is also observed from the competing methods (MOSEG and CHARCOAL) tailored for the sparse regime. 
The OcScan (Section~\ref{sec:adaptive}) which combines the outputs from McScan and QcScan, and its refinement OcScan.R (Section~\ref{sec:sim}), are able to adaptively enjoy the optimal performance of McScan and QcScan regardless of the regimes determined by $\mathfrak{s}$.
In particular, OcScan.R improves upon the estimators of McScan and OcScan in the sparse regime, which demonstrates the efficacy of the refinement proposed in Section~\ref{sec:refine}. 

\subsection{General covariance}
\label{sec:general}

For the general case of the model~\eqref{eq:amoc} with $\Cov(\mbf x_t) = \bm\Sigma$ for a non-negative definite matrix $\bm\Sigma \in \R^{p \times p}$, it holds (in the mean squared sense if $\text{rank}(\bm\Sigma) < p$) that 
\begin{align}
\label{eq:mp}
\E(Y_t \mid \mbf x_t) = \l[ (\bm\Sigma^{1/2})^\dagger \mbf x_t \r]^\top \bm\Sigma^{1/2} \l( \bm\beta_0 \cdot \mathbb{I}_{\{t \le \cp\}} + \bm\beta_1 \cdot \mathbb{I}_{\{t > \cp\}} \r),
\end{align}
where $(\bm\Sigma^{1/2})^\dagger$ denotes the (Moore--Penrose) pseudoinverse of $\bm\Sigma^{1/2}$.
Hence the lower bound reported in Table~\ref{tab:minimax} is straightforward generalisation of Theorem~\ref{th:gaussid} to the model class $\mc P^{\mathfrak{s}, n, p}_{\bm\Sigma, \sigma}(\tau)$ with general $\bm\Sigma$.
It is readily seen that the phase transition is determined by the sparsity of the covariance weighted differential parameter, i.e.\ the inherent sparsity $\mathfrak{s} = \vert \bm\Sigma^{1/2} \bm\delta \vert_0$. Also, the lower bound is formulated with the size of the change measured by $\vert \bm\Sigma^{1/2} \bm\delta \vert_2$, complementing the following lemma:

\begin{lem}[Lemma~1 of \citeauthor{cho2024detection}, \citeyear{cho2024detection}]
\label[lem]{lem:tv}
Let $\p_n(\bm\beta)$ denote the joint distribution of $\{(Y_t, \mbf x_t): \,  t \in [n] \}$ such that $Y_t = \mbf x_t^\top\bm\beta + \vep_t$ for all $t$. 
We also denote by $\p_{\cp}(\bm\beta_0, \bm\beta_1)$ the joint distribution of $\{(Y_t, \mbf x_t): \, t \in [n] \}$ under~\eqref{eq:amoc}. 
Then the total variation distance between $\p_n(\bm\beta)$ and $\p_{\cp}(\bm\beta_0, \bm\beta_1)$, denoted by $\mathrm{TV}\bigl( \p_n(\bm\beta),\, \p_{\cp}(\bm\beta_0, \bm\beta_1) \bigr)$, satisfies
\begin{align*}
\frac{1}{100} \le \frac{\min_{\bm\beta, \bm\beta_0, \bm\beta_1: \, \bm\beta_1 - \bm\beta_0 = \bm\delta} \mathrm{TV} \bigl( \p_n(\bm\beta),\, \p_{\cp}(\bm\beta_0, \bm\beta_1) \bigr)}{\min\left\{ 1, \, \sqrt{\frac{\cp (n - \cp)}{n \sigma^2}\bm\delta^\top \bm\Sigma \bm\delta}\right\}} \le \frac{3\sqrt{3}}{2}.     
\end{align*}
\end{lem}

%\todo[inline]{Maybe we should redesign the simulation setup in \Cref{f:toep_0p6_p200_rho2}; We need an example where $\vert\bm\Sigma\bm\delta\vert_{\infty}$ scales with the inherent sparsity $\mathfrak{s}$, rather than $\mathfrak{s}_{\bm\delta}$.}

% When $\bm\Sigma$ is known, McScan and QcScan (see Section~\ref{sec:upper}) are shown to be minimax optimal (up to logarithmic factors) in the sparse and the dense regimes, respectively, via pre-rotation of $\mbf x_t$.
% However, it is imperative to note that McScan and QcScan do not require any knowledge on $\bm\Sigma$ or $\mathfrak{s}$ and, provided that the smallest non-zero eigenvalue $\sigma_*$ of $\bm\Sigma$ is bounded away from zero, they remain minimax rate (near-)optimal even when $\bm\Sigma$ is unknown and rank deficient, see the discussions below Theorems~\ref{thm:mcscan} and~\ref{thm:qcscan}.

\begin{figure}[h]
\centering
\includegraphics[width = \linewidth]{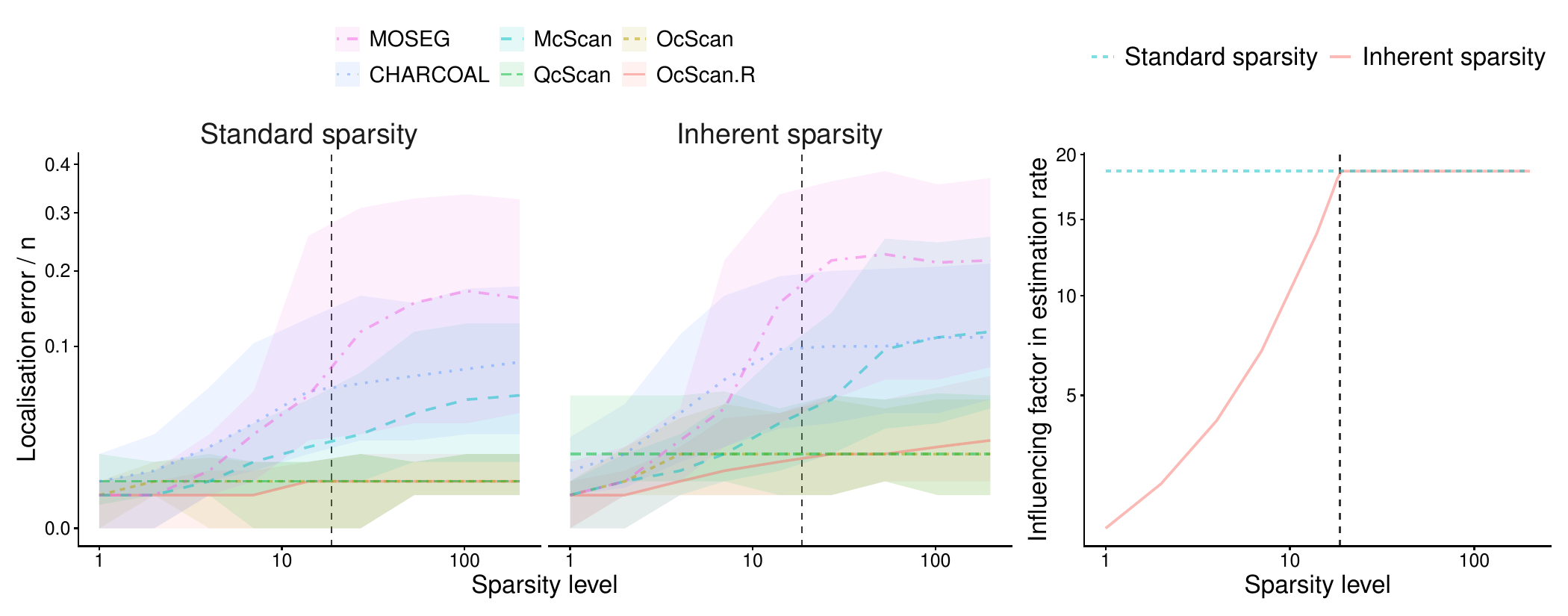}
\caption{Estimation performance of MOSEG \citep{cho2022high}, CHARCOAL \citep{gao2022sparse}, and the proposed methods:  McScan (Section~\ref{sec:sparse}; a variant of \citealp{cho2024detection}),  QcScan (Section~\ref{sec:dense}), OcScan (Section~\ref{sec:adaptive}) and OcScan.R (Section~\ref{sec:refine}), when $\bm\Sigma = [0.6^{\vert i - j \vert}]_{i, j}$. Taken from (M2) in \Cref{sec:sim}, we set $n = 300$, $p = 200$, $\bm\beta_0 + \bm\beta_1 = \mbf 0$ and $\vert\bm\Sigma^{1/2}\bm\delta\vert_2 = 2$. 
In the left, the $x$-axis denotes the standard sparsity $\mathfrak{s}_\delta = \vert \bm\delta \vert_0$ while in the middle, it is the inherent sparsity $\mathfrak{s} = \vert \bm\Sigma^{1/2}\bm\delta \vert_0$. 
In each scenario, results are based on 1000 repetitions, with median error curves shown alongside shaded regions representing the interquartile range. % (upper and lower quartiles). 
The right panel plots the quantity $\min( \vert\bm\Sigma^{1/2}\bm\delta\vert_0, \sqrt{p{\log\log (n)}} )$ as either $\mathfrak{s}_\delta$ or $\mathfrak{s}$ varies. This quantity is the influencing factor, save logarithmic factors, appearing in the minimax rate of estimation (see \Cref{tab:minimax}), and represents a proxy of the difficulty of the estimation problem. 
The vertical dashed lines mark where $\mathfrak{s}_\delta = \sqrt{p {\log\log (n)}}$ or $\mathfrak{s} = \sqrt{p{\log\log (n)}}$. 
The $x$-axis is shown on a log scale, and the $y$-axis on a squared root scale. 
% For further comparisons, see \Cref{f:toep_0p6_p200_rho2,f:toep_0p6_p200,f:toep_0} and \Cref{f:toep_0p6_p400,f:toep_n0p6_p200,f:toep_n0p6_p400} in the appendix.
\label{f:toep_0p6_p200_rho2}}
\end{figure}

Figure~\ref{f:toep_0p6_p200_rho2} numerically confirms these observations when $\bm\Sigma = [0.6^{\vert i - j \vert}]_{i, j}$. Specifically, the best performance of various change point estimation methods exhibits a steeper dependence on the intrinsic sparsity (middle panel) than on the standard sparsity (left panel), which is in alignment with the theoretical results reported in Table~\ref{tab:minimax} (see also the right panel). Besides, we notice that McScan, a method attaining minimax rate-optimality in the sparse regime, has its estimation error exceed that of QcScan, a method tailored to achieve optimal performance in the dense regime, close to where the inherent sparsity $\mathfrak{s}$ meets $\sqrt{p{\log\log (n)}}$ (middle panel). On the other hand, this phenomenon occurs far before the sparsity of $\bm\delta$ reaches $\mathfrak{s}_\delta = \sqrt{p{\log\log (n)}}$, which is accounted for by that $\vert \bm\Sigma^{1/2} \bm\delta \vert_0 = p$ in all the scenarios considered in the left panel.
% Another noteworthy point is the better estimation performance of OcScan.R, which combines the refinement of McScan with QcScan, over that of OcScan (combining McScan and QcScan without refinement) when $\bm\Sigma^{1/2} \bm\delta$ is very sparse, demonstrating the efficacy of the refinement proposed in Section~\ref{sec:refine}. 

% \todo[inline]{Housen checks \Cref{f:toep_0p6_p200_rho2}; re-define the influence factor with $\sqrt p$? There are two options:\\
% 1. We stay with the quantity $\min( \vert\bm\Sigma^{1/2}\bm\delta\vert_0, \sqrt{p{\log\log (n)}})$, and the saying \enquote{save logarithmic factors} is fine (as this is the case for sparse regime). The advantage is that we have consistent location for the vertical dashed lines.  \\ 
% 2. The figure \texttt{Toeplitz\_Extract\_10Gamma\_6\_p\_200\_rho\_2\_nu\_0\_new.pdf} has the same left two panels, but with updated rightmost panel with quantity $\min( \vert\bm\Sigma^{1/2}\bm\delta\vert_0, \sqrt{p})$, together with vertical line at $\sqrt{p}$. This seems to be a bit inconsistent in vertical lines.\\
% I would prefer to stay with option 1 (the current version), but open to your opinions.}

\begin{rem} % [Phase transition in change point problems]
In the high-dimensional change point detection literature, the phase transition with respect to sparsity has been noted in the context of mean change detection, see \cite{enikeeva2019high}, \cite{liu2021minimax}, \cite{pilliat2023optimal} and \cite{kovacs2024optimistic}, all considering the case of cross-sectional independence.
Theorem~\ref{th:gaussid} is a first of its kind for the linear model setting where the cross-sectional dependence is present in $\mbf x_t Y_t$ (whose expectation is piecewise constant as $\E(\mbf x_t Y_t) = \bm\Sigma \bm\beta_0 + \bm\Sigma \bm\delta \cdot \mathbb{I}_{\{t \ge \cp + 1\}}$), even if $\Cov(\mbf x_t) = \mbf I$.
Moreover, the sparsity inherent to this change point problem is characterised by $\vert \bm\Sigma^{1/2} \bm\delta \vert_0$ rather than the $\ell_0$-norm of $\bm\delta$, which also sets our findings apart from the existing work.
\end{rem}

\section{Upper bounds}
\label{sec:upper}

% This section proposes change point method that attains minimax optimality in both sparse and dense regimes, matching the lower bounds derived in Theorem~\ref{th:gaussid} (up to a logarithmic term).
% Section~\ref{sec:sparse} reviews the McScan method proposed by \cite{cho2024detection} and demonstrates its optimality in the sparse regime. 
% Furthermore, focusing on the at-most-one-change setting (AMOC, $q \le 1$), we propose a new method that achieves near-optimality in the dense regime in Section~\ref{sec:dense}.
% Then in Section~\ref{sec:adaptive}, we propose a strategy to combine the two methods that is adaptively optimal in both regimes.

\subsection{Sparse regime}
\label{sec:sparse}

\citet{cho2024detection} propose McScan, a data segmentation method for the multiple change point problem in linear models.
Motivated by the observation that $\E(\mbf x_tY_t)$ contains the information about local parameter vector, it searches for change points by scanning the local covariance between $Y_t$ and $\mbf x_t$.
Avoiding the costly computation of local parameter estimates via $\ell_1$-regularision, it is highly efficient computationally and moreover, attains better statistical efficiency over a wider class of change point problems when compared to the state-of-the-art, by relaxing the conditions of exact sparsity imposed on the regression or differential parameters. 

In the AMOC setting of~\eqref{eq:amoc}, we consider the detector statistic 
\begin{align}
\bar{T}_k = \vert \mbf X^\top \wt{\mbf Y}(k) \vert_\infty &= \sqrt{\frac{k (n - k)}{n}} \l\vert \frac{1}{n - k} \sum_{t = k + 1}^n \mbf x_tY_t - \frac{1}{k} \sum_{t = 1}^k \mbf x_tY_t \r\vert_\infty, 
\label{eq:mcscan}
\end{align}
where $\mbf X = [\mbf x_1, \ldots, \mbf x_n]^\top$ and $\wt{\mbf Y}(k) = \bigl(\wt{Y}_1(k), \ldots, \wt{Y}_n(k)\bigr)^\top$ with
\begin{align}
\wt{Y}_t(k) &= - \sqrt{\frac{n - k}{nk}} Y_t \cdot \mathbb{I}_{\{t \le k\}} + \sqrt{\frac{k}{n (n - k)}} Y_t \cdot \mathbb{I}_{\{t > k\}}.
\label{eq:tilde}
\end{align}
We propose to evaluate $\bar{T}_k$ at strategically selected locations for the detection and estimation of the change point (if any). 
For this, we adopt the advanced optimistic search proposed by \cite{kovacs2024optimistic}, see Algorithm~\ref{alg:os} for the generic pseudo code applicable with any detector statistic denoted by $V_k$.
It first evaluates $V_k$ on a dyadic grid denoted by $\mc D$, for a given trimming parameter $\varpi$, and identifies the maximiser $k^*$.
If $V_{k^*}$ does not exceed some threshold $\zeta$, we quit the algorithm while if $V_{k^*} > \zeta$, we perform a divide-and-conquer-type search through calling the function \texttt{OS} over the interval determined with $k^*$ and its neighbours on the dyadic grid.
Through this, we iteratively compare $V_k$ at dynamically selected search locations based on prior evaluations and discard the sections of the data that are not likely to contain any change point.
Replacing the conventional exhaustive search, this adaptive strategy significantly reduces the total computational complexity from linear to logarithmic, while maintaining the statistical efficiency.
% The idea behind the algorithm is to replace the conventional exhaustive search for change points with an adaptive search strategy which, following a divide-and-conquer principle, dynamically selects subsequent search locations based on prior evaluations. In doing so, the total computational complexity is significantly reduced from linear to logarithmic, while the statistical efficiency is maintained.
We apply Algorithm~\ref{alg:os} with McScan-specific detector statistic $\bar{T}_k$, threshold $\zeta_{\Mc}$ and trimming parameter $\varpi_{\Mc}$, to obtain the estimators of the number and location of the change point, $(\wh\cp_{\Mc}, \wh q_{\Mc})$.  With a slight abuse of terminology, we refer to the resulting method also as McScan (\underline{m}aximum \underline{c}ovariance \underline{scan}ning).

\begin{algorithm}[h!t!b!]
\caption{Advanced optimistic search (adapted from \citealp{kovacs2024optimistic}).}
\label{alg:os}
\DontPrintSemicolon
% \SetAlgoLined
\SetKwProg{Fn}{Function}{}{}
\SetKwFunction{OS}{OS}

\KwIn{Data $(Y_t, \mbf x_t)_{t \in [n]}$, detector statistic $V_k$, threshold $\zeta > 0$, trimming $\varpi \ge 0$}
\BlankLine

\Fn{\OS{$s$, $t$, $e$} with integers $0\le s < t <e \le n$}{ 
    \uIf{$e - s \le 5$}{ 
        $\wh\cp \leftarrow \mathop{\arg\max}_{k \in \{s + 1, \ldots, e -1\}} V_k$
        
        \Return $\wh\cp$
    }
    \BlankLine
    
    \uIf{$e - t > t - s$}{ 
        $w \leftarrow \lceil e - (e - t)/2 \rceil$
        
        \lIf{$V_w \ge  V_t$}{
            $\OS(t, w, e)$
        } \lElse{
            $\OS(s, t, w)$
        }
    }
    \BlankLine

    \uIf{$e - t \le t - s$}{ 
        $w \leftarrow \lfloor s + (t - s)/2 \rfloor$
        
        \lIf{$ V_w \ge  V_t$}{
            $\OS(s, w, t)$
        } \lElse{
            $\OS(w, t, e)$
        }
    }   
}
\BlankLine

$\ell \leftarrow \lfloor \log_2(n/(2\varpi)) \rfloor$
and $\mc D \leftarrow \{ \lfloor 2^{-l}n \rfloor, \lceil n - 2^{-l} n \rceil \, : \, l \in [\ell] \}$ 

$k^* \leftarrow \mathop{\arg\max}_{k \in \mc D} V_{k}$
\BlankLine

\uIf{$V_{k^*} > \zeta$}{

\lIf{$k^* \le n/2$}{
$a \leftarrow \lfloor k^*/2 \rfloor$ and $b \leftarrow \lceil 2 k^* \rceil$
} \lElse{
$a \leftarrow \lfloor k^* - (n - k^*) \rfloor$ and $b \leftarrow \lceil k^* + (n - k^*)/2 \rceil$
}

$\wh\cp \leftarrow \OS(a, k^*, b)$ and $\wh q \leftarrow 1$
}
\lElse{$\wh\cp \leftarrow n$ and $\wh q \leftarrow 0$}
\BlankLine

\KwOut{$(\wh\cp, \, \wh q)$}

\end{algorithm}

\begin{thm}% [Theorem~2 of \citealp{cho2024detection}] 
\label[thm]{thm:mcscan}
Assume the AMOC model~\eqref{eq:amoc}, and set $\varpi_{\Mc} \asymp \log(p \log(n))$. Suppose either that there is a change point, with $q = 1$, $\Delta \ge 2 \varpi_{\Mc}$ and
\begin{align}
\vert \bm\Sigma \bm\delta \vert_\infty^2 \Delta \ge \bar{c}_0 \sigma_X^2 \Psi^2 \log(p \log(n)),
\label{eq:sparse:dlb}
\end{align}
or that there is no change point, with $q = 0$ and $\cp = n$. Then, McScan applied with $\zeta_{\Mc} = \bar{c} \sigma_X \Psi \sqrt{\log(p \log(n))}$ and $\varpi_{\Mc}$, outputs $(\wh\cp_{\Mc}, \wh q_{\Mc})$ which satisfy
\begin{align*}
\p\l\{ \wh q_{\Mc} = q \, ; \, \vert \wh\cp_{\Mc} - \cp \vert \le \bar{c}_1 \vert \bm\Sigma\bm\delta \vert_\infty^{-2} \sigma_X^2 \Psi^2 \log(p \log(n)) \r\} \ge 1 - \bar{c}_2 (p \log(n))^{-1}
\end{align*}
for some universal constants $\bar{c}$ and $\bar{c}_i \in (0, \infty), \, 0 \le i \le 2$. 
\end{thm}

With $\sigma_*$ denoting the smallest non-zero eigenvalue of $\bm\Sigma$, we have by H\"{o}lder's inequality, 
\begin{align*}
\vert \bm\Sigma \bm\delta \vert_\infty \ge
\frac{ \bm\delta^\top \bm\Sigma^{3/2} \bm\delta }{ \vert \bm\Sigma^{1/2} \bm\delta \vert_1} \ge
\frac{ \sigma_*^{1/2} \vert \bm\Sigma^{1/2} \bm\delta \vert_2 } {\mathfrak{s}^{1/2}}.
\end{align*}
Thus, as long as $\sigma_*$ is bounded away from zero, McScan is minimax rate optimal up to a possible logarithmic factor in the sparse regime $\mathfrak{s} < \sqrt{p\log\log(n)}$ (cf.\ Table~\ref{tab:minimax}), without requiring any prior knowledge of $\bm\Sigma$, or ruling out strong collinearity among the variables that results in rank-deficient $\bm\Sigma$ --- a situation that often arises with increasing number of variables.
This is in stark contrast to the most existing methods, which require all eigenvalues of $\bm\Sigma$ to be uniformly bounded above and away from zero; see Section~S1.1 of \cite{cho2024detection} and also Figure~\ref{f:rank} in Section~\ref{sec:sim} for numerical validation. 
When $\bm\Sigma$ is known, pre-rotating $\mbf x_t$ as in~\eqref{eq:mp} allows for re-formulating the condition in~\eqref{eq:sparse:dlb}, to
\begin{align*}
\vert \bm\Sigma^{1/2} \bm\delta \vert_\infty^2 \Delta \gtrsim \Psi^2 \log(p \log(n)).
\end{align*}
Since $\vert \bm\Sigma^{1/2} \bm\delta \vert_\infty \ge \mathfrak{s}^{-1/2} \vert \bm\Sigma^{1/2} \bm\delta \vert_2$ (recall that $\mathfrak{s} = \vert \bm\Sigma^{1/2} \bm\delta \vert_0$), the above displayed condition is independent of $\sigma_*$, and matches the minimax lower bound on the detection boundary in the sparse regime provided that $\log(n) \lesssim p$ and  $\mathfrak{s}$ is not very close to $\sqrt{p\log\log(n)}$.

\begin{rem} % [McScan for data segmentation]
\label{rem:mcscan}
The original McScan algorithm of \cite{cho2024detection} combines the detector statistic in~\eqref{eq:mcscan} with the narrowest-over-threshold principle \citep{baranowski2019narrowest} and the seeded intervals \citep{kovacs2020seeded}, evaluating the detector statistics over strategically selected intervals and selecting the significant local maximisers for the detection of multiple change points. 
Theorem~\ref{thm:mcscan} is a variant of the more general Theorem~2 in \cite{cho2024detection} which, permitting serial dependence and sub-Weibull tail behaviour in the data, shows the consistency of McScan in multiple change point detection, all without placing any condition on the exact sparsity of regression or differential parameters, or that on $p$ beyond what is implied by the detection lower bound comparable to~\eqref{eq:sparse:dlb}.
There, it is further demonstrated that McScan allows for a strictly wider class of change point models in comparison with the state-of-the-art. %, see Section~S1 in the supplementary material of \cite{cho2024detection}.
\end{rem}

\subsection{Dense regime}
\label{sec:dense}

Recalling the definition of $\wt{Y}_t(k)$ in~\eqref{eq:tilde}, we consider the detector statistic
\begin{align}
T_k :=& \, \wt{\mbf Y}(k)^\top \l( \mbf X\mbf X^\top - \frac{1}{n} \tr(\mbf X\mbf X^\top) \mbf I_n \r) \wt{\mbf Y}(k) = \vert \mbf X^\top \wt{\mbf Y}(k) \vert_2^2 - \frac{1}{n} \Vert \mbf X \Vert_F^2 \vert \wt{\mbf Y}(k) \vert_2^2.
\label{eq:qcscan}
\end{align}
The construction of $T_k$ is motivated by the estimator of the explained variance proposed by \cite{Dicker14} and \cite{VeGa18} and 
% However, $T_k$ infers the $\ell_2$-aggregation of $\bm\Sigma \bm\delta$, the covariance weighted differential parameter
indeed, $T_\cp$ approximates $\tfrac{\cp(n - \cp)}{n} \vert \bm\Sigma\bm\delta \vert_2^2$, which follows from the following general result.
% \begin{align*}
% f_\cp = \frac{\cp(n - \cp)}{n} \vert \bm\Sigma\bm\delta \vert_2^2 + \frac{(n - \cp) \vert \bm\Sigma\bm\beta_0 \vert_2^2 + \cp \vert \bm\Sigma\bm\beta_1 \vert_2^2}{n},
% \end{align*}
% see Proposition~\ref{prop:dense} below.
\begin{prop}
\label[prop]{prop:dense}
Under the AMOC model~\eqref{eq:amoc}, suppose that $p \ge n$ and that there exists some $\eps_n \to \infty$ satisfying $n^{-1} \eps_n \to 0$, such that $\cp \in (\eps_n, n - \eps_n)$ if $\bm\beta_0 \neq \bm\beta_1$, and $\cp = n$ if $\bm\beta_0 = \bm\beta_1$.
Then, with 
\begin{align}
f_k & := \frac{k (n - k)}{n} \l[ 
\l( \frac{(n - \cp)^2}{(n - k)^2} \mathbb{I}_{\{k \le \cp\}} + \frac{\cp^2}{k^2} \mathbb{I}_{\{k > \cp\}} \r) \bm\delta^\top \bm\Sigma^2 \bm\delta
\r.
\nn \\
& \l. - \l( \frac{\cp - k}{(n - k)^2} \mathbb{I}_{\{k \le \cp\}} - \frac{k - \cp}{k^2} \mathbb{I}_{\{k > \cp\}} \r) \bm\delta^\top \bm\Sigma^2 (\bm\beta_0 + \bm\beta_1) 
+ \frac{1}{k} \bm\beta_0^\top \bm\Sigma^2 \bm\beta_0 + \frac{1}{n - k} \bm\beta_1^\top \bm\Sigma^2 \bm\beta_1
\r],
\label{eq:f}
\end{align}
there exist some constants $C_0, C_1 \in (0, \infty)$ such that % (\vert \bm\Sigma^{1/2} \bm\beta_0 \vert_2^2 + \vert \bm\Sigma^{1/2} \bm\beta_1 \vert_2^2 + \sigma^2)
\begin{align*}
& \p\l\{ \l\vert T_k -  f_k \r\vert \ge C_0 \Vert \bm\Sigma \Vert \Psi^2 \sqrt{pz} \r\} \le C_1 e^{-z}
\end{align*}
for any $z \in (0, \eps_n^{1/3})$ and $k \in (\eps_n, n - \eps_n)$.
\end{prop}

In contrast to $\bar{T}_k$ in~\eqref{eq:mcscan}, the statistic $T_k$ involves the $\ell_2$-aggregation of the contrast between the local sample covariance before and after~$k$, and thus lends itself naturally to the detection of the change point in the dense regime. 
We can write $T_k = \sum_{i \in [n]} (\gamma_i - \bar{\gamma}) \bigl( \mbf w_i^\top \wt{\mbf Y}(k) \bigr)^2$, where $(\gamma_i, \mbf w_i), \, i \in [n]$, denote the pairs of eigenvalues and eigenvectors of $\mbf X\mbf X^\top$ and $\bar{\gamma} = n^{-1} \sum_{i \in [n]} \gamma_i$. 
This representation  reveals the importance of {centring} of $\mbf X\mbf X^\top$ in~\eqref{eq:qcscan}, which enables $T_k$ to amplify when $\E(\wt{\mbf Y}(k) \mid \mbf X)$ is more aligned with $\mbf w_i$'s associated with larger $\gamma_i$'s. 
% as $T_k$ makes use of that $\E(\wt{\mbf Y}(k) \vert \mbf X)$ is more aligned with $\mbf w_i$'s associated with larger $\gamma_i$'s.
% An alternative representation of $T_k$ reveals the importance of the `centering' of $\mbf X \mbf X^\top$ in its construction, as 
% \begin{align*}
% T_k = \sum_{i \in [n]} (\sigma_i - \bar{\sigma}) \bigl( \mbf w_i^\top \wt{\mbf Y}(k) \bigr)^2,
% \end{align*}
% where $(\sigma_i, \mbf w_i), \, i \in [n]$, denote the pairs of eigenvalues and eigenvectors of $\mbf X\mbf X^\top$ and $\bar{\sigma} = n^{-1} \sum_{i \in [n]} \sigma_i$.
% In other words, $T_k$ makes use of that $\E(\wt{\mbf Y}(k) \vert \mbf X)$ is more aligned with $\mbf w_i$'s associated with larger $\sigma_i$'s.

We propose to perform change point detection by scanning $T_k$ as the detector statistic in combination with the advanced optimistic search (Algorithm~\ref{alg:os}), which we refer to as QcScan (\underline{q}uadratic \underline{c}ovariance \underline{scan}ning). % the procedure returns the estimators $(\wh\cp_{\Qc}, \wh q_{\Qc})$ that achieve consistency with suitable choices for the threshold $\zeta_{\Qc}$ and the trimming parameter $\varpi_{\Qc}$, as shown below.

\begin{thm}
\label[thm]{thm:qcscan}
Under the AMOC model~\eqref{eq:amoc}, let $p \ge n$ and $\varpi_{\Qc} \asymp \bigl(\log\log(n)\bigr)^3$. Suppose either that there is a change point, with $q = 1$, $\Delta \ge 2 \varpi_{\Qc}$ and
\begin{align}
\vert \bm\Sigma \bm\delta \vert_2^2 \Delta &\ge c_0 \Vert \bm\Sigma \Vert \Psi^2 \sqrt{p \log\log(n)},
\label{eq:dense:dlb} 
\end{align}
or that there is no change point, with $q = 0$ and $\cp = n$. Then, QcScan applied with $\zeta_{\Qc} = c \Vert \bm\Sigma \Vert \Psi^2 \sqrt{p\log\log(n)}$ and $\varpi_{\Qc}$, outputs $(\wh\cp_{\Qc}, \wh q_{\Qc})$ satisfying
\begin{align*}
\p\l\{ \wh q_{\Qc} = q \, ; \, \vert \wh\cp_{\Qc} - \cp \vert \le c_1 \vert \bm\Sigma \bm\delta \vert_2^{-2} \Vert \bm\Sigma \Vert \Psi^2 \sqrt{p \log\log(n)}\r\} \ge 1 -  c_2 \bigl(\log(n)\bigr)^{-1},
\end{align*}
where $c$ and $c_i, \, 0 \le i \le 2$, are constants that depend only on $C_0$ and $C_1$ in Proposition~\ref{prop:dense}.
\end{thm}

With $\sigma_*$ denoting the smallest non-zero eigenvalue of $\bm\Sigma$, we have
$\vert \bm\Sigma \bm\delta \vert_2^2 \ge \sigma_* \vert \bm\Sigma^{1/2} \bm\delta \vert_2^2$.
Thus, QcScan is minimax optimal in the dense regime where $\mathfrak{s} \ge \sqrt{p{\log\log (n)}}$ (cf.\ Table~\ref{tab:minimax}), provided that $\sigma_*^{-1} \Vert \bm\Sigma \Vert$ is bounded from the above (see also the discussions below Theorem~\ref{thm:mcscan}).
When $\bm\Sigma$ is known, pre-rotating $\mbf x_t$ allows us to rewrite~\eqref{eq:dense:dlb} without the dependence on $\bm\Sigma$, as
\begin{align*}
\vert \bm\Sigma^{1/2} \bm\delta \vert_2^2 \Delta \ge c_0 \Psi^2 \sqrt{p\log\log(n)}.
\end{align*}
Note, however, the lower bound (\Cref{th:gaussid}) requires \eqref{e:spn} and hence $p < n$ in the dense regime. See the remark below for a possible matching upper bound in this setting. 

\begin{rem}\label{r:qcscan}
\begin{enumerate}[label = (\alph*)]
\item\label{i:r:np} As in \cite{VeGa18}, we focus on the high-dimensional setting by assuming that $p \ge n$ in Proposition~\ref{prop:dense}.
In the low-dimensional setting with $p < n$, the conclusion of the proposition becomes $\p\{ \vert T_k - f_k \vert \ge C_0 \Vert \bm\Sigma \Vert \Psi^2 \sqrt{nz} \} \le C_1 e^{-z}$ and accordingly, the results in Theorem~\ref{thm:qcscan} no longer depend on $p$. 
%\textcolor{red}{
The condition~\eqref{eq:dense:dlb} implies that $p \lesssim n^2$, which sets the dense regime apart from the sparse regime as in the latter, consistent change point detection is feasible provided that $p = O(\exp(n))$.
%}\footnote{HC: discussion on this point was lost while we were editing the paper -- I think we originally discussed this point with the (incorrect) minimax LB on estimation rate}
%When\footnote{In fact, \eqref{eq:dense:dlb} is feasible only when $\Delta \gtrsim p^2$. 
%\\
%shall we explicitly state it? better in relation to the minimax lower bound?} $p \gtrsim n^2$, the minimax rate of estimation in the dense regime becomes a trivial one, which is distinguished from the sparse regime where consistent change point detection is feasible even in ultra-high dimensions.

\item The proof of Proposition~\ref{prop:dense} is based on the exponential concentration inequalities for Gaussian chaos \citep{adamczak2015concentration}.
We may relax the Gaussianity assumption to permit sub-Weibull distributions by means of Theorems~1.3 and~1.5 of \cite{GSS21}, and thus extend the guarantees of QcScan to non-Gaussian settings. % NOTE: linear serial dependence may be handled following the approach of \cite{chen2022inference} 
\end{enumerate}
\end{rem}

\subsection{Post-detection estimation of the differential parameter}
\label{sec:post}

Once a change point is detected and its location estimated, it is of interest to directly learn the differential parameter $\bm\delta$ and its contribution to the overall variability.
We explore the problem of post-detection estimation of~$\bm\delta$ and related quantities.

\subsubsection{When the differential parameter is sparse}
\label{sec:post:sparse}

Directly imposing sparsity on $\bm\delta$, \cite{cho2024detection} propose $\ell_1$-penalised estimator for its estimation, which is obtained as
\begin{align}
\wh{\bm\delta} \in \mathop{\arg\min}_{\mbf a \in \R^p} \frac{1}{2n}
\l\vert \bmx - \frac{n}{\wh\cp} \mbf Y_{0, \wh\cp} \\ \frac{n}{n - \wh\cp} \mbf Y_{\wh\cp, n} \emx - \mbf X\mbf a \r\vert_2^2 + \lambda \sqrt{\frac{n}{\wh\cp (n - \wh\cp)}} \vert \mbf a \vert_1
\label{eq:lope}
\end{align}
for some $\lambda > 0$ and $\wh\cp$ estimating $\cp$.
Their Proposition~3 translates to our setting as follows:
\begin{prop}
\label[prop]{prop:lope}
Suppose that \textit{either} (i)~$\wh{\bm\delta}$ in~\eqref{eq:lope} is obtained with $\wh\cp = \wh\cp_{\Mc}$ and the conditions in Theorem~\ref{thm:mcscan} are met, \textit{or} (ii) $\wh{\bm\delta}$ is obtained with $\wh\cp = \wh\cp_{\Qc}$ and the conditions in Theorem~\ref{thm:qcscan} are met. 
Let $\underline{\sigma}$ denote the smallest eigenvalue of $\bm\Sigma$. 
Provided that $n \ge 4 C_{\text{RE}} \underline{\sigma}^{-2} \log(p \vee n)$ and $\lambda = C_\lambda \sigma_X \Psi \sqrt{\log(p \vee n)}$ for some constants $C_{\text{RE}}, C_\lambda \in (0, \infty)$, we have under~(i), 
\begin{align*}
\bigl\vert \wh{\bm\delta} - \bm\delta\bigr \vert_2 &\lesssim \frac{\sigma_X \Psi \sqrt{\mathfrak{s}_\delta \log(p \vee n)}}{\underline{\sigma} \sqrt{\Delta}} 
\text{ \ and \ }
\bigl\vert \wh{\bm\delta} - \bm\delta \bigr\vert_1 \lesssim \frac{\sigma_X \Psi \mathfrak{s}_\delta \sqrt{\log(p \vee n)}}{\underline{\sigma} \sqrt{\Delta}},
\end{align*}
with probability at least $1 - \bar{c}_2(p\log(n))^{-1} - c'(p \vee n)^{-1}$, and  under~(ii),
\begin{align*}
\bigl\vert \wh{\bm\delta} - \bm\delta\bigr \vert_2 &\lesssim \frac{ \Psi \max\{ \sigma_X \sqrt{\mathfrak{s}_\delta \log(p \vee n)}, \Vert \bm\Sigma \Vert^{1/2} \bigl(p\log\log(n)\bigr)^{1/4} \}}{\underline{\sigma} \sqrt{\Delta}} 
\text{ \ and \ }
\\
\bigl\vert \wh{\bm\delta} - \bm\delta \bigr\vert_1 &\lesssim \frac{ \Psi \sqrt{\mathfrak{s}_\delta} \max\{ \sigma_X \sqrt{\mathfrak{s}_\delta \log(p \vee n)}, \Vert \bm\Sigma \Vert^{1/2} \bigl(p\log\log(n)\bigr)^{1/4} \}}{\underline{\sigma} \sqrt{\Delta}},
\end{align*}
with probability at least $1 - c_2 \bigl(\log(n)\bigr)^{-1} - c'(p \vee n)^{-1}$, for some universal constant $c' \in (0, \infty)$ and $\bar{c}_2$ and $c_2$ in Theorems~\ref{thm:mcscan} and~\ref{thm:qcscan}.
\end{prop}
For the consistency of $\wh{\bm\delta}$, we make the stronger assumption requiring $\underline{\sigma} > 0$, compared to that required for the consistency in change point detection involving $\sigma_*$, the smallest non-zero eigenvalue of $\bm\Sigma$.
In the inherently sparse regime with $\wh\cp = \wh\cp_{\Mc}$, the consistency of $\wh{\bm\delta}$ is analogous to that of the Lasso estimator of the regression parameter (see e.g.\ \citeauthor{bickel2009simultaneous}, \citeyear{bickel2009simultaneous}), with $\Delta = \min(\cp, n - \cp)$ taking the role of sample size.
\cite{cho2024detection} additionally propose a bias-corrected estimator of $\bm\delta$ and establish a non-asymptotic Gaussian approximation result, which enables simultaneous inference about all $p$ coordinates of $\bm\delta$. 
Later in Section~\ref{sec:refine}, we devise a change point estimator which, based on scanning the projected covariance as $\wh{\bm\delta}^\top \mbf X^\top \wt{\mbf Y}(k)$, attains a refined rate of estimation that is inversely proportional to $\bm\delta^\top \bm\Sigma \bm\delta$. %, improving upon that of~$\wh\cp_{\Mc}$.

\subsubsection{For the differential parameter of arbitrary sparsity}
\label{sec:post:dense}

Generally, it is infeasible to infer $\bm\delta$, or even the regression parameter under stationarity without any assumption on its sparsity.
Yet it is still of interest to learn the strength of the change measured by $V_\circ := \bm\delta^\top \bm\Sigma \bm\delta$, a quantity closely related to the signal-to-noise ratio or the explained variance in linear models. 
In particular, $\sigma^{-2} V_\circ$ determines the total variation distance between the joint distributions of $\{(Y_t, \mbf x_t)\}_{t \in [n]}$ with and without a change, see Lemma~\ref{lem:tv}.
Motivated by Proposition~\ref{prop:dense} and Theorem~\ref{thm:qcscan}, we propose the following estimator 
\begin{align*}
\wh{V} = \frac{n}{\wh\cp_{\Qc} (n - \wh\cp_{\Qc})} T_{\wh\cp_{\Qc}}
\end{align*}
for $V := \vert \bm\Sigma \bm\delta \vert_2^2$, the $\ell_2$-aggregation of the covariance weighted differential parameter, as a proxy of $V_\circ$.
%The trick of taking the median is similarly motivated as the modification to the advanced optimistic search procedure in Algorithm~\ref{alg:os}, namely from that the sub-exponential tail bound in Proposition~\ref{prop:dense} excludes $k$'s too close to $\cp$.

\begin{prop}
\label[prop]{prop:v}
Suppose that the conditions in Theorem~\ref{thm:qcscan} hold. Then, we have
\begin{align}
\label{eq:v}
\frac{\vert \wh{V} - V \vert}{V} \lesssim \frac{ \Vert \bm\Sigma \Vert \Psi^2 \sqrt{p\log\log(n)}}{\vert \bm\Sigma \bm\delta \vert_2^2 \Delta}
% + \frac{\varpi_{\Qc}}{\Delta} + \frac{1}{\sqrt{p\log\log(n)}}
\end{align}
with probability at least $1 - c'' \bigl(\log(n)\bigr)^{-1}$, with $c'' \in (0, \infty)$ depending only on $C_1$ in Proposition~\ref{prop:dense}.
\end{prop}

The upper bound in~\eqref{eq:v} tends to zero as the gap between the two sides of~\eqref{eq:dense:dlb} increases.
If $\bm\Sigma$ is known, we can pre-rotate the data so that $\wh{V}$ estimates $V_\circ$ directly. 
In the absence of such knowledge, we conjecture that it is not possible to estimate $V_\circ$ without imposing any constraint on the sparsity of $\bm\delta$, which is supported by the impossibility result derived in Theorem~4.5 of \cite{VeGa18} for the related problem of estimating the explained variance in stationary linear models.

\section{Adaption to unknown sparsity, refinement and computational complexity}
\label{sec:four}

\subsection{Adaptively optimal change point detection}
\label{sec:adaptive}

Section~\ref{sec:upper} demonstrates that McScan and QcScan are minimax (near-)optimal in sparse and dense regimes, respectively. Each regime is determined by whether $\mathfrak{s} < \sqrt{p{\log\log (n)}}$ or not, where $\mathfrak{s}$ is typically unknown in practice. This calls for a method that adapts to unknown sparsity level. We propose to combine the outputs from McScan and QcScan, which leads to a procedure that detects a change whenever either of the detection requisites in~\eqref{eq:sparse:dlb} and~\eqref{eq:dense:dlb} is met. Further, the resulting change point estimator inherits the localisation property from either $\wh\cp_{\Mc}$ or $\wh\cp_{\Qc}$ that attains the better rate of estimation. In short, referred to as OcScan (\underline{o}ptimal \underline{c}ovariance \underline{scan}ning), the combined method adaptively enjoys the optimality of McScan and QcScan without the explicit knowledge of $\mathfrak{s}$.

It is natural to determine the existence of the change point by $\wh q = \max\{ \mathbb{I}_{\{ \bar{T}_{\wh\cp_{\Mc}} > \zeta_{\Mc} \}}, \, \mathbb{I}_{\{ T_{\wh\cp_{\Qc}} > \zeta_{\Qc} \}} \}$. %, which ensures that a change is detected regardless of the regime determined by $\mathfrak{s}$.
If $\wh q = 0$, we set ${\wh\cp_{\Oc}} = n$, while if $\wh q = 1$ (with the convention that $\bar{T}_n = T_{n} = 0$),
\begin{align}
& \wh\cp_{\Oc} = \begin{cases}
\wh\cp_{\Mc} & \text{if \ } \mathbb{I}_{\{ \bar{T}_{\wh\cp_{\Mc}} > \zeta_{\Mc} \}} \cdot \mathbb{I}_{\{ T_{\wh\cp_{\Qc}} \le \zeta_{\Qc} \}} = 1 \text{ \ or \ } C_{\Mc/\Qc} > 1, \\
\wh\cp_{\Qc} & \text{if \ } \mathbb{I}_{\{ \bar{T}_{\wh\cp_{\Mc}} \le \zeta_{\Mc} \}} \cdot \mathbb{I}_{\{ T_{\wh\cp_{\Qc}} > \zeta_{\Qc} \}} = 1 \text{ \ or \ } C_{\Mc/\Qc} \le 1, \\
\end{cases} 
\label{eq:ocscan}
\\
& \text{with \ } C_{\Mc/\Qc} = \l( \frac{T_{\wh\cp_{\Qc}}}{{ \Vert \bm\Sigma \Vert} \sqrt{p \log\log(n)}} \r)^{-1} \frac{\bar{T}_{\wh\cp_{\Mc}}^2}{\sigma_X^2 \log(p \log(n))},
\nn
\end{align}
from the observation that $\bar{T}^2_{\wh\cp_{\Mc}}$ (resp.\ $T_{\wh\cp_{\Qc}}$) serves as a good proxy for the magnitude of $\vert \bm\Sigma \bm\delta \vert_\infty^2 \Delta$ (resp.\ $\vert \bm\Sigma \bm\delta \vert_2^2 \Delta$), and thus $C_{\Mc/\Qc}$ is indicative of which estimator is to be preferred if $\mathbb{I}_{\{ \bar{T}_{\wh\cp_{\Mc}} > \zeta_{\Mc} \}} \cdot \mathbb{I}_{\{ \bar{T}_{\wh\cp_{\Qc}} > \zeta_{\Qc} \}} = 1$.

\begin{prop}
\label[prop]{prop:adapt} 
With large enough constants $\bar{c}, c \in (0, \infty)$, set $\zeta_{\Mc} = \bar{c} \sigma_X \Psi \sqrt{\log(p \log(n))}$ and $\zeta_{\Qc} = c \Vert \bm\Sigma \Vert \Psi^2 \sqrt{p\log\log(n)}$, and also $\varpi_{\Mc} \asymp \log(p \log(n))$ and $\varpi_{\Qc} \asymp \bigl(\log\log(n)\bigr)^3$.
{Suppose that $p \ge n$, and} in the case of a change point ($q = 1$), we suppose that at least one of~\eqref{eq:sparse:dlb} and~\eqref{eq:dense:dlb} holds.  
If no change point is present ($q = 0$), we regard that $\cp = n$.  
Then, there exist some constants $\bar{c}_i$ and $c_i$, $i = 0, 1, 2$ ($\bar{c}_0$ and $c_0$ appearing in~\eqref{eq:sparse:dlb} and~\eqref{eq:dense:dlb}), which may depend on $\bar{c}$ and $c$, such that, 
\begin{multline*}
\p\l\{ \wh q = q \, ; \, \vert \wh\cp_{\Oc} - \cp \vert \le \Psi^2
\min\l[ \bar{c}_1 \sigma_X^2 \vert \bm\Sigma \bm\delta \vert_\infty^{-2} \log(p\log(n)), \, c_1 \Vert \bm\Sigma \Vert \vert \bm\Sigma \bm\delta \vert_2^{-2} \sqrt{p \log\log(n)} \r] \r\} 
\\
\ge 1 - \bar{c}_2(p \log(n))^{-1} - c_2 \bigl(\log(n)\bigr)^{-1}.
\end{multline*}
\end{prop}

\begin{rem}
Adapting to unknown sparsity in estimation and inference has been actively explored in the high-dimensional statistics literature.
On the change point front, \cite{wang2023computationally} establish the distributions of max- and sum-type statistics for testing for a change in the mean of high-dimensional panel data, as well as their asymptotic independence, which enables sparsity-adaptive testing. 
For this problem, aggregation based on a cross-sectional grid has also been considered for change point detection \citep{enikeeva2019high, liu2021minimax, pilliat2023optimal, kovacs2024optimistic, moen2024efficient}.
For linear models, \cite{VeGa18} propose a sparsity-adaptive estimator of the explained variance which achieves optimality by combining a quadratic form of the response and the square-root Lasso estimator. 
\end{rem}

\begin{rem}\label{r:comp:covar}
The change point problem in linear models can be viewed as a special case of the change point problem in covariances, which allows us to place the performance of our proposed methods in a context. This is particularly fruitful in elucidating the role of quantities such as $\sigma_X$, $\Vert \bm\Sigma \Vert$ and $\vert \bm\Sigma^{1/2} \bm\beta_j\vert_2$ (through $\Psi$) that are often treated as constants in the literature on change point detection in linear models.
Compared to the state-of-the-art methods \citep{wang2021optimal, DePaYa22, moen2024minimax} that target different types of changes in covariances, our OcScan consistently achieves equal or superior performance across a variety of structural settings. For instance, OcScan attains a sharper detection boundary than the projection-based method of \citet{wang2021optimal}, with an improvement factor of $p\min(\Vert \bm\Sigma \Vert \sigma_*^{-1} \mathfrak{s}, \vert \bm\Sigma \bm\delta \vert_0, \sqrt{p})^{-1}$. This gain arises from effectively leveraging the specific structure inherent in the linear model. 
In addition, scanning for (univariate) changes in the marginal variance of the response can be beneficial, particularly in regimes with highly unbalanced signal strengths (namely, when $\lvert \bm\Sigma^{1/2}\bm\beta_0\rvert_2$ and $\lvert \bm\Sigma^{1/2}\bm\beta_1\rvert_2$ differ substantially).
Further discussion and comparison can be found in Appendix~\ref{sec:cov}.
\end{rem}

We finally discuss adaptation to the unknown quantities including $\sigma_X^2$, $\Vert \bm\Sigma \Vert$ and $\Psi$.
Let us write $\wh{\bm\Sigma} = n^{-1} \sum_{t \in [n]} \mbf x_t \mbf x_t^\top$, from which we derive the approximation of $\sigma_X^2$ and $\Vert \bm\Sigma \Vert$ as $\wh\sigma_X^2 = \max_{i \in [p]} \wh\Sigma_{ii}$ and $\Vert \wh{\bm\Sigma} \Vert$, respectively, with $\wh\Sigma_{ii}, \, i \in [p]$, denoting the diagonal entries of $\wh{\bm\Sigma}$.
Also, let $\wh\Psi^2_{1, t} = t^{-1} \sum_{u = 1}^t Y_u^2$ and $\wh\Psi^2_{2, t} = t^{-1} \sum_{u = n - t + 1}^n Y_u^2$, and define a dyadic grid $\mc T = \{2^\ell \, : \, \lceil \log_2\log\log(n) \rceil \le \ell \le \lfloor \log_2(n/2) \rfloor\}$, from which we obtain $\wh\Psi^2 = \max_{t \in \mc T} (\wh\Psi^2_{1, t} \vee \wh\Psi^2_{2, t})$.
With probability tending to one, we can show that $\wh\sigma_X^2 \asymp \sigma_X^2$, $\Vert \wh{\bm\Sigma} \Vert \asymp \Vert \bm\Sigma \Vert$ and $\wh\Psi^2 \asymp \Psi$, leading to:
\begin{cor}
\label[cor]{cor:adapt:short} 
Suppose that $n \le p \le c_*n$ for some constant $c_* \in (0, \infty)$. 
Then, Proposition~\ref{prop:adapt} continues to hold with $\zeta_{\Mc} = \bar{c} \wh\sigma_X \wh\Psi \sqrt{\log(p \log(n))}$ and $\zeta_{\Qc} = c \Vert \wh{\bm\Sigma} \Vert \wh\Psi^2 \sqrt{p\log\log(n)}$.
\end{cor}
The restrictive condition that $p \le c_*n$ arises from controlling $\Vert \wh{\bm\Sigma} - \bm\Sigma \Vert$, and it can be lifted when deriving the consistency of McScan with the adaptively selected $\zeta_{\Mc}$.
Alternatively, we may adopt bootstrap techniques \citep{chernozhukov2023nearly} in combination with McScan; a full exploration of such an approach lies beyond the scope of the present paper but represents an interesting direction for future research. 
%\footnote{potentially better alternative would be to adopt bootstrapping or even resampling under independence. self-normalisation tends to reduce power I think.. { Honestly, I did not know how self-normalisation empirically works. Shall we mention both alternatives? Do you have a proper reference for the bootstrapping or subsampling? } \\ if I recall correctly Xiaofeng's paper \url{https://arxiv.org/pdf/1905.08446} demonstrates that this approach does not achieve the same level of efficiency as the minimax optimal method. Re. bootstrapping we could add \cite{chernozhukov2023nearly}? at least for McScan their result on the multiplier bootstrapping could work, see \url{https://www.sciencedirect.com/science/article/pii/S0047259X25000442} (just searched with keywords). As for re-sampling, it's used e.g. in Matteson and James but they don't discuss how it performs under the alternative. }

\subsection{Refined change point estimation}
\label{sec:refine}

When $\bm\delta$ itself is sparse, further refinement of $\wh\cp_{\Oc}$ output by OcScan as in~\eqref{eq:ocscan} is feasible by utilising the estimator of $\bm\delta$.
Specifically, recall from~\eqref{eq:tilde} the definition of $\wt{Y}_t(k)$, and $\wh{\bm\delta}$ obtained as in~\eqref{eq:lope} with $\wh\cp = \wh\cp_{\Oc}$ which estimates $\bm\delta$.
With some $\varpi_{\Rf} > 0$, we define OcScan.R (a refinement of OcScan)  as
\begin{align*}
\wh{\cp}_{\Rf} = \mathop{\arg\max}_{\varpi_{\Rf} < k < n - \varpi_{\Rf}} \, \wh{\bm\delta}^\top \mbf X^\top \wt{\mbf Y}(k).
\end{align*}

\begin{prop}
\label[prop]{prop:refine}
Suppose that the conditions in Proposition~\ref{prop:adapt} are met with $q = 1$, and set $\varpi_{\Rf} \asymp \log(p \vee n)$. Recall that $\underline{\sigma}$ denotes the smallest eigenvalue of $\bm\Sigma$.
Then, provided that
\begin{align}
\bm\delta^\top \bm\Sigma \bm\delta \cdot \Delta \ge \wt{c}_0  \underline{\sigma}^{-2} \Vert \bm\Sigma \Vert^2 \Psi^2 \mathfrak{s}_\delta^2 \log(p \vee n),
\label{eq:prop:refine}
\end{align}
there exist some universal constants $\wt{c}_i \in (0, \infty), \, i = 0, 1$, such that
\begin{align*}
\p\l\{ \vert \wh{\cp}_{\Rf} - \cp \vert \le \wt{c}_1 (\bm\delta^\top \bm\Sigma \bm\delta)^{-1} \Psi^2 \log(p \vee n) \r\} \ge 1 - \bar{c}_2 (p \log(n))^{-1} - c_2 \bigl(\log(n)\bigr)^{-1} - c' (p \vee n)^{-1},
\end{align*} 
where $\bar{c}_2$ and $c_2$ are as in Proposition~\ref{prop:adapt} and $c'$ as in Proposition~\ref{prop:lope}.
\end{prop}
The condition in~\eqref{eq:prop:refine} requires the exact sparsity of $\bm\delta$ and, if $p \ge n$ (as assumed in Theorem~\ref{thm:qcscan} and Proposition~\ref{prop:adapt}), it indicates that $\mathfrak{s}_\delta \lesssim \sqrt{p}$.
At the price of the stronger condition, the rate attained by $\wh\cp_{\Rf}$ is sharper than that of $\wh\cp_{\Oc}$ up to a multiplicative factor as large as $\sqrt p$ (save logarithmic terms), thereby yielding substantial improvements in high dimensions. This can be seen from the following inequalities
% \footnote{
% We use the inequalities below $+$ note `max'
% \begin{align*}
% \mathfrak{s}_\delta^{-1} \underline{\sigma} \vert \bm\delta \vert_1^2 \le \underline{\sigma} \vert \bm\delta \vert_2^2 \le \bm\delta^\top \bm\Sigma \bm\delta \le \vert \bm\Sigma\bm\delta \vert_\infty \vert \bm\delta \vert_1
% \end{align*}
% {
% For the lower bound, we use $\lvert\bm\Sigma\bm\delta\rvert_\infty \le\lVert\bm\Sigma\rVert_{\infty}\lvert\bm\delta\rvert_{\infty}$ and then
% \[
%  \max_{\bm\delta \neq 0, \vert\bm\delta\vert_0 \le \mathfrak{s}_\delta} \frac{\bm\delta^\top\bm\Sigma\bm\delta}{\vert \bm\Sigma \bm\delta \vert_\infty^2} \ge  \max_{\bm\delta \neq 0, \vert\bm\delta\vert_0 \le \mathfrak{s}_\delta} \frac{\bm\delta^\top\bm\Sigma\bm\delta}{\lVert\bm\Sigma\rVert_{\infty}^2\lvert\bm\delta\rvert_{\infty}^2}  \ge \max_{\bm\delta \neq 0, \vert\bm\delta\vert_0 \le \mathfrak{s}_\delta} \frac{\underline{\sigma}\lvert\bm\delta\rvert_2^2}{\lVert\bm\Sigma\rVert_{\infty}^2\lvert\bm\delta\rvert_{\infty}^2} = \frac{\underline{\sigma}\mathfrak{s}_\delta}{\Vert \bm\Sigma \Vert_\infty^2}.
% \]}
% }
\begin{align*}
\frac{\underline{\sigma}\mathfrak{s}_\delta}{\Vert \bm\Sigma \Vert_\infty^2} \,\le\, \max_{\bm\delta: 0 <  \vert\bm\delta\vert_0 \le \mathfrak{s}_\delta} \frac{\bm\delta^\top\bm\Sigma\bm\delta}{\vert \bm\Sigma \bm\delta \vert_\infty^2} \,\le\, \frac{\mathfrak{s}_\delta}{\underline{\sigma}} \text{ \ and \ }
\frac{\sqrt{p}}{\Vert \bm\Sigma \Vert}\, \le\, \frac{\bm\delta^\top\bm\Sigma\bm\delta \sqrt{p}}{\vert \bm\Sigma \bm\delta \vert_2^2}\, \le\, \frac{\sqrt{p}}{\underline{\sigma}},
\end{align*}
where $\Vert\bm\Sigma\Vert_\infty$ denotes the maximum absolute row sum of $\bm\Sigma$.  
It is straightforward to verify that the localisation rate in Proposition~\ref{prop:refine} is minimax optimal up to a logarithmic factor, see \Cref{ss:opt:loc}. This optimality, however, does not imply that the localisation rates achieved by McScan, QcScan or OcScan (cf.\ Theorems~\ref{thm:mcscan}, \ref{thm:qcscan} and \Cref{prop:adapt}) are sub-optimal, as the condition in~\eqref{eq:prop:refine} is substantially stronger than the minimax detection boundary. Determining the optimal localisation rates under the weakest detectable signal strengths remains an open problem. Notably, this question is unresolved even in the comparatively simpler setting of high-dimensional mean change detection, see \citet[Section~7.2]{pilliat2023optimal} and \citet[Section~7]{kovacs2024optimistic}.

\begin{rem}
% \begin{enumerate}[label = (\alph*)]
% \item 
Results comparable to Proposition~\ref{prop:refine} are reported in Theorem~3 of \cite{xu2024change}.
Their change point estimator attains the rate $O_P(\vert \bm\delta \vert_2^{-2})$, which allows for the derivation of its limiting distribution, % in both the vanishing ($\vert \bm\delta \vert_2 \to 0$) and non-vanishing regimes. 
provided that $\vert \bm\delta \vert_2^2 \Delta \gg \mathfrak{s}_\delta^2 \log^3(p \vee n)$, $\bm\Sigma$ has bounded eigenvalues (i.e.\ $\underline{\sigma}^{-1} \Vert \bm\Sigma \Vert$ is fixed) and $\Psi \lesssim 1$.
Under comparable conditions, from that $\bm\delta^\top \bm\Sigma \bm\delta \ge \underline{\sigma} \vert \bm\delta \vert_2^2$, the requirement in~\eqref{eq:prop:refine} is weaker by the factor of $\log^2(p \vee n)$, and clearly specifies the effects of $\Vert \bm\Sigma \Vert$, $\underline{\sigma}$ and~$\Psi$.
We also remark that Lemma~\ref{lem:tv} shows the inverse of $\bm\delta^\top \bm\Sigma \bm\delta$ to be a better measure of the inherent difficulty of the detection problem in comparison with $\vert \bm\delta \vert_2^{-2}$. 
% Indeed, they support this claim by an example demonstrating that a change is not detectable even when $\vert \bm\delta \vert_2^2$ is large if $\bm\delta^\top \bm\Sigma \bm\delta$ is small. 
% \item TODO; something about relaxation of $\underline{\sigma} > 0$ using \cite{van2018tight}?
% \end{enumerate}
\end{rem}

%{
%May-Be-Doable: Under the setup of Proposition~\ref{prop:refine}, we may establish the limiting distribution of $\wh{\cp}_{\Rf}$ similarly as in \citet{xu2024change}. Specifically:
%\begin{itemize}
%\item Proposition~\ref{prop:refine} implies that $\vert \wh{\cp}_{\Rf} - \cp \vert \cdot (\bm\delta^\top \bm\Sigma \bm\delta)/ \Psi^2 $ is tight.
%\item We may replace $k \mapsto \wh{\bm\delta}^\top \mbf X^\top \wt{\mbf Y}(k)$ by a process $k \mapsto  {\bm\delta}^\top \mbf X^\top \wt{\mbf Y}(k)-{\bm\delta}^\top \mbf X^\top \wt{\mbf Y}(\cp) $. 
%\item Then we can apply the argmax continuous mapping theorem. 
%\end{itemize}
%Question: How does the limit look like? Does it depend on $\bm\beta_0 + \bm\beta_1$?
%\\
%Thanks to the projection via $\bm\delta$, the problem is cast as that of detecting a change in the mean of the univariate sequence $z_t := \bm\delta^\top \mbf x_t Y_t$.
%Note that
%\begin{align*}
%z_t = \sum_{j = 1}^2 \bm\delta^\top \bm\Sigma \bm\beta_{j-1} \cdot \mathbb{I}_{\{\cp_{j - 1} < t \le \cp_j\}} +  \underbrace{\sum_{j = 1}^2 \bm\delta^\top (\mbf x_t \mbf x_t^\top - \bm\Sigma) \bm\beta_{j-1}\cdot \mathbb{I}_{\{\cp_{j - 1} < t \le \cp_j\}} + \bm\delta^\top \mbf x_t \vep_t}_{\text{noise}}.
%\end{align*}
%Looking at the literature on the limiting distribution of the change point estimator, e.g.\ \cite{cho2022bootstrap}, \cite{antoch1995change} and \cite{antoch1999estimators} (the latter two are better since it is for CUSUM estimator), may be relevant. 
%}

\subsection{Runtime and memory}\label{ss:comp}

We analyse the computational complexity of our methods in terms of both runtime and memory. As a pre-processing step, we compute and store the following quantities:
\begin{subequations}
\begin{align}
a_0 & :=\frac{1}{n} \sum_{t = 1}^n\vert\mbf x_t \vert_2^2, \label{e:prea}\\
r_k & := \sum_{t = 1}^k Y_t^2 \text{ \ and \ } {\bm S}_k := \sum_{t = 1}^k \mbf x_t Y_t \quad \text{for}\quad k \in [n], \label{e:preb}
\end{align}
\end{subequations}
which require $O(np)$ runtime and $O(np)$ memory. 
We recommend that practitioners store their data directly in the form of \eqref{e:prea}--\eqref{e:preb} at the data collection stage whenever possible.
Then, the detector statistics in~\eqref{eq:mcscan} and~\eqref{eq:qcscan} can be written as 
\begin{align*}
\bar{T}_k & =  \sqrt{\frac{n}{k(n-k)}}\l\vert\bm S_k - \frac{k}{n}\bm S_n\r\vert_\infty \text{ and} \\
T_k & = \frac{n}{k(n-k)}\l\vert\bm S_k - \frac{k}{n}\bm S_n\r\vert_2^2 - a_0 \l(\frac{n-2k}{k(n-k)}r_k +\frac{k}{n(n-k)}r_n\r).
\end{align*}
Each evaluation of these statistics takes $O(p)$ time, assuming that the quantitates in \eqref{e:prea}--\eqref{e:preb} are available. Since Algorithm~\ref{alg:os} requires $O\bigl(\log (n)\bigr)$ such evaluations, we derive the following computational complexity for our proposed covariance scanning-based approaches.
\begin{prop}
Consider data $\{(Y_t, \mbf x_t) \in \R^{1+p} : \, t \in [n]\}$, e.g.\ from the AMOC model~\eqref{eq:amoc}. 
\begin{enumerate}[label = (\roman*)]
\item Both runtime and memory requirements of McScan, QcScan and OcScan are $O(np)$. 
\item If the data are pre-processed and stored in the form given by \eqref{e:prea}--\eqref{e:preb}, the runtime of McScan, QcScan and OcScan is reduced to $O\bigl(p \log(n) \bigr)$.
\item The refinement method described in Section~\ref{sec:refine} requires $O(np)$ memory and has a runtime of $O\bigl(\mathrm{Lasso}(n, p)\bigr)$, where $\mathrm{Lasso}(n, p)$ denotes the runtime of solving a Lasso problem of size $n \times p$.
\end{enumerate}
\end{prop}
As a consequence, the runtime of McScan, QcScan, and OcScan is approximately equivalent to a single full pass over the data, up to a small constant factor. See Figure~\ref{f:idtm} in Section~\ref{sec:sim} for empirical evidence supporting this observation.

\section{Simulation study}
\label{sec:sim}

We evaluate the performance of the covariance scanning methods, namely {McScan} (\Cref{sec:sparse}), {QcScan} (\Cref{sec:dense}), {OcScan} (\Cref{sec:adaptive}), and its refined variant OcScan.R (\Cref{sec:refine}). For comparison, we consider {MOSEG} \citep{cho2022high} and {CHARCOAL} \citep{gao2022sparse}, the latter of which is applicable only when $p < n$. {The implementations of MOSEG and CHARCOAL are obtained from the GitHub repositories \url{https://github.com/Dom-Owens-UoB/moseg} and \url{https://github.com/gaofengnan/charcoal}, respectively.} We have also explored other variants of covariance scanning approaches, such as using the full grid search in place of the advanced optimistic search (Algorithm~\ref{alg:os}), and different combinations of refinement and adaptation; these alternatives exhibit similar, though often slightly inferior, performance, see \Cref{app:comp} for details. 
Overall, we recommend OcScan and OcScan.R as they satisfactorily combine the good performance of McScan and QcScan in sparse and dense regimes, respectively, and OcScan.R is shown to enhance the estimation accuracy when the signal-to-noise is large and/or $\bm\delta$ is highly sparse.

{
\subsection{Change point estimation performance}
\label{ss:loc:err}
First, we focus on evaluating the accuracy of the estimated change point location, by assuming that exactly one change point is present. 
% To circumvent the model selection effect, we assume that exactly one change point is present.  
Change point detection performance is considered separately in \Cref{ss:num:detect}.}

\begin{figure}[ht]
\centering
\includegraphics[width=\linewidth]{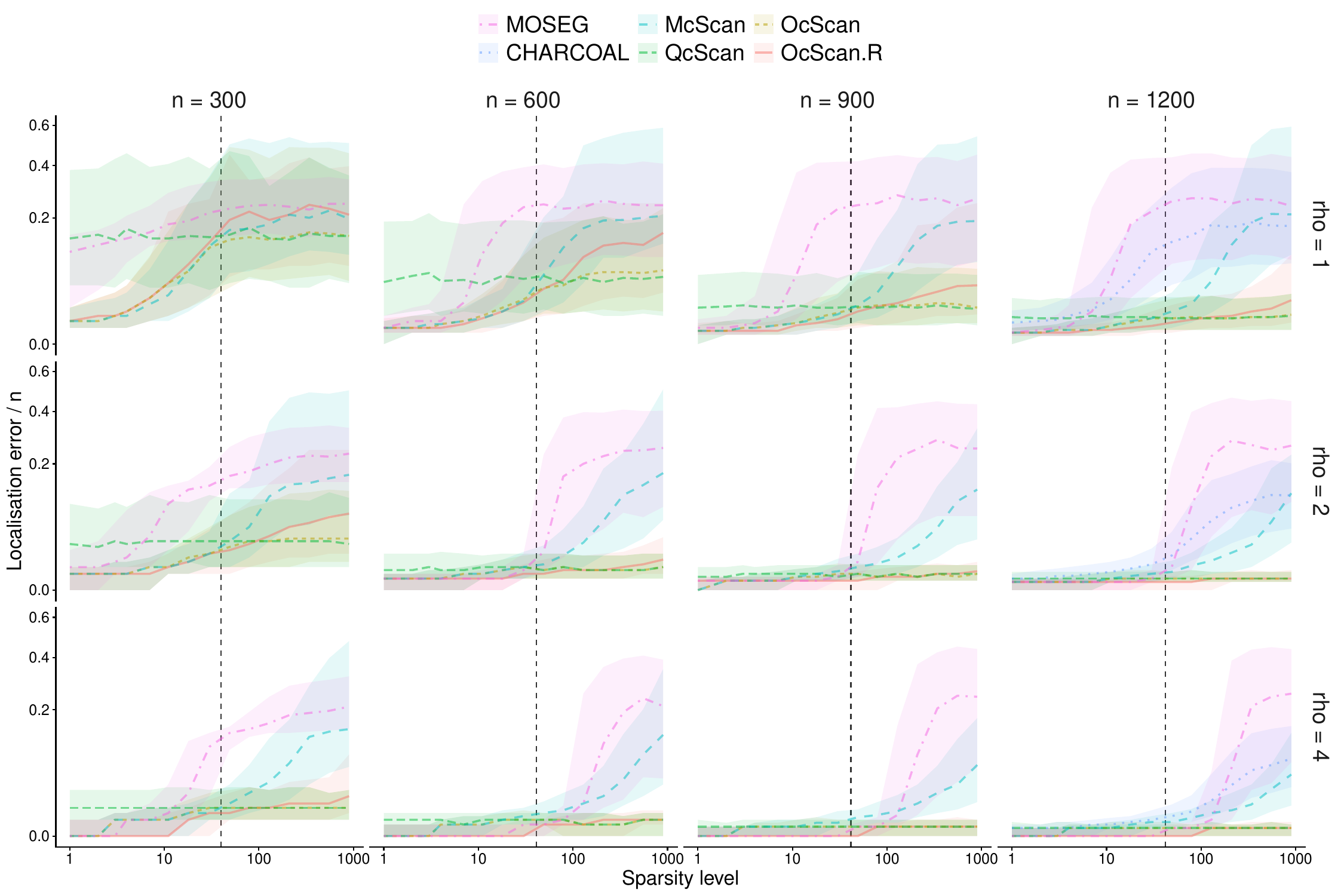}
\caption{Estimation performance of MOSEG, CHARCOAL, McScan, QcScan, OcScan and OcScan.R in (M1); CHARCOAL is only applicable in the last column where $n = 1200 > p = 900$. 
The $x$-axis denotes the inherent sparsity $\mathfrak{s}$ which coincides with $\mathfrak{s}_\delta$. 
In each scenario, the results are based on 1000 repetitions, with median error curves shown alongside shaded regions representing the interquartile range. 
The vertical dashed lines mark where $\mathfrak{s} = \sqrt{p{\log\log(n)}}$. 
The $x$-axis is shown on a log scale, and the $y$-axis on a squared root scale. \label{f:iderr}}
\end{figure}

\begin{figure}[h!t!b!]
\centering
\includegraphics[width=\linewidth]{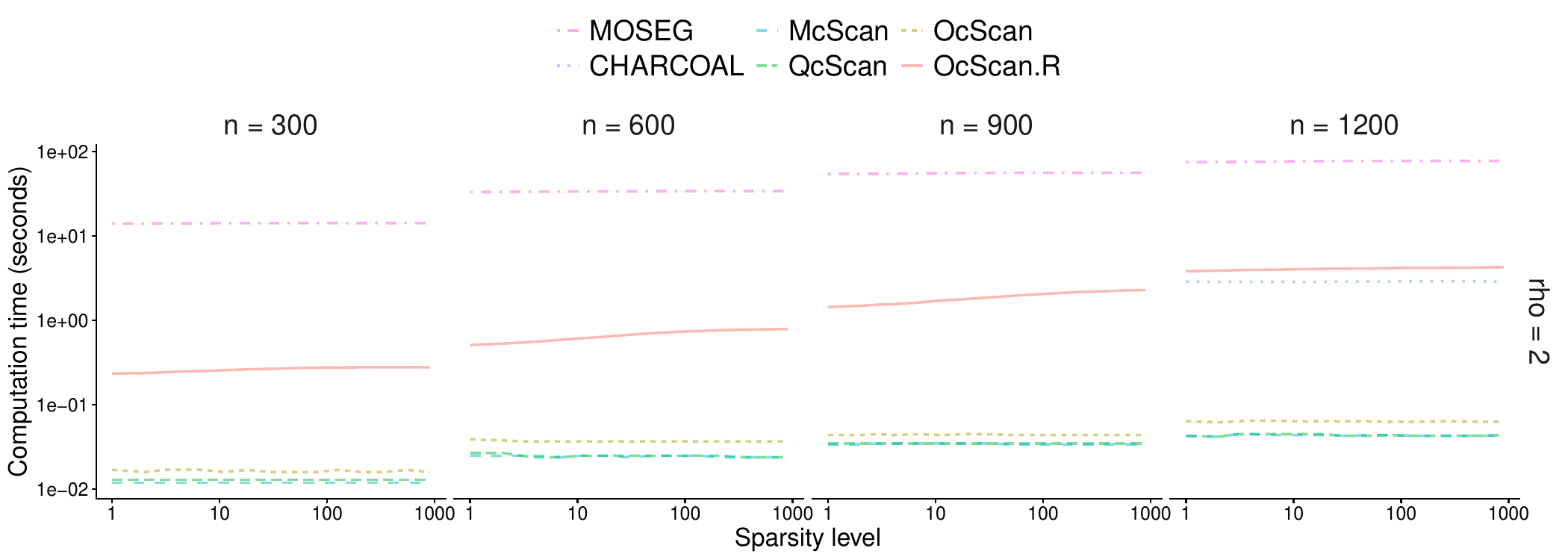}
\caption{Computation time of MOSEG, CHARCOAL, McScan, QcScan, OcScan and OcScan.R in (M1); CHARCOAL is only applicable in the last column where $n = 1200 > p = 900$. The curves plot the median runtimes (on 1.0 GHz processors) over 1000 repetitions when $\rho = 2$. The runtimes for $\rho \in\{1,4\}$ are very similar. Both axes are on a log scale. 
\label{f:idtm}}
\end{figure}

\paragraph{(M1) Isotropic Gaussian design with varying sparsity in regression coefficients.} 
We set $\bm\Sigma = \mbf I$ with $p = 900$, and $\bm{\beta}_0 = -\bm{\beta}_1 = \rho \bm\delta$ with $\rho \in \{1, 2, 4\}$, where $\bm\delta = (\delta_1, \ldots, \delta_p)^\top$ with $\delta_i = 0$ for $i \notin \mc S \subseteq [p]$. The support $\mc{S}$ is sampled uniformly with cardinality $|\mc{S}| = \mathfrak{s}$, and $(\delta_i, \, i \in \mc S)$ is drawn uniformly from the unit sphere in $\mathbb{R}^{\mathfrak{s}}$. The sparsity level~$\mathfrak{s}$ ranges over $ \{1,  2,   3,   4,   7,  11,  18,  30,  49,  79, 129, 210, 341, 554, 900\}$, a subset of $[p]$ with values approximately evenly spaced on a log scale. We consider sample sizes $n \in \{300, 600, 900, 1200\}$ with the change point located at $\cp = n/4$; see Figure~\ref{f:iderr} that expands on Figure~\ref{f:iderrex}. 
As noted in Section~\ref{sec:isotropic}, we observe a clear phase transition around $\mathfrak{s} = \sqrt{p{\log\log (n)}}$, which marks where the performance of the methods designed for handling the sparse regime (McScan, MOSEG and CHARCOAL), is overtaken by that of QcScan that performs uniformly across all $\mathfrak{s}$; OcScan and OcScan.R are able to enjoy the good performance of McScan and QcScan adaptively to the unknown sparsity, with the latter performing marginally better in the sparse regime where $\bm\delta$ itself is also sparse. 
Figure~\ref{f:idtm} additionally confirms that the good performance of the proposed covariance scanning methods is achieved with a fraction of time taken by MOSEG or CHARCOAL. OcScan.R involves a Lasso-type estimation of $\bm\delta$ and thus is relatively slower, but still is much faster than MOSEG that requires the computation of Lasso estimators of local regression parameters over time.

\begin{figure}[h!t!b!]
\centering
\includegraphics[width=\linewidth]{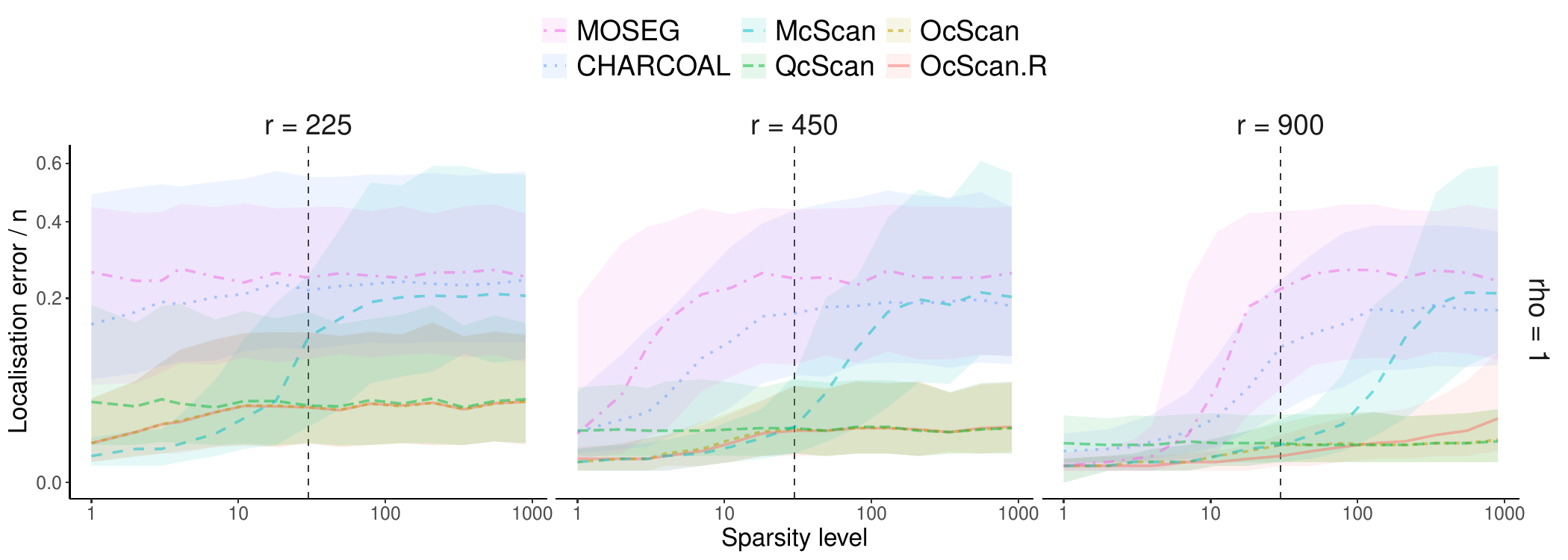}
\caption{Estimation performance of MOSEG, CHARCOAL, McScan, QcScan, OcScan and OcScan.R in modified (M1) where $n = 1200$, $p = 900$ and $\bm\Sigma = \mbf U_r \mbf U_r^\top$ with varying $r \in \{p/4, p/2, p\}$, where $\mbf U_r \in \R^{p \times r}$ is randomly generated with orthonormal columns. 
The $x$-axis denotes $\mathfrak{s}_\delta = \vert \bm\delta \vert_0$. 
In each scenario, the results are based on 1000 repetitions, with median error curves shown alongside shaded regions representing the interquartile range. 
The vertical dashed lines mark where $\mathfrak{s}_\delta = \sqrt{p{\log\log(n)}}$. 
The $x$-axis is shown on a log scale, and the $y$-axis on a squared root scale.}
\label{f:rank}
\end{figure}

We further consider the case where $\bm\Sigma$ is rank deficient, which is accommodated by our theoretical framework but is typically excluded by the standard assumptions in the existing literature. For this, we retain the same data generation process as described above, with the modification pertaining to the construction of $\mbf x_t$ only: We generate $\mbf x_t = \mbf U_r \wt{\mbf x}_t$ for some $r \le p$, where $\mbf U_r \in \R^{p \times r}$ is a randomly generated matrix with orthonormal columns, and independently $\wt{\mbf x}_t \sim_{\iid} \mathcal{N}_r(\bm 0, \mbf I)$. It follows that $\bm\Sigma = \mbf U_r \mbf U_r^\top$ which is of rank $r$. As illustrated in \Cref{f:rank}, the proposed covariance scanning methods exhibit stable performance across varying $r \in \{p/4, p/2, p\}$ and their regime-specific behaviour is comparable to that observed when $\bm\Sigma = \mbf I$.
On the other hand, the performance of MOSEG and CHARCOAL grows worse as $r$ decreases, even when $\bm\delta$ is highly sparse and $n = 1200 > p = 900$. Additional results are presented in \Cref{app:m1}.

\paragraph{(M2) Toeplitz Gaussian design with varying sparsity in differential parameters.} 
We set $\bm\Sigma = [\gamma^{\vert i - j \vert}]_{i, j = 1}^p$ for $\gamma \in \{-0.6,\, 0,\, 0.6\}$, with the convention $0^0 = 1$, and define $\bm\beta_0 = \bm\mu - \bm\delta/2$ and $\bm\beta_1 = \bm\mu + \bm\delta/2$. For each realisation, we generate a non-sparse average regression vector $\bm\mu = \nu \cdot \bm\mu_\circ / \vert\bm\Sigma^{1/2}\bm\mu_\circ\vert_2$ for $\bm\mu_\circ \sim \mc N_p(\mbf 0, \mbf I)$ and $\nu \in \{0,0.5,1\}$,  and a \enquote{sparse} differential parameter vector $\bm\delta = \rho \cdot \bm\delta_0/\vert \bm\Sigma^{1/2}\bm\delta_0\vert_2$ for $\rho \in \{2, 4\}$, where we consider two types of sparse structure in $\bm\delta_0$ as follows: 
\begin{description}[wide, leftmargin=0pt]
\item[Standard sparsity.] The vector $\bm\delta_0$ itself has $\mathfrak{s}_{\bm\delta}$ non-zero elements taking values from $\{1, - 1\}$ at random locations.
\item[Inherent sparsity.] The covariance weighted vector $\bm\Sigma^{1/2}\bm\delta_0$ has $\mathfrak{s}$ non-zero elements taking values from $\{1, - 1\}$ at random locations.
\end{description}
The sparsity levels range over $ \{1,  2,  4,  7, 14,  27,  53, 103, 200\}$, a subset of $[200]$ with values approximately evenly spaced on a log scale. In the special case of $\gamma = 0$, the two types of sparsity coincide and, if further $\nu = 0$, the scenario is similar to (M1) with the only difference in the generation of non-zero entries. 
We fix $n = 300$ and $\cp = 75$, and vary $p \in \{200, 400\}$. See \Cref{f:toep_0p6_p200_rho2,f:toep_0p6_p200,f:toep_0}, and also those given in Appendix~\ref{app:m2}. % \Cref{f:toep_0p6_p400,f:toep_n0p6_p200,f:toep_n0p6_p400} in the appendix. 
As discussed in Section~\ref{sec:general}, we observe that the inherent sparsity better captures the phase transition demonstrated by the performance of methods geared towards the sparse and dense regimes, respectively.
Regardless of how the sparsity level is determined, OcScan is shown to adaptively enjoy the good performance of McScan or QcScan. % depending on which method outperforms the other. 
The refinement step (Section~\ref{sec:refine}) adopted by OcScan.R appears to be effective in the sparse regime (both in terms of $\mathfrak{s}_\delta$ and $\mathfrak{s}$) particularly as the size of change increases with $\rho$.

\begin{figure}[h!t!b!]
\centering
\includegraphics[width=\textwidth]{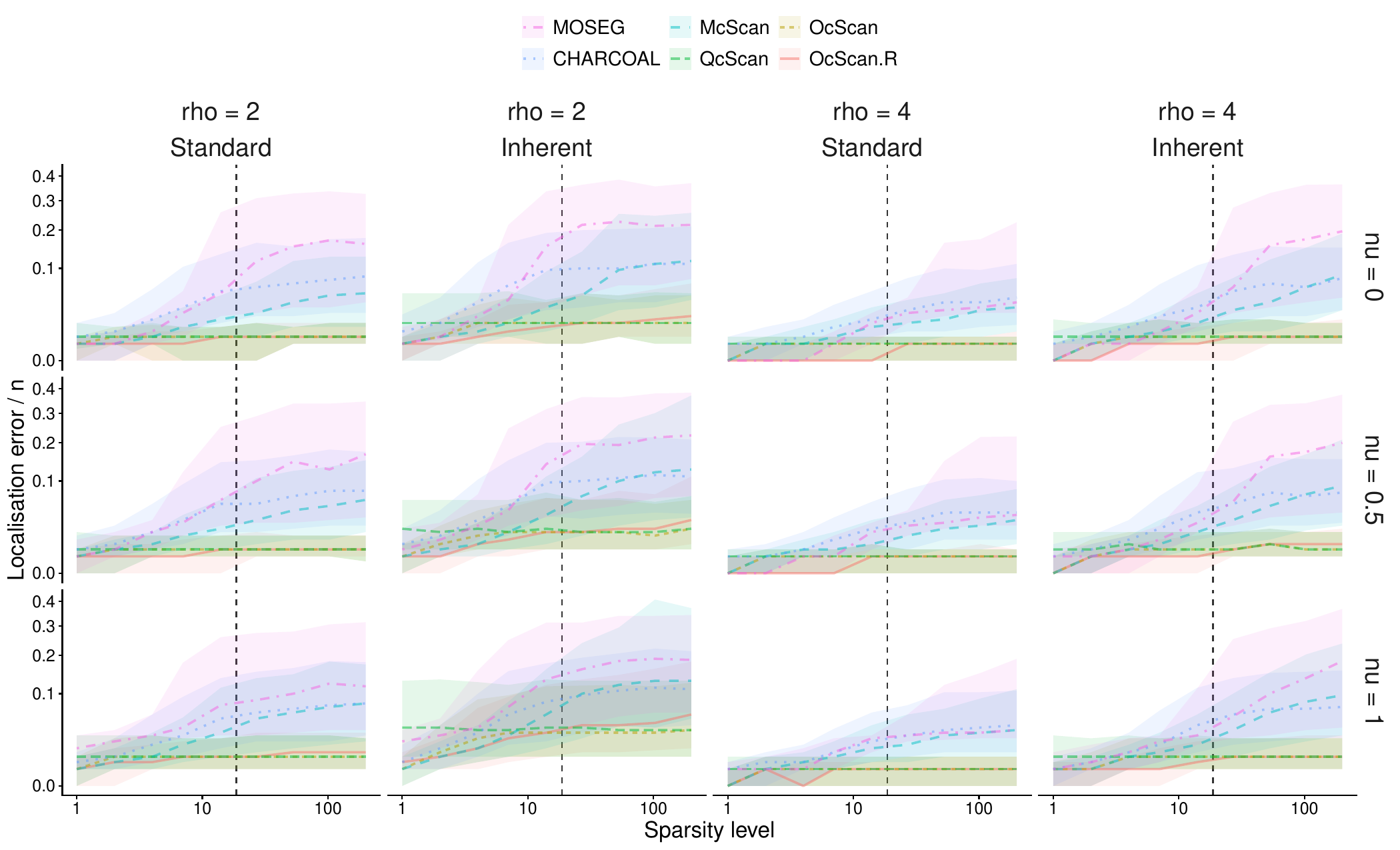}
\caption{Estimation performance of MOSEG, CHARCOAL, McScan, QcScan, OcScan and OcScan.R in (M2) with $\gamma = 0.6$ and $p = 200$.
In the left two columns, the $x$-axis denotes the standard sparsity $\mathfrak{s}_\delta = \vert \bm\delta \vert_0$ while in the right two columns, it is the inherent sparsity $\mathfrak{s} = \vert \bm\Sigma^{1/2}\bm\delta \vert_0$. 
In each scenario, the results are based on 1000 repetitions, with median error curves shown alongside shaded regions representing the interquartile range. 
The vertical dashed lines mark where $\mathfrak{s}_\delta = \sqrt{p{\log\log(n)}}$ or $\mathfrak{s} = \sqrt{p{\log\log(n)}}$. 
The $x$-axis is shown on a log scale, and the $y$-axis on a squared root scale. 
\label{f:toep_0p6_p200}}
\end{figure}

\begin{figure}[h!t!b!]
\centering
\includegraphics[width=\textwidth]{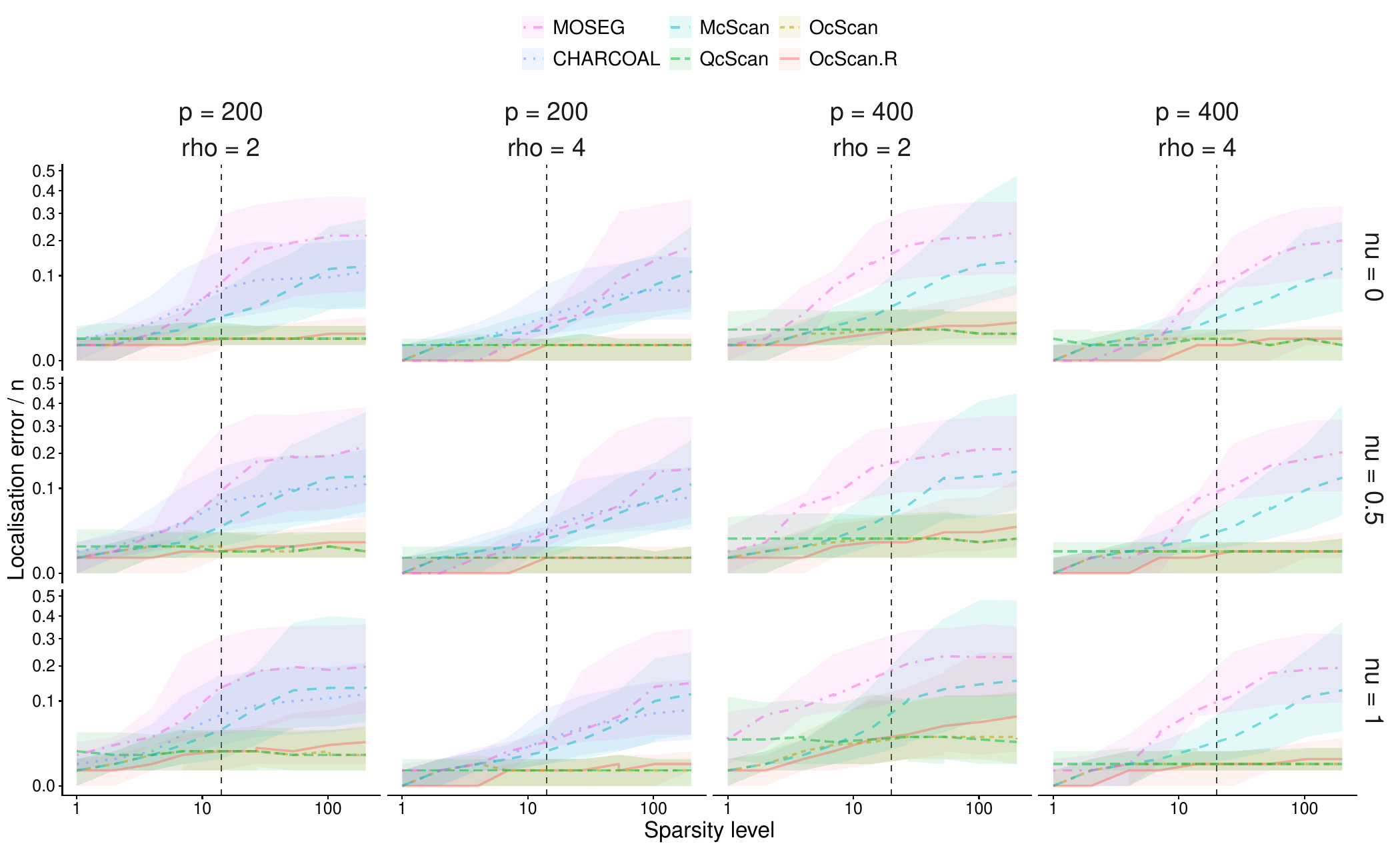}
\caption{
Estimation performance of MOSEG, CHARCOAL, McScan, QcScan, OcScan and OcScan.R in (M2) with $\gamma = 0$; CHARCOAL is only applicable in the left two columns where $p = 200 < n = 300$.
The $x$-axis denotes the inherent sparsity which coincides with the standard sparsity as $\bm\Sigma = \mbf I$ and thus $\mathfrak{s} = \vert\bm\Sigma^{1/2}\bm\delta\vert_0 = \vert\bm\delta\vert_0$. 
In each scenario, the results are based on 1000 repetitions, with median error curves shown alongside shaded regions representing the interquartile range. 
The vertical dashed lines mark where $\mathfrak{s} = \sqrt{p{\log\log(n)}}$.
The $x$-axis is shown on a log scale, and the $y$-axis on a squared root scale. 
\label{f:toep_0}}
\end{figure}

\subsection{Change point detection performance}
\label{ss:num:detect}

We investigate the performance of the proposed methods for testing the existence of a change point under the AMOC model~\eqref{eq:amoc}. 
We set the data-adaptive thresholds as $\zeta_{\Mc} = \bar{c} \wh\sigma_X \wh\Psi \sqrt{\log(p \log(n))}$ for McScan and $\zeta_{\Qc} = c \Vert \wh{\bm\Sigma} \Vert \wh\Psi^2 \sqrt{p\log\log(n)}$ for QcScan (see Corollary~\ref{cor:adapt:short}). 
The constants $\bar{c} = 1.3$ and $c = 0.7$ are chosen based on the simulation studies reported in \Cref{app:thd}. The same thresholds are used for the McScan and QcScan components within OcScan.
With the input bandwidth set at $50$, MOSEG selects the final model via cross validation with the maximum number of change points set to be one. 
For CHARCOAL, we apply the method to determine the existence of at most one change point with the threshold $0.5\tilde{\sigma}\log(p)$, where $\tilde{\sigma}$ is an estimator of the noise level $\sigma$ proposed by \citet{Dicker14}, as recommended by the authors.

% For MOSEG, we use \texttt{moseg.cv(..., G = 50, max.cps = 1)}, where the number of change points is selected via cross validation with bandwidth $50$, under the assumption that at most one change point is present. For CHARCOAL, we apply \texttt{cpreg(...)} to detect a single change point and determine its existence by comparing the detector statistic to the threshold $0.5\tilde{\sigma}\log(p)$, where $\tilde{\sigma}$ is an estimator of the noise level $\sigma$ proposed by \citet{Dicker14}, as recommended in the original paper.

We focus on model (M2) with $\nu = 0$, $n = 300$ and $p = 200$, %\footnote{previously $n = 200$, $p = 300$ { In (M2) in \Cref{ss:loc:err}, we used $n = 300$ and $p \in \{200, 400\}$, so it is consistent. Or you meant somewhere else?} I just meant that there were typos :) and wanted your confirmation}
and consider both the null ($q= 0, \cp = n$) and the alternative ($q = 1, \cp= n/4$) scenarios.
We vary the intrinsic sparsity as $\mathfrak{s} = \vert \bm\Sigma^{1/2} \bm\delta \vert_0 = \vert \bm\Sigma^{1/2} \bm\beta_0 \vert_0 \in \{1, 4, 14, 53, 200\}$, with the grid chosen to be approximately evenly spaced on the logarithmic scale.

\begin{figure}[h!t!]
\centering
\includegraphics[width = \linewidth]{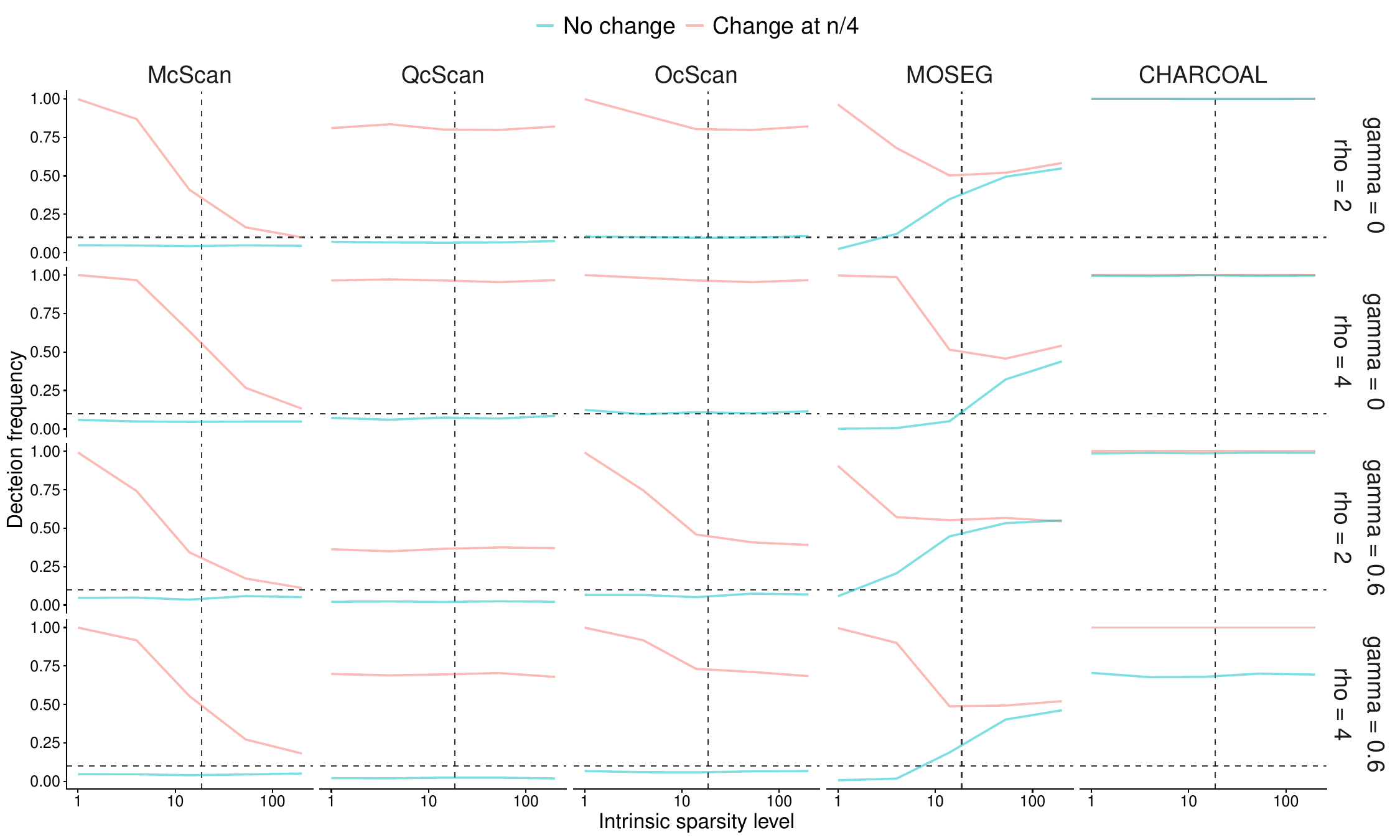}
\caption{Detection frequencies of McScan, QcScan, OcScan, MOSEG and CHARCOAL under (M2) with $\nu = 0$, $n = 300$, and $p = 200$. In each scenario, detection frequencies are computed from 1000 repetitions. The $x$-axis, shown on a log scale, represents the intrinsic sparsity $\mathfrak{s} = \vert \bm\Sigma^{1/2}\bm\delta \vert_0$, and the vertical dashed line marks $\mathfrak{s} = \sqrt{p\log\log(n)}$. The horizontal dashed lines correspond to the $10\%$ level.}
\label{f:detect:freq}
\end{figure}

The detection frequencies are reported in \Cref{f:detect:freq}. As expected, McScan performs well in sparse regimes, whereas QcScan is more effective in dense settings. OcScan combines the strengths of both procedures and achieves strong performance across all sparsity levels. MOSEG performs well in sparse regimes, although it is slightly outperformed by McScan; in dense regimes, it tends to produce many false positives. 
CHARCOAL tends to struggle in distinguishing between the presence and absence of a change point, consistently reporting high detection frequencies regardless of whether a change point is present.

\subsection*{Acknowledgements}

We would like to thank Prof.\ Richard Samworth for a valuable comment.
% and gratefully acknowledge the editor, the associate editor and the two anonymous reviewers for their constructive comments on an earlier version of the paper. 
HC's work is supported by Engineering and Physical Sciences Research Council (EP/Z531327/1).
HL is supported by the Deutsche Forschungsgemeinschaft (DFG, German Research Foundation) under Germany's Excellence Strategy (EXC 2067/1-39072994).
\clearpage

\bibliographystyle{apalike}
\bibliography{ref}

\clearpage

\appendix

\numberwithin{equation}{section}
\numberwithin{figure}{section}
\numberwithin{table}{section}
\numberwithin{thm}{section}

\section{Connection to the covariance change detection problems}
\label{sec:cov}

% In Introduction, we have extensively discussed existing results on change point analysis in linear models. 
The change point problem in linear models can be viewed as a special case of changes in covariances. Specifically, the AMOC model~\eqref{eq:amoc} is equivalent to  
\[
\bmx
\mbf x_t\\
Y_t
\emx
\sim_{\mathrm{ind.}} 
\begin{cases}
\mc N_{p+1}\l(\mbf 0, \bm\Sigma_0\r) & \text{for \ } 1 \le t \le \cp, \\
\mc N_{p+1}\l(\mbf 0, \bm\Sigma_1\r) & \text{for \ } \cp + 1 \le t \le n, 
\end{cases}
\quad\text{ with }\quad \bm\Sigma_i =\bmx\bm\Sigma & \bm\Sigma\bm\beta_i \\
\bm\beta_i^\top\bm\Sigma & \bm\beta_i^\top \bm\Sigma\bm\beta_i  + \sigma^2\emx.
\]
This perspective allows us to draw connections with the literature on covariance change detection in high dimensions.

\subsection{Change size measured in operator norm}

\citet{wang2021optimal} proposed a projection-based binary segmentation method and established its consistency.
Their Assumption~3 on detection requirement translates to our setting~as % follows:
\begin{align}
\label{eq:cov:dlb}
\Delta \kappa^2 \gtrsim  B^2 p\log(n)
\end{align}
where $B = \max(\Vert \bm\Sigma \Vert, \Psi^2)$ and $\kappa = \Vert \bm\Sigma_1 - \bm\Sigma_0 \Vert$, which satisfies 
\begin{align*}
\max(\vert \bm\Sigma \bm\delta \vert_2^2, \vert \bm\beta_0^\top \bm\Sigma \bm\beta_0 - \bm\beta_1^\top \bm\Sigma \bm\beta_1 \vert^2) \le \kappa^2 &\le 2\max(\vert \bm\Sigma \bm\delta \vert_2^2, \vert \bm\beta_0^\top \bm\Sigma \bm\beta_0 - \bm\beta_1^\top \bm\Sigma \bm\beta_1 \vert^2).
% \\
% &\le 2\max(\vert \bm\Sigma \bm\delta \vert_2^2, 2 \Psi \vert \bm\Sigma^{1/2} \bm\delta \vert_2).
\end{align*}
To ease comparison, assume that $\Var(Y_t)$ remains invariant before and after the change, i.e.\ 
\begin{equation}\label{e:ev}
\bm\beta_0^\top \bm\Sigma \bm\beta_0 = \bm\beta_1^\top \bm\Sigma \bm\beta_1.
\end{equation} 
Ignoring logarithmic factors, our detection requirement~\eqref{eq:sparse:dlb} is weaker than~\eqref{eq:cov:dlb} by a factor of (at least) $p / \min(\Vert \bm\Sigma \Vert \sigma_*^{-1} \mathfrak{s}, \mathfrak{s}_*)$, with $\mathfrak{s}_* = \vert\bm\Sigma \bm\delta\vert_0$, and~\eqref{eq:dense:dlb} is weaker by a factor of $\sqrt{p}$. 
Thus, with OcScan (\Cref{sec:adaptive}), we achieve an improvement by a factor of 
\begin{align*}
\frac{p}{\min(\Vert \bm\Sigma \Vert \sigma_*^{-1} \mathfrak{s}, \, \mathfrak{s}_*, \, \sqrt{p})}.
\end{align*}
It is important to note that the requirement~\eqref{eq:cov:dlb} is minimax optimal, up to $\log(n)$, for the covariance change point problem in general \citep[Lemma~3.1]{wang2021optimal}.
The improved detection performance of OcScan arises from leveraging the specific structure of the linear model, which also enables the careful analysis of the phase transition with respect to the inherent sparsity $\mathfrak{s}$.
Additionally, the estimation rate achieved by our refinement method described in \Cref{sec:refine}, nearly matches the minimax lower bound of order $\kappa^{-2} B^2$ reported in \citet[Lemma~3.2]{wang2021optimal}.

\subsection{Element-wise sparsity of the change}

\citet{DePaYa22} investigated the covariance change point estimation problem by imposing the sparsity directly on $\bm\Sigma_1 - \bm\Sigma_0$.
Let $(\bm\Sigma)_{i, i'}$ denotes the element in the $i$th row and $i'$th column of $\bm\Sigma$.
Then for estimation consistency, their proposed method requires that 
$\vert\bm\Sigma_j\vert_{\infty}\lesssim 1$, $(\bm\Sigma_j)_{i,i} \lesssim 1$, for $i \in[p+1]$ and $j\in\{0,1\}$, and 
$$
\min\bigl\{\vert(\bm\Sigma_0 - \bm\Sigma_1)_{i,j}\vert :  (\bm\Sigma_0 - \bm\Sigma_1)_{i,j} \neq 0\bigr\} \Delta^2  \; \gtrsim\; n^{3/2}\log(p\vee n),
$$
which translates to $\min\bigl\{\vert(\bm\Sigma\bm\delta)_{i}\vert : (\bm\Sigma\bm\delta)_{i} \neq 0\bigr\} \Delta^2 \gtrsim n^{3/2} \log(p\vee n)$, in our setting. 
This is stricter than our requirement~\eqref{eq:sparse:dlb}, which is $\vert\bm\Sigma\bm\delta\vert_{\infty}^2 \Delta \gtrsim \log(p \log(n))$, since  $1/\sqrt{\Delta} \le n^{3/2}/\Delta^2$ due to $\Delta \le n$.

\subsection{Sparsity in the principal direction of the change}

\cite{moen2024minimax} investigated the minimax optimality in covariance change point testing over the class of parameters that are \textit{$d$-sparse} for any $d \in [p]$, which fulfil
\begin{align*}
\Vert \bm\Sigma_1 - \bm\Sigma_0 \Vert = \Lambda^d_{\max}(\bm\Sigma_1 - \bm\Sigma_0) := \sup_{\substack{\mbf v: \, \vert \mbf v \vert_2 = 1, \\ \vert \mbf v \vert_0 = d}} \vert \mbf v^\top (\bm\Sigma_1 - \bm\Sigma_0) \mbf v \vert,
\end{align*}
i.e.\ the leading eigenvector of $\bm\Sigma_1 - \bm\Sigma_0$ is sparse with $d$ non-zero elements.
Their (computationally intractable) test matches the minimax lower bound on the detection boundary, which requires
\begin{align}
\min\l\{\frac{\Lambda^d_{\max}(\bm\Sigma_1 - \bm\Sigma_0)}{(\Vert \bm\Sigma_0 \Vert \vee \Vert \bm\Sigma_1 \Vert) - \Lambda^d_{\max}(\bm\Sigma_1 - \bm\Sigma_0)}, \l( \frac{\Lambda^d_{\max}(\bm\Sigma_1 - \bm\Sigma_0)}{(\Vert \bm\Sigma_0 \Vert \vee \Vert \bm\Sigma_1 \Vert) - \Lambda^d_{\max}(\bm\Sigma_1 - \bm\Sigma_0)} \r)^2 \r\} \Delta
\nn \\
\gtrsim\; d\log(ep/d) \vee \log\log(8n)
\nn
\end{align}
for \textit{some} $d$ which, combined with that $\bm\Sigma_j, \, j = 0, 1$, are positive definite, simplifies to
\begin{align}
\l( \frac{\Lambda^d_{\max}(\bm\Sigma_1 - \bm\Sigma_0)}{\Vert \bm\Sigma_0 \Vert \vee \Vert \bm\Sigma_1 \Vert} \r)^2 \Delta\; \gtrsim\; d\log(ep/d) \vee \log\log(8n).
\label{eq:moen}
\end{align}
Some elementary calculations show that
\begin{align*}
\Vert \bm\Sigma_0 \Vert \vee \Vert \bm\Sigma_1 \Vert & \ge \max(\Psi^2, \Vert \bm\Sigma \Vert), \text{ \ and }
\\
\Lambda^d_{\max}(\bm\Sigma_1 - \bm\Sigma_0) & \le \Vert \bm\Sigma_1 - \bm\Sigma_0 \Vert \le \Vert \bm\Sigma_1 - \bm\Sigma_0 \Vert_F \lesssim \vert \bm\Sigma \bm\delta \vert_2 + \vert \bm\beta_1^\top \bm\Sigma \bm\beta_1 - \bm\beta_0^\top \bm\Sigma \bm\beta_0 \vert.
% \\
% & \Lambda^d_{\max}(\bm\Sigma_1 - \bm\Sigma_0) \le d \vert \bm\Sigma \bm\delta \vert_\infty,
\end{align*}
To facilitate the comparison between~\eqref{eq:moen} and our results, we relate the notion of $d$-sparsity and the sparsity of $\bm\Sigma \bm\delta$, as $d = \mathfrak{s}_* + 1$.
Then, again assuming~\eqref{e:ev}, our requirement that either~\eqref{eq:sparse:dlb} or~\eqref{eq:dense:dlb} holds, is weaker than~\eqref{eq:moen} by a factor of $\max(d/\sqrt{p}, 1)$, save logarithmic factors. 
That is, the detection boundary of OcScan is more favourable in the scenario where $d \ge \sqrt{p}$. 

\subsection{Variance changes in responses}\label{ss:var}
Under the AMOC model~\eqref{eq:amoc}, we observe that
\[
Y_t \sim_{\mathrm{ind.}}
\begin{cases}
\mc N(0, \sigma_0^2), & 1 \le t \le \cp, \\
\mc N(0, \sigma_1^2), & \cp + 1 \le t \le n,
\end{cases}
\qquad \text{where} \qquad
\sigma_i^2 = \bm\beta_i^\top \bm\Sigma \bm\beta_i + \sigma^2,
\]
which indicates that the change point $\cp$ is detectable as a change in $\Var(Y_t)$. However, this strategy is not generally recommended, as it is completely powerless when
$\lvert \bm\Sigma^{1/2}\bm\beta_0\rvert_2 = \lvert \bm\Sigma^{1/2}\bm\beta_1\rvert_2$,
even though $\lvert \bm\Sigma^{1/2}\bm\delta\rvert_2 \neq 0$, where $\bm\delta = \bm\beta_1 - \bm\beta_0$.

Nevertheless, we relate our detection conditions in~\eqref{eq:sparse:dlb} and~\eqref{eq:dense:dlb} to the  minimax detection boundary for univariate variance change point detection (see \citealp[Theorem~1]{moen2024minimax}),
\begin{equation}\label{e:var:bnd}
\Delta \Biggl\{
\Bigl(\frac{\lvert \sigma_0^2 - \sigma_1^2\rvert}{\min\{\sigma_0^2, \sigma_1^2\}}\Bigr)
\wedge
\Bigl(\frac{\lvert \sigma_0^2 - \sigma_1^2\rvert}{\min\{\sigma_0^2, \sigma_1^2\}}\Bigr)^2
\Biggr\}
\asymp \log\log(n).
\end{equation}
For ease of comparison, assume that
$\lvert \bm\Sigma^{1/2}\bm\beta_0\rvert_2 \neq \lvert \bm\Sigma^{1/2}\bm\beta_1\rvert_2$
(or equivalently, $\sigma_0^2 \neq \sigma_1^2$) and that $\bm\Sigma$ has its spectrum bounded away from zero and infinity. Under the latter condition, OcScan detects the change point $\theta$ whenever
\begin{equation}\label{e:simple}
\Delta\lvert\bm\Sigma^{1/2}\bm\delta\rvert_2^2 \gtrsim \Psi^2\min\l\{\mathfrak{s}\log\bigl(p\log(n)\bigr),\,\sqrt{p \log\log(n)}\r\},
\end{equation}
where $\Psi = \max\{
\lvert \bm\Sigma^{1/2}\bm\beta_0\rvert_2,
\lvert \bm\Sigma^{1/2}\bm\beta_1\rvert_2,
\sigma
\}$ and $\mathfrak{s} = \lvert\bm\Sigma^{1/2}\bm\delta\rvert_0$.
Then, we distinguish the following two regimes.

\begin{enumerate}[label=(\alph*)]
\item\label{i:same:var} \emph{Comparable variances}, i.e.\ $\sigma_0^2 \asymp \sigma_1^2$.  
In this case, \eqref{e:var:bnd} is equivalent to
\[
\Delta
\Bigl(\frac{\lvert \sigma_0^2 - \sigma_1^2\rvert}{\sigma_0^2 + \sigma_1^2}\Bigr)^2
\asymp \log\log(n).
\]

We relate $\lvert \sigma_0^2 - \sigma_1^2\rvert$ to $\lvert \bm\Sigma^{1/2}\bm\delta\rvert_2$ as follows:
\begin{align}
\lvert \sigma_0^2 - \sigma_1^2\rvert
&=
\Bigl\lvert
\lvert \bm\Sigma^{1/2}\bm\beta_1\rvert_2^2
-
\lvert \bm\Sigma^{1/2}\bm\beta_0\rvert_2^2
\Bigr\rvert \notag \\
&=
\lvert \bm\delta^\top \bm\Sigma (\bm\beta_0 + \bm\beta_1)\rvert \notag \\
&\le
\lvert \bm\Sigma^{1/2}\bm\delta\rvert_2
\bigl(
\lvert \bm\Sigma^{1/2}\bm\beta_0\rvert_2
+
\lvert \bm\Sigma^{1/2}\bm\beta_1\rvert_2
\bigr)
\label{e:ineq} \\
&\le
2 \lvert \bm\Sigma^{1/2}\bm\delta\rvert_2 \, \Psi. \notag
\end{align}

Since $\sigma_0^2 + \sigma_1^2 \asymp \Psi^2$, we obtain
\[
\Delta
\Bigl(\frac{\lvert \sigma_0^2 - \sigma_1^2\rvert}{\sigma_0^2 + \sigma_1^2}\Bigr)^2
\lesssim
\Delta \frac{\lvert \bm\Sigma^{1/2}\bm\delta\rvert_2^2}{\Psi^2}.
\]
Consequently, our requirement in~\eqref{e:simple} is always weaker than~\eqref{e:var:bnd} when
$p \lesssim \log\log(n)$, and is comparable when $\mathfrak{s} \asymp 1$ and $\log(p) \asymp \log\log(n)$. Otherwise, the boundary in~\eqref{e:var:bnd} may be weaker, although this ultimately depends on the tightness of the inequality in~\eqref{e:ineq}.

\item \emph{Unbalanced variances}, i.e.\ $\sigma_0^2 \gg \sigma_1^2$ or $\sigma_0^2 \ll \sigma_1^2$.  
WLOG, consider $\sigma_0^2 \ll \sigma_1^2$. In this regime, \eqref{e:var:bnd} simplifies to
\[
\Delta \frac{\lvert \sigma_0^2 - \sigma_1^2\rvert}{\sigma_0^2}
\asymp \log\log(n).
\]

Because $\sigma_0^2 \ll \sigma_1^2$ implies
$\lvert \sigma_0^2 - \sigma_1^2\rvert \asymp \sigma_1^2
\asymp \lvert \bm\Sigma^{1/2}\bm\delta\rvert_2^2$,
we obtain
\[
\Delta \frac{\lvert \sigma_0^2 - \sigma_1^2\rvert}{\sigma_0^2}
\asymp
\Delta \frac{\lvert \bm\Sigma^{1/2}\bm\delta\rvert_2^2}{\sigma_0^2}
\gg
\Delta \frac{\lvert \bm\Sigma^{1/2}\bm\delta\rvert_2^2}{\Psi^2}.
\]
Hence, in this setting, the boundary~\eqref{e:var:bnd} is weaker than~\eqref{e:simple} whenever
$p \gtrsim \log\log(n)$.
\end{enumerate}

Overall, this comparison suggests that additionally scanning for changes in $\Var(Y_t)$ can be advantageous, particularly when $\lvert \bm\Sigma^{1/2}\bm\beta_0\rvert_2$ and $\lvert \bm\Sigma^{1/2}\bm\beta_1\rvert_2$ differ substantially. Note, however, that the change point literature for high-dimensional linear models (including the present paper) primarily focuses on the more challenging regime described in~\ref{i:same:var}.

\section{Proofs}

\subsection{Proofs for the results in Section~\ref{sec:minimax}}

\subsubsection{Proof of Theorem~\ref{th:gaussid}}

We define 
$$
\mc{M} = \bigl\{\lfloor\log_2(n^{2\nu})\rfloor, \lfloor\log_2(n^{2\nu})\rfloor+1,\ldots, \lfloor \log_2(n/4) \rfloor\bigr\}
$$ 
and choose, with $\underline{c}_1 = 1/6$,
\[
\kappa^2 = \underline{c}_1\sigma^2\mathfrak{s} \log\Bigl(1 + \tfrac{p {\log\log(n)}}{\mathfrak{s}^2}\vee\tfrac{\sqrt{p {\log\log(n)}}}{\mathfrak{s}}\Bigr).
\]
Let $\mathbf{1}^{\pm}_S$ denote a $p$-dimensional vector with its $i$th entry being $1$ or $-1$ if $i \in S$ and being zero if $i \not\in S$, and introduce $\mathcal{U} = \bigl\{\bm{u} = ({\kappa}/{\sqrt{\mathfrak{s}}})\mathbf{1}^{\pm}_S \in \R^p:  S \subseteq [p],  \; \vert S \vert = \mathfrak{s} \bigr\}$. Define
\[
\p^{(1)} = \frac{1}{\vert\mc U\vert\,\vert\mc M\vert} \sum_{\bm{u} \in \mathcal{U}} \sum_{m \in \mc{M}}\p_{2^m,  2^{-{m}/{2}}\bm u , \bm 0}\,.
\]
Since $\p_{2^m,  2^{-{m}/{2}}\bm u , \bm 0}$ lies in $\mc P^{\mathfrak{s}, n, p}_{\mbf I, \sigma}(\tau)$, by the method of fuzzy hypotheses (see e.g.\ Section~2.7.4 in \citealp{Tsy09}), we have
$$
\inf_{\psi}\biggl\{  \sup_{\p \in \mc P_{\mbf I, \sigma}^{0, n, p}} \E_{\p}\l[\psi\l((Y_t, \mbf x_t)_{t\in[n]}\r)\r]  +  \sup_{\p \in \mc P_{\mbf I, \sigma}^{\mathfrak{s}, n, p}(\tau)} \E_{\p}\l[1-\psi\l((Y_t, \mbf x_t)_{t\in[n]}\r)\r]\biggr\} \ge \frac{\exp\l(-\chi^2(\p^{(1)},\p^{(0)})\r)}{4},
$$
where  $\p^{(0)} \in \mc P_{\mbf I, \sigma}^{0, n, p} $ is the distribution of $\{(\mbf x_t,  Y_t)\}_{t\in [n]}$ such that  $\mbf x_t \sim_{\iid} \mc N_p(\mbf 0, \mbf I)$ and
independently  $Y_t \sim_{\iid}\mc N (0,\sigma^2)$.  
Define $\p_{\bm u}$ as in Lemma~\ref{l:cross}, and $Z_k$ as the sum of $k$ independent Rademacher random variables. Then, by Lemma~\ref{l:cross},  we have 
\begin{align*}
& \chi^2(\p^{(1)},\p^{(0)}) + 1 \\
=\, &\frac{1}{\vert\mc U\vert^2\,\vert\mc M\vert^2}\sum_{m, \ell \in \mc M}\sum_{\bm u, \bm v \in \mc U} \l[\E_{\p_{\bm0}}\l(\frac{\mathrm{d}\p_{2^{-{m}/{2}}\bm u}}{\mathrm{d}\p_{\bm 0}}\frac{\mathrm{d}\p_{2^{-{\ell}/{2}}\bm v}}{\mathrm{d}\p_{\bm 0}}\r)\r]^{2^m\wedge 2^\ell}  \\
% =\, &\frac{1}{\vert\mc M\vert^2{p \choose \mathfrak{s}}^{2}}\sum_{m, \ell \in \mc M} \sum_{S, T \subseteq [p],\, \vert S\vert = \vert T\vert = \mathfrak{s}} \E \l[\l(\l(1-\frac{Z_{\vert S \cap T\vert}}{\mathfrak{s}}\frac{\kappa^2}{\sigma^2 2^{(m+\ell)/2}}\r)^2 - \l(\frac{\kappa^2}{\sigma^2 2^{(m+\ell)/2}}\r)^2\r)^{-2^{(m\wedge \ell) -1}}\r] \\
 =\, & \frac{1}{\vert\mc M\vert^2{p \choose \mathfrak{s}}}\sum_{m, \ell \in \mc M} \sum_{S \subseteq [p],\, \vert S\vert = \mathfrak{s}} \E\l[\l(\l(1-\frac{Z_{\vert S \cap [\mathfrak s]\vert}}{\mathfrak{s}}\frac{\kappa^2}{\sigma^22^{(m+\ell)/2}}\r)^2 - \l(\frac{\kappa^2}{\sigma^22^{(m+\ell)/2}}\r)^2\r)^{-2^{(m\wedge \ell) -1}}\r].
\end{align*}
Let $\rho := \kappa^2 / \sigma^2$ and then %\footnote{Here we require \(4\kappa^2 \le \sigma^2\).} 
$\rho \le 2^{\lfloor\log_2( n^{2\nu})\rfloor}/3$. Then, by Lemma~\ref{l:basic}, for $m,\ell\in\mc M$,
\begin{align*}
& \E\l[\l(\l(1-\frac{Z_{\vert S \cap [\mathfrak s]\vert}}{\mathfrak{s}}\frac{\kappa^2}{\sigma^22^{(m+\ell)/2}}\r)^2 - \l(\frac{\kappa^2}{\sigma^22^{(m+\ell)/2}}\r)^2\r)^{-2^{(m\wedge \ell) -1}}\r] \\
\le\, &\E\l[\l(1-\frac{2\rho}{2^{(m+\ell)/2}}\cdot\frac{Z_{\vert S \cap [\mathfrak s]\vert}}{\mathfrak{s}}- \frac{\rho^2}{2^{m+\ell}}\r)^{-2^{(m\wedge \ell) -1}}\r] \\
\le\, & \exp\l(\frac{\rho^2}{2^{m\vee \ell}}\r)\E\l[\l(1-\frac{2^{1-(m+\ell)/2}\rho}{1-2^{-(m+\ell)}\rho^2}\cdot\frac{ Z_{\vert S \cap [\mathfrak s]\vert}}{\mathfrak{s}}\r)^{-2^{(m\wedge \ell) -1}}\r].
\end{align*}
By the symmetry of Rademacher distribution, we can find an event $A$ with $\p(A) = 1/2$ such that $A \subseteq \{Z_{\vert S \cap [\mathfrak s]\vert} \ge 0\}$ and $A^c \subseteq \{Z_{\vert S \cap [\mathfrak s]\vert} \le 0\}$. As $2^{1-(m+\ell)/2}\rho/(1-2^{-(m+\ell)}\rho^2) \le 3/4$, we apply again Lemma~\ref{l:basic} and obtain, for $m,\ell\in\mc M$,
\begin{align*}
&\E\l[\l(1-\frac{2^{1-(m+\ell)/2}\rho}{1-2^{-(m+\ell)}\rho^2}\cdot\frac{ Z_{\vert S \cap [\mathfrak s]\vert}}{\mathfrak{s}}\r)^{-2^{(m\wedge \ell) -1}}\r] \\ 
\le \; &  \E\l[\mathbb{I}_{A}\cdot\exp\l( 2^{m\wedge \ell} \cdot \frac{2^{1-(m+\ell)/2}\rho}{1-2^{-(m+\ell)}\rho^2}\cdot\frac{Z_{\vert S \cap [\mathfrak s]\vert}}{\mathfrak{s}}\r)\r] \\
&\qquad\quad+  \E\l[\mathbb{I}_{A^c}\cdot\exp\l( -\frac{ 2^{m\wedge \ell} }{4} \cdot \frac{2^{1-(m+\ell)/2}\rho}{1-2^{-(m+\ell)}\rho^2}\cdot\frac{Z_{\vert S \cap [\mathfrak s]\vert}}{\mathfrak{s}}\r)\r]\\
\le \; &   \l(\frac{1}{2} \E\l[\exp\l( 2^{(m\wedge \ell)+1}  \cdot \frac{2^{1-(m+\ell)/2}\rho}{1-2^{-(m+\ell)}\rho^2}\cdot\frac{Z_{\vert S \cap [\mathfrak s]\vert}}{\mathfrak{s}}\r)\r]\r)^{1/2} \\
& \qquad\quad + \l(\frac{1}{2} \E\l[\exp\l(-2^{(m\wedge \ell)-1} \cdot \frac{2^{1-(m+\ell)/2}\rho}{1-2^{-(m+\ell)}\rho^2}\cdot\frac{Z_{\vert S \cap [\mathfrak s]\vert}}{\mathfrak{s}}\r)\r]\r)^{1/2}\\
%=  \; & \frac{\sqrt{2}}{2} \l[\cosh\l( \frac{2^{(m\wedge \ell) +1} }{\mathfrak{s}}\cdot\frac{2^{1-(m+\ell)/2}\rho}{1-2^{-(m+\ell)}\rho^2}\r)\r]^{\vert S \cap [\mathfrak s]\vert/2} + \frac{\sqrt{2}}{2}\l[\cosh\l( \frac{2^{(m\wedge \ell) - 1} }{\mathfrak{s}}\cdot\frac{2^{1-(m+\ell)/2}\rho}{1-2^{-(m+\ell)}\rho^2}\r)\r]^{\vert S \cap [\mathfrak s]\vert/2}\\
\le \; & \sqrt{2} \l[\cosh\l( \frac{9\rho  }{2^{\vert m-\ell\vert/2+1}\mathfrak{s}}\r)\r]^{\vert S \cap [\mathfrak s]\vert/2},
\end{align*}
where the second inequality is by Hölder's inequality. Note that the summation involving $\vert S \cap [\mathfrak s]\vert$ can be seen as the expectation of the hypergeometric distribution with parameters $(p, \mathfrak{s}, \mathfrak{s})$. Let $\xi \sim \mathrm{Bin}(\mathfrak{s}, \mathfrak{s}/p)$. Applying Lemma~\ref{l:dila}, we obtain
\begin{align*}
\frac{1}{{p \choose \mathfrak{s}}}\sum_{S \subseteq [p],\, \vert S\vert = \mathfrak{s}} \l[\cosh\l( \frac{9\rho  }{2^{\vert m-\ell\vert/2+1}\mathfrak{s}}\r)\r]^{\vert S \cap [\mathfrak s]\vert/2} 
& \le  \E \l[\cosh\l(  \frac{9\rho  }{2^{\vert m-\ell\vert/2+1}\mathfrak{s}}\r)^{\xi /2}\r]\\
& = \l\{ 1 + \frac{\mathfrak{s}}{p}\l( \cosh\l( \frac{9\rho  }{2^{\vert m-\ell\vert/2+1}\mathfrak{s}}\r)^{1/2} - 1 \r)\r\}^{\mathfrak{s}}.
\end{align*}
Let $\bar{m}$ and $\bar{\ell}$ be independent uniform random variables on $\mc M$. Then, 
\begin{align*}
\chi^2(\p^{(1)},\p^{(0)}) + 1 \le  \sqrt{2}\E\l[\exp\l(\frac{\rho^2}{2^{\bar{m}\vee \bar{\ell}}}\r)\l\{ 1 + \frac{\mathfrak{s}}{p}\l( \cosh\l( \frac{9\rho  }{2^{\vert \bar{m}-\bar{\ell}\vert/2+1}\mathfrak{s}}\r)^{1/2} - 1 \r)\r\}^{\mathfrak{s}}\r].
\end{align*}
Note that, for $m,\ell\in\mc M$, $\sqrt{2}\exp\bigl(\rho^2/2^{m\vee\ell}\bigr) \le \sqrt{2} \exp\bigl(2\underline{c}_1^2\bigr) \le 3/2.$ Next we consider separately two scenarios:
\begin{description}
\item[{Dense regime:}] $s \ge \sqrt{p\log\log (n)}$. 
It implies that $\rho \le c_1\sqrt{p\log\log(n)}$. By the basic inequalities $\cosh(x)\le \exp(x^2/2)$, $\exp(x) - 1 \le x\exp(x)$, and $1 +x\le \exp(x)$ for $x\ge 0$, we have 
\begin{align*}
\E\l[\l\{ 1 + \frac{\mathfrak{s}}{p}\l( \cosh\l( \frac{9\rho  }{2^{\vert \bar{m}-\bar{\ell}\vert/2+1}\mathfrak{s}}\r)^{1/2} - 1 \r)\r\}^{\mathfrak{s}}\r]
\le \E\l[\exp\l\{ \frac{81\rho^2}{p2^{\vert \bar{m}-\bar{\ell}\vert+4}}\exp\l( \frac{81\rho^2  }{2^{\vert \bar{m}-\bar{\ell}\vert+4}\mathfrak{s}^2}\r)\r\}\r].
\end{align*}
We split the expectation above into two parts corresponding to $\vert\bar{m}-\bar{\ell}\vert \le \log\log(n)$ and  $\vert\bar{m}-\bar{\ell}\vert > \log\log (n)$. Then, the first part satisfies
\begin{align*}
& \E\l[ \exp\l\{ \frac{81\rho^2}{p2^{\vert \bar{m}-\bar{\ell}\vert+4}}\exp\l( \frac{81\rho^2  }{2^{\vert \bar{m}-\bar{\ell}\vert+4}\mathfrak{s}^2}\r)\r\}\mathbb{I}_{\l\{\vert\bar{m}-\bar{\ell}\vert \le \log\log(n)\r\}}\r] \\
\le\; & \exp\l\{\frac{81\rho^2}{16 p}\exp\l(\frac{81\rho^2}{16\mathfrak{s}^2}\r)\r\} \p\l\{\vert\bar{m}-\bar{\ell}\vert \le \log\log(n)\r\} \\
\le\; &\log(n)^{\frac{81\underline{c}_1^2}{16}\exp\l(\frac{81\underline{c}_1^2}{16}\r)}\cdot\frac{2\log\log (n) +1}{(1-2\nu)\log_2(n) - 2 }
\end{align*}
which goes to zero as $n\to\infty$, as $\frac{81\underline{c}_1^2}{16}\exp\l(\frac{81\underline{c}_1^2}{16}\r) < \frac{1}{5}$. For the second part, we have 
\begin{align*}
& \E\l[ \exp\l\{\frac{81\rho^2}{p2^{\vert \bar{m}-\bar{\ell}\vert+4}}\exp\l( \frac{81\rho^2  }{2^{\vert \bar{m}-\bar{\ell}\vert+4}\mathfrak{s}^2}\r)\r\}\mathbb{I}_{\l\{\vert\bar{m}-\bar{\ell}\vert > \log\log(n)\r\}}\r] \\
\le\;& \exp\l\{\frac{81\rho^2}{p2^{\log\log(n)+4}}\exp\l(\frac{81\rho^2}{16\mathfrak{s}^2}\r)\r\} \p\l\{\vert\bar{m}-\bar{\ell}\vert > \log\log(n)\r\} \\
\le \; & \exp\l\{\frac{81\underline{c}_1^2\log\log(n)}{16\log(n)^{\log(2)}}\exp\l(\frac{81\underline{c}_1^2}{16}\r)\r\},
\end{align*}
which goes to one, as $n\to\infty$.
\item[{Sparse regime:}] $s < \sqrt{p\log\log (n)}$. By the basic inequalities $\cosh(x)\le \exp(x)$ and $1+x \le \exp(x)$ for $x\ge 0$, we have 
\begin{align*}
\E\l[\l\{ 1 + \frac{\mathfrak{s}}{p}\l( \cosh\l( \frac{9\rho  }{2^{\vert \bar{m}-\bar{\ell}\vert/2+1}\mathfrak{s}}\r)^{1/2} - 1 \r)\r\}^{\mathfrak{s}}\r]
\le \E\l[ \exp\l( \frac{\mathfrak{s}^2}{p}\l\{\exp\l( \frac{9\rho  }{2^{\vert \bar{m}-\bar{\ell}\vert/2+2}\mathfrak{s}}\r)-1\r\}\r)\r].
\end{align*}
Note that $\rho \le \underline{c}_1 \mathfrak{s}\log\bigl(1 +  p\log\log(n)/\mathfrak{s}^2\bigr)$. Similarly as in the dense regime, we first bound  
\begin{align*}
& \E\l[ \exp\l( \frac{\mathfrak{s}^2}{p}\l\{\exp\l( \frac{9\rho  }{2^{\vert \bar{m}-\bar{\ell}\vert/2+2}\mathfrak{s}}\r)-1\r\}\r)\mathbb{I}_{\l\{\vert\bar{m}-\bar{\ell}\vert \le \log\log(n)\r\}} \r]\\
\le & \exp\l( \frac{\mathfrak{s}^2}{p}\l\{\exp\l( \frac{9\rho  }{4\mathfrak{s}}\r)-1\r\}\r)\p\l\{\vert\bar{m}-\bar{\ell}\vert \le \log\log(n)\r\} \\
\le & \exp\l(\log\log(n)\cdot\frac{\mathfrak{s}^2}{p\log\log(n)}\l\{\l(1 + \frac{p\log\log(n)}{\mathfrak{s}^2}\r)^{{9\underline{c}_1}/{4}} - 1\r\}\r) \p\l\{\vert\bar{m}-\bar{\ell}\vert \le \log\log(n)\r\}\\
\le & \log(n)^{\alpha_1}\cdot\frac{2\log\log (n) +1}{(1-2\nu)\log_2(n) - 2 }\; \to\; 0, \qquad\text{as }\; n \to\infty,
\end{align*}
where $\alpha_1 = \max_{0 < x \le 1} x \bigl((1+x^{-1})^{9\underline{c}_1/4} -1\bigr) = 2^{9\underline{c}_1/4} - 1 < 0.3$. We next bound
\begin{align*}
& \E\l[ \exp\l( \frac{\mathfrak{s}^2}{p}\l\{\exp\l( \frac{9\rho  }{2^{\vert \bar{m}-\bar{\ell}\vert/2+2}\mathfrak{s}}\r)-1\r\}\r)\mathbb{I}_{\l\{\vert\bar{m}-\bar{\ell}\vert > \log\log(n)\r\}} \r]\\
\le & \exp\l(\frac{\mathfrak{s}^2}{p}\l\{\exp\l( \frac{9\rho  }{4\log(n)^{\log(2)/2}\mathfrak{s}}\r)-1\r\}\r) \\
\le & \exp\l(\log\log(n)\cdot\frac{\mathfrak{s}^2}{p\log\log(n)}\l\{\l(1 +\frac{p\log\log(n)}{\mathfrak{s}^2}\r)^{\frac{9\underline{c}_1 }{4\log(n)^{\log(2)/2}}}-1\r\}\r) \\
\le & \exp\l(\log\log(n)\cdot\l\{ 2^{\frac{9\underline{c}_1 }{4\log(n)^{\log(2)/2}}}-1\r\}\r)\;\to\; 1, \qquad\text{as }\; n\to \infty,
\end{align*}
where the last inequality follows from the fact that $\max_{0< x \le 1} x\bigl((1+x^{-1})^{\alpha_2} - 1\bigr) = 2^{\alpha_2} - 1$ for any $\alpha_2\in(0,1)$. 
\end{description}
Combining both scenarios, we can have, for sufficiently large $n$,
\[
\E\l[\l\{ 1 + \frac{\mathfrak{s}}{p}\l( \cosh\l( \frac{9\rho  }{2^{\vert \bar{m}-\bar{\ell}\vert/2+1}\mathfrak{s}}\r)^{1/2} - 1 \r)\r\}^{\mathfrak{s}}\r] \le \frac{4}{3},
\]
which implies $\chi^2(\p^{(1)},\p^{(0)}) + 1 \le \frac{3}{2} \cdot \frac{4}{3}$, that is,  $\chi^2(\p^{(1)},\p^{(0)}) \le 1$. Thus
\[
\frac{\exp\l(-\chi^2(\p^{(1)},\p^{(0)})\r)}{4} \ge \frac{1}{4e},
\]
which is the intended lower bound.  

\subsubsection{Supporting lemmas}

\begin{lem}\label[lem]{l:basic}
For $0 \le x \le 3/4$, it holds that $1/(1 -x)\le \exp(2x)$ and $1/(1+x) \le \exp(x/2)$.
\end{lem}

\begin{lem}\label[lem]{l:cross}
 For $\bm u \in \R^p$, let $\p_{\bm u}$ denote the centred $(p+1)$-dimensional Gaussian distribution with covariance matrix
\[
\bmx
    \bm\Sigma & \bm u \\
    \bm u^\top & \bm u^\top\bm u + \sigma^2
\emx \in \R^{(p+1)\times(p +1)}
\]
for some positive definite matrix $\bm\Sigma$. 
If $\lvert\bm\Sigma^{1/2}\bm u\rvert_2 \cdot \lvert\bm\Sigma^{1/2}\bm v\rvert_2 + \langle \bm\Sigma^{1/2}\bm u, \bm\Sigma^{1/2}\bm v\rangle  < \sigma^{2}$, then%\footnote{This reveals that the computation in \citet{VeVi10} is wrong, see equation (8.5) therein. } 
\begin{equation*}%\label{eq:uv}
\E_{\p_{\bm 0}}\l(\frac{\mathrm{d}\p_{\bm u}}{\mathrm{d}\p_{\bm 0}}\frac{\mathrm{d}\p_{\bm v}}{\mathrm{d}\p_{\bm 0}}\r)= \l(\l(1-\frac{\langle\bm\Sigma^{1/2}\bm u, \bm\Sigma^{1/2}\bm v\rangle}{\sigma^2}\r)^2 - \l(\frac{\lvert\bm\Sigma^{1/2}\bm u\rvert_2\cdot\lvert\bm\Sigma^{1/2}\bm v\rvert_2}{\sigma^2}\r)^2\r)^{-1/2}.
\end{equation*}
\end{lem}
\begin{proof}
It follows from a general result for centred multivariate Gaussians (see e.g.\ Lemma~11 in \citealp{CaGu17}). For completeness, we give below a direct computation: 
\begin{align*}
\E_{\p_{\bm 0}}\l(\frac{\mathrm{d}\p_{\bm u}}{\mathrm{d}\p_{\bm 0}}\frac{\mathrm{d}\p_{\bm v}}{\mathrm{d}\p_{\bm 0}}\r) & = \E_{\p_{\bm0}}\l(\frac{\mathrm{d}\p_{\bm u}^{\bm Y \vert \mbf x}\mathrm{d}\p_{\bm u}^{\mbf x}}{\mathrm{d}\p_{\bm 0}^{\bm Y \vert \mbf x}\mathrm{d}\p_{\bm 0}^{\mbf x}}\frac{\mathrm{d}\p_{\bm v}^{\bm Y \vert \mbf x}\mathrm{d}\p_{\bm v}^{\mbf x}}{\mathrm{d}\p_{\bm 0}^{\bm Y \vert \mbf x}\mathrm{d}\p_{\bm 0}^{\mbf x}}\r) =\E_{\p_{\bm0}}\l(\frac{\mathrm{d}\p_{\bm u}^{\bm Y \vert \mbf x}}{\mathrm{d}\p_{\bm 0}^{\bm Y \vert \mbf x}}\frac{\mathrm{d}\p_{\bm v}^{\bm Y \vert \mbf x}}{\mathrm{d}\p_{\bm 0}^{\bm Y \vert \mbf x}}\r) \nonumber \\
& = \mathrm{det}\l(\mbf I - \frac{1}{\sigma^2}\l(\tilde{\bm u} \tilde{\bm v}^{\top} + \tilde{\bm v}\tilde{\bm u}^{\top}\r)\r)^{-1/2} \nonumber \\
& = \l(\l(1-\frac{\langle\tilde{\bm u}, \tilde{\bm v}\rangle}{\sigma^2}\r)^2 - \l(\frac{\vert\tilde{\bm u}\vert_2\cdot\vert\tilde{\bm v}\vert_2}{\sigma^2}\r)^2\r)^{-1/2},
\end{align*}
where $\tilde{\bm u} := \bm\Sigma^{1/2}\bm u$, $\tilde{\bm v} := \bm\Sigma^{1/2}\bm v$ and the last equality can be seen from the fact that 
\begin{multline*}
\l(\mbf I - \frac{1}{\sigma^2}\l(\tilde{\bm u} \tilde{\bm v}^{\top} + \tilde{\bm v}\tilde{\bm u}^{\top}\r)\r) (\vert\tilde{\bm v}\vert_2\cdot\tilde{\bm u} \pm \vert\tilde{\bm u}\vert_2\cdot\tilde{\bm v}) \\
=  \l(1-\frac{\langle\tilde{\bm u}, \tilde{\bm v}\rangle}{\sigma^2} \mp  \frac{\vert\tilde{\bm u}\vert_2\cdot\vert\tilde{\bm v}\vert_2}{\sigma^2}\r)(\vert\tilde{\bm v}\vert_2\cdot\tilde{\bm u} \pm \vert\tilde{\bm u}\vert_2\cdot\tilde{\bm v}).\qedhere
\end{multline*}
\end{proof}
\begin{rem}
Note that in Lemma~\ref{l:cross} the condition $\vert\bm\Sigma^{1/2}\bm u\vert_2 \cdot\vert\bm\Sigma^{1/2}\bm v\vert_2 + \langle \bm\Sigma^{1/2}\bm u, \bm\Sigma^{1/2}\bm v\rangle  < \sigma^{2}$ is needed to ensure the invertibility of the matrix $\mbf I - \frac{1}{\sigma^2}\bm\Sigma^{1/2}\l(\bm u \bm v^{\top} + \bm v \bm u^{\top}\r)\bm\Sigma^{1/2}$. In the special case of $\vert\bm\Sigma^{1/2}\bm u\vert_2 \cdot\vert\bm\Sigma^{1/2}\bm v\vert_2+\langle \bm\Sigma^{1/2}\bm u,\bm\Sigma^{1/2}\bm v\rangle = \sigma^2$, by the same calculation, we have 
\[
\E_{\p_{\bm 0}}\l(\frac{\mathrm{d}\p_{\bm u}}{\mathrm{d}\p_{\bm 0}}\frac{\mathrm{d}\p_{\bm v}}{\mathrm{d}\p_{\bm 0}}\r) = \l(1-\frac{\langle\bm\Sigma^{1/2}\bm u, \bm\Sigma^{1/2}\bm v\rangle}{\sigma^2} + \frac{\vert\bm\Sigma^{1/2}\bm u\vert_2 \cdot\vert\bm\Sigma^{1/2}\bm v\vert_2}{\sigma^2}\r)^{-1/2}.
\]
However, if $\vert\bm\Sigma^{1/2}\bm u\vert_2 \cdot\vert\bm\Sigma^{1/2}\bm v\vert_2+\langle \bm\Sigma^{1/2}\bm u,\bm\Sigma^{1/2}\bm v\rangle > \sigma^2$, we have $\E_{\p_{\bm 0}}\l(\tfrac{\mathrm{d}\p_{\bm u}}{\mathrm{d}\p_{\bm 0}}\tfrac{\mathrm{d}\p_{\bm v}}{\mathrm{d}\p_{\bm 0}}\r) = \infty$.
\end{rem}

\begin{lem}[Dilatation]\label{l:dila}
If $X \sim \mathrm{Hypergeometric}(N, K, n)$ and $Y \sim \mathrm{Bin}(n, K/N)$, then $\E{f(X)}\le\E{f(Y)}$ for every convex function $f:\R \to \R$. 
\end{lem}

\begin{proof}
See either the first proof in \citet[Theorem~4 on pages 21--23]{Hoe63}, or a more explicit one in \citet[Proposition~(20.6) on pages 171--172]{Ald85}, or a systematic treatment in \citet[Section~3]{Kem73}. 
\end{proof}

\begin{rem}
Sampling with and without replacement can yield similar distributions. In particular, as $n \to \infty$ and $n/N \to 0$,  uniformly over $K \in [N]$,  
\[
\mathrm{TV} \biggl(\mathrm{Hypergeometric}(N, K, n),  \mathrm{Bin}\Bigl(n, \frac{K}{N}\Bigr)\biggr) = \frac{1}{\sqrt{2\pi e}}\frac{n}{N} + o\l(\frac{n}{N}\r)
\]
as shown by \citet[Theorem~1.6]{DiFr88}.
\end{rem}

\subsection{Proofs for the results in Section~\ref{sec:upper}}

Let $\mathbb{B}_d(r) = \{\mbf a: \, \vert \mbf a \vert_d \le 1\}$ denote the $\ell_d$-ball of radius $r$ with the dimension of $\mbf a$ determined within the context. 
Also $\mbf e_\ell, \, \ell \in [p]$ denote the canonical basis vectors of $\R^p$.

\subsubsection{Proof of Theorem~\ref{thm:mcscan}}

Recall that $\sigma_X^2 = \max_{i \in [p]} \Var(X_{it})$, where $\mbf x_t = (X_{1t}, \ldots, X_{pt})^\top$.

\begin{prop}
\label{prop:sparse}
Under the AMOC model~\eqref{eq:amoc}, suppose that there exists some $\eps_n \to \infty$ that satisfies $n^{-1} \eps_n \to 0$ such that $\cp \in (\eps_n, n - \eps_n)$ if $\bm\beta_0 \neq\bm\beta_1$, and $\cp = n$ if $\bm\beta_0 = \bm\beta_1$. 
Define \begin{align}
\bar{\mbf f}_k := \sqrt{\frac{k (n - k)}{n}} \l( \frac{n - \cp}{n - k} \mathbb{I}_{\{k \le \cp\}} + \frac{\cp}{k} \mathbb{I}_{\{k > \cp\}} \r) \bm\Sigma \bm\delta
\label{eq:f:bar}
\end{align}
and recall $\wt{\mbf Y}(k)$ from~\eqref{eq:tilde}. Then there exist some constants $\bar{C}_0, \bar{C}_1 \in (0, \infty)$ such that
\begin{align*}
\p\l\{ \vert \mbf X^\top \wt{\mbf Y}(k) - \bar{\mbf f}_k \vert_\infty \ge \bar{C}_0 \sigma_X \Psi z \r\} \le \bar{C}_1 p \exp(-z^2)
\end{align*}
for any $z \in (0, \eps_n^{1/2})$ and $k \in (\eps_n, n - \eps_n)$.
\end{prop}

\begin{proof}
Let us write the eigendecomposition of $\bm\Sigma$ by $\bm\Sigma = \sum_{j \in [p]} \mu_j \mbf v_j \mbf v_j^\top$, where $\mu_1 \ge \cdots \ge \mu_p \ge 0$ and $\mbf v_j^\top \mbf v_{j'} = \mathbb{I}_{\{j = j'\}}$. 
Then, $\bm\Sigma^{1/2} = \sum_{j \in [p]} \mu_j^{1/2} \mbf v_j \mbf v_j^\top$ and its (generalised) inverse is $(\bm\Sigma^{1/2})^\dagger = \sum_{j \in [m]} \mu_j^{-1/2} \mbf v_j \mbf v_j^\top$ with $m = \text{rank}(\bm\Sigma)$, such that
\begin{align}
\label{eq:x:circ}
\mbf x^\circ_t = (\bm\Sigma^{1/2})^\dagger \mbf x_t \sim_{\iid} \mc N_p\l( \mbf 0, (\bm\Sigma^{1/2})^\dagger \bm\Sigma (\bm\Sigma^{1/2})^\dagger \r),
\end{align}
where $(\bm\Sigma^{1/2})^\dagger \bm\Sigma (\bm\Sigma^{1/2})^\dagger$ is a diagonal matrix with $m$ ones and $p - m$ zeros on its diagonal.
By Lemma~\ref{lem:hw} and that $\vert \bm\Sigma^{1/2} \mbf e_i \vert_2 = \sqrt{\Var(X_{it})} \le \sigma_X$,
% $\Vert \bm\Sigma^{1/2}\Vert = \mu_1^{1/2} = \Vert \bm\Sigma\Vert^{1/2} $, 
we have for any $0 \le s < e \le n$ with $e - s \ge \eps_n$,
\begin{align}
& \p\l\{ \frac{1}{\sqrt{e - s}} \l\vert \sum_{t = s + 1}^e \mbf x_t \vep_t \r\vert_\infty \ge C \sigma_X \sigma z \r\}
\nn \\
=& \, \p\l\{ \frac{1}{\sqrt{e - s}} \max_{i \in [p]} \l\vert \sum_{t = s + 1}^e ( \bm\Sigma^{1/2} \mbf e_i )^\top \mbf x^\circ_t \vep_t \r\vert \ge C \sigma_X \sigma z \r\}
\le 2p \exp( - z^2),
\label{eq:xe:bound}
\end{align}
and
\begin{align}
& \p\l\{ \frac{1}{\sqrt{e - s}} \l\vert \sum_{t = s + 1}^e (\mbf x_t \mbf x_t^\top - \bm\Sigma) \bm\beta_j \r\vert_\infty \ge C \sigma_X \Psi z \r\}
\nn \\
=& \, \p\l\{ \frac{1}{\sqrt{e - s}} \max_{i \in [p]} \l\vert \sum_{t = s + 1}^e (\bm\Sigma^{1/2} \mbf e_i)^\top \l[ \mbf x^\circ_t (\mbf x^\circ_t)^\top - (\bm\Sigma^{1/2})^\dagger \bm\Sigma (\bm\Sigma^{1/2})^\dagger \r] \bm\Sigma^{1/2} \bm\beta_j \r\vert \ge C \sigma_X \Psi z \r\}
\nn \\
\le& \, 2p \exp( - z^2)
\label{eq:xx:bound}
\end{align}
for any $j \in \{0, 1\}$, where the constant $C$ is from Lemma~\ref{lem:hw}.
Let $\wh{\bm\gamma}_{a, b} = (b - a)^{-1} \sum_{t = a + 1}^b \mbf x_tY_t$ and $\bm\gamma_{a, b} = \E(\wh{\bm\gamma}_{a, b})$ for any $0 \le a < b \le n$.
For any $k \le \cp$, we can write
\begin{align}
& \p\l\{ \vert \mbf X^\top \wt{\mbf Y}(k) - \bar{\mbf f}_k \vert_\infty \ge \bar{C}_0 \sigma_X \Psi z \r\}
\nn \\
= & \, \p\l\{ \sqrt{\frac{k(n - k)}{n}} \Biggl\vert \frac{1}{n - k} \sum_{t = k + 1}^n \mbf x_t \mbf x_t^\top \bm\beta_1  - \frac{1}{n - k} \sum_{t = k + 1}^{\cp} \mbf x_t \mbf x_t^\top \bm\delta - \frac{1}{k} \sum_{t = 1}^k \mbf x_t \mbf x_t^\top \bm\beta_0 - \frac{n - \cp}{n - k} \bm\Sigma \bm\delta \r.
\nn \\
& \qquad \qquad \qquad \qquad \qquad \l. - \frac{1}{k} \sum_{t = 1}^k \mbf x_t \vep_t + \frac{1}{n - k} \sum_{t = k + 1}^n \mbf x_t \vep_t \Biggr\vert_\infty  \ge \bar{C}_0 \sigma_X \Psi z \r\}
\nn \\
\le & \, \p\l\{ \frac{1}{\sqrt{n - k}} \l\vert \sum_{t = k + 1}^n (\mbf x_t \mbf x_t^\top - \bm\Sigma) \bm\beta_1 \r\vert_\infty \ge \frac{\bar{C}_0}{5} \sigma_X \Psi z \r\}
\nn \\
& \, + \p\l\{ \frac{1}{\sqrt{n - k}} \l\vert \sum_{t = k + 1}^\cp (\mbf x_t \mbf x_t^\top - \bm\Sigma) \bm\delta \r\vert_\infty \ge \frac{\bar{C}_0}{5} \sigma_X \Psi z \r\}
\nn \\
& \, + \p\l\{ \frac{1}{\sqrt{k}} \l\vert \sum_{t = 1}^k (\mbf x_t \mbf x_t^\top - \bm\Sigma) \bm\beta_0 \r\vert_\infty \ge \frac{\bar{C}_0}{5} \sigma_X \Psi z \r\}
\nn \\
& \, + \p\l\{ \frac{1}{\sqrt{k}} \l\vert \sum_{t = 1}^k \mbf x_t \vep_t^\top \r\vert_\infty \ge \frac{\bar{C}_0}{5} \sigma_X \sigma z \r\} 
\nn \\
& \, + \p\l\{ \frac{1}{\sqrt{n - k}} \l\vert \sum_{t = k + 1}^n \mbf x_t \vep_t^\top \r\vert_\infty \ge \frac{\bar{C}_0}{5} \sigma_X \sigma z \r\}.
\label{eq:mcscan:decomp}
\end{align}
Note that by Lemma~\ref{lem:hw},
\begin{align*}
& \p\l\{ \frac{1}{\sqrt{n - k}} \l\vert \sum_{t = k + 1}^\cp (\mbf x_t \mbf x_t^\top - \bm\Sigma) \bm\delta \r\vert_\infty \ge C \sigma_X \Psi  z \r\} \\
= \; &  \p\l\{ \frac{1}{\sqrt{\cp - k}} \l\vert \sum_{t = k + 1}^\cp (\mbf x_t \mbf x_t^\top - \bm\Sigma) \bm\delta \r\vert_\infty \ge C \sigma_X \Psi z  \frac{\sqrt{n - k}}{\sqrt{\cp - k}} \r\} \\
\le \; & 2p \exp\l[-\min\l( z^2\frac{n-k}{\cp-k}, z\sqrt{n-k} \r) \r]  \le 2p\exp(-z^2),
\end{align*}
where the last inequality is due to $\tfrac{n-k}{\cp-k} \ge 1$ and $z \le \sqrt{\eps_n} \le \sqrt{n-k}$.
This, together with~\eqref{eq:xe:bound} and~\eqref{eq:xx:bound}, leads to the following upper bound on the RHS of~\eqref{eq:mcscan:decomp}, by $10 p \exp( - z^2)$, provided that $\bar{C}_0 \ge 5C$. The case of $k > \cp$ is handled analogously.
\end{proof}

The rest of the proof follows from the arguments analogous to those adopted in the proof of Theorem~\ref{thm:qcscan}, in combination with Proposition~\ref{prop:sparse}, Lemma~\ref{lem:f:bar} and the observation that
\begin{align*}
\p\l( \bar{T}_m \ge \bar{T}_k \r) 
\le \p\l( \vert \mbf X^\top \wt{\mbf Y}(k) - \bar{\mbf f}_k \vert_\infty + \vert \mbf X^\top \wt{\mbf Y}(m) - \bar{\mbf f}_m \vert_\infty \ge
\vert \bar{\mbf f}_k - \bar{\mbf f}_m \vert_\infty \r);
\end{align*}
we omit the detailed arguments.

\subsubsection{Proof of Proposition~\ref{prop:dense}}

Let us write $\mbf Y = (Y_1, \ldots, Y_n)^\top$ and $\bm\vep = (\vep_1, \ldots, \vep_n)^\top$.
WLOG, we regard $\eps_n$ as an integer and focus on $k$ satisfying $k \le \cp$; the case with $k > \cp$ is handled analogously.
Let us define 
\begin{align*}
\mbf C(k) = \bmx - \frac{1}{k} \mbf I_k & \mbf O \\
\mbf O & \frac{1}{n - k} \mbf I_{n - k} \emx.
\end{align*}
Then we may write
\begin{align*}
\frac{n}{k(n - k)} T_k =& \, \mbf Y^\top \mbf C(k) \l( \mbf X\mbf X^\top - \frac{1}{n}\tr(\mbf X\mbf X^\top) \mbf I_n \r) \mbf C(k) \mbf Y
\\
= & \, \mbf Y^\top \mbf C(k) \l( \tr(\bm\Sigma) \mbf I_n - \frac{1}{n}\tr(\mbf X\mbf X^\top) \mbf I_n \r) \mbf C(k) \mbf Y
\\
& +  \E(\mbf Y \vert \mbf X)^\top \mbf C(k) \l( \mbf X\mbf X^\top - \tr(\bm\Sigma) \mbf I_n \r) \mbf C(k) \E(\mbf Y \vert \mbf X)
\\
& +  \bm\vep^\top \mbf C(k) \l( \mbf X\mbf X^\top - \tr(\bm\Sigma) \mbf I_n \r) \mbf C(k) \bm\vep
\\
& +  2 \E(\mbf Y \vert \mbf X)^\top \mbf C(k) \l( \mbf X\mbf X^\top - \tr(\bm\Sigma) \mbf I_n \r) \mbf C(k) \bm\vep 
\\
=: & \, U_1(k) + U_2(k) + U_3(k) + U_4(k).
\end{align*}

\paragraph{Bound on $U_1(k)$.} With some $\mbf Z \sim \mc N_n(\mbf 0, \mbf I_n)$, we have $\mbf Y^\top \mbf C(k)^2 \mbf Y = \mbf Z^\top \wt{\mbf C}(k)^2 \mbf Z$ where
\begin{align*}
\wt{\mbf C}(k) = \text{diag}\l( \sqrt{ \vert \bm\Sigma^{1/2} \bm\beta_0 \vert_2^2 + \sigma^2 } \, \mbf I_\cp, \, \sqrt{ \vert \bm\Sigma^{1/2} \bm\beta_1 \vert_2^2 + \sigma^2} \, \mbf I_{n - \cp} \r) \mbf C(k).
\end{align*}
Noting that
% \footnote{
% HC:
% \begin{align*}
% \tr(\wt{\mbf C}(k)^2) &\le \frac{\vert \bm\Sigma^{1/2} \bm\beta_0 \vert_2^2 + \sigma^2}{k} + \frac{(\vert \bm\Sigma^{1/2} \bm\beta_0 \vert_2^2 + \sigma^2)(\cp - k)}{(n - k)^2} + \frac{(\vert \bm\Sigma^{1/2} \bm\beta_1 \vert_2^2 + \sigma^2)(n - \cp)}{(n - k)^2}
% \\
% &\le (\vert \bm\Sigma^{1/2} \bm\beta_0 \vert_2^2 + \vert \bm\Sigma^{1/2} \bm\beta_1 \vert_2^2 + \sigma^2) \l( \frac{1}{k} + \frac{1}{n - k} \r) \le (\vert \bm\Sigma^{1/2} \bm\beta_0 \vert_2^2 + \vert \bm\Sigma^{1/2} \bm\beta_1 \vert_2^2 + \sigma^2) \frac{n}{k(n - k)},
% \\
% \Vert \wt{\mbf C}(k)^2 \Vert &\le \max\l( \frac{\vert \bm\Sigma^{1/2} \bm\beta_0 \vert_2^2 + \sigma^2}{k^2}, \frac{\vert \bm\Sigma^{1/2} \bm\beta_0 \vert_2^2 + \sigma^2}{(n - k)^2}, \frac{\vert \bm\Sigma^{1/2} \bm\beta_1 \vert_2^2 + \sigma^2}{(n - k)^2} \r)
% \\
% &\le (\vert \bm\Sigma^{1/2} \bm\beta_0 \vert_2^2 + \vert \bm\Sigma^{1/2} \bm\beta_1 \vert_2^2 + \sigma^2) \frac{1}{\min(k, n - k)^2},
% \\
% \Vert \wt{\mbf C}(k)^2 \Vert_F^2 &\le \frac{(\vert \bm\Sigma^{1/2} \bm\beta_0 \vert_2^2 + \sigma^2)^2}{k^3} + \frac{(\vert \bm\Sigma^{1/2} \bm\beta_0 \vert_2^2 + \sigma^2)^2(\cp - k)}{(n - k)^4} + \frac{(\vert \bm\Sigma^{1/2} \bm\beta_1 \vert_2^2 + \sigma^2)^2(n - \cp)}{(n - k)^4}
% \\
% &\le (\vert \bm\Sigma^{1/2} \bm\beta_0 \vert_2^2 + \vert \bm\Sigma^{1/2} \bm\beta_1 \vert_2^2 + \sigma^2)^2 \l( \frac{1}{k^3} + \frac{1}{(n - k)^3} \r) \le (\vert \bm\Sigma^{1/2} \bm\beta_0 \vert_2^2 + \vert \bm\Sigma^{1/2} \bm\beta_1 \vert_2^2 + \sigma^2)^2 \frac{2}{\min(k, n - k)^3},
% \end{align*}
% }
\begin{align*}
\tr(\wt{\mbf C}(k)^2) &\le (\vert \bm\Sigma^{1/2} \bm\beta_0 \vert_2^2 + \vert \bm\Sigma^{1/2} \bm\beta_1 \vert_2^2 + \sigma^2) \frac{2}{\min(k, n - k)},
\\
\Vert \wt{\mbf C}(k)^2 \Vert &\le (\vert \bm\Sigma^{1/2} \bm\beta_0 \vert_2^2 + \vert \bm\Sigma^{1/2} \bm\beta_1 \vert_2^2 + \sigma^2) \frac{1}{\min(k, n - k)^2},
\\
\Vert \wt{\mbf C}(k)^2 \Vert_F &\le (\vert \bm\Sigma^{1/2} \bm\beta_0 \vert_2^2 + \vert \bm\Sigma^{1/2} \bm\beta_1 \vert_2^2 + \sigma^2) \frac{\sqrt{2}}{\min(k, n - k)^{3/2}},
\end{align*}
we have from Lemma~\ref{lem:chisq},
\begin{align*}
\p\l[  \mbf Z^\top \wt{\mbf C}(k)^2 \mbf Z  \ge \frac{7 (\vert \bm\Sigma^{1/2} \bm\beta_0 \vert_2^2 + \vert \bm\Sigma^{1/2} \bm\beta_1 \vert_2^2 + \sigma^2)}{\min(k, n - k)}
% \l( \frac{2}{\min(k, n - k)} + \frac{2 \sqrt{2n}}{\min(k, n - k)^{3/2}} + \frac{2n}{\min(k, n - k)^2} \r)
\r] \le e^{-\eps_n}.
\end{align*}
This, together with Lemma~\ref{lem:tr:xx}, gives
\begin{align*}
\p\l[ \vert U_1(k) \vert \ge 56 \Vert \bm\Sigma \Vert (\vert \bm\Sigma^{1/2} \bm\beta_0 \vert_2^2 + \vert \bm\Sigma^{1/2} \bm\beta_1 \vert_2^2 + \sigma^2) \frac{\sqrt{pz}}{\sqrt{n} \min(k, n - k)} \r] \le 3e^{-z}
\end{align*}
for any $z \in (0, \eps_n)$.

\paragraph{Bound on $U_2(k)$.}
For any $0 \le a < b \le n$ and $0 \le c < d \le n$ such that  $(a, b] \cap (c, d] = \emptyset$, we define $\mbf A^{a, b} := \mbf X_{a, b}^\top (\mbf X_{a, b} \mbf X_{a, b}^\top - \tr(\bm\Sigma) \mbf I_{b - a}) \mbf X_{a, b}$ and $\mbf B_{a, b}^{c, d} := \mbf X_{a, b}^\top \mbf X_{a, b} \mbf X_{c, d}^\top \mbf X_{c, d}$.
Then, 
\begin{align*}
U_2(k) =& \, \bmx - \frac{1}{k} \mbf X_{0, k} \bm\beta_0 \\ \frac{1}{n - k} \mbf X_{k, \cp} \bm\beta_0 \\ \frac{1}{n - k} \mbf X_{\cp, n} \bm\beta_1 \emx^\top \l( \mbf X\mbf X^\top - \tr(\bm\Sigma) \mbf I_n \r) \bmx - \frac{1}{k} \mbf X_{0, k} \bm\beta_0 \\ \frac{1}{n - k} \mbf X_{k, \cp} \bm\beta_0 \\ \frac{1}{n - k} \mbf X_{\cp, n} \bm\beta_1 \emx
\\
=& \, \frac{1}{k^2} \bm\beta_0^\top \mbf A^{0, k} \bm\beta_0 + \frac{1}{(n - k)^2} \bm\beta_0^\top \mbf A^{k, \cp} \bm\beta_0 + \frac{1}{(n - k)^2} \bm\beta_1^\top \mbf A^{\cp, n} \bm\beta_1
\\
& - \frac{2}{k (n - k)} \bm\beta_0^\top \mbf B_{0, k}^{k, \cp} \bm\beta_0
- \frac{2}{k (n - k)} \bm\beta_0^\top \mbf B_{0, k}^{\cp, n} \bm\beta_1
+ \frac{2}{(n - k)^2} \bm\beta_0^\top \mbf B_{k, \cp}^{\cp, n} \bm\beta_1
\\
=:& \, \sum_{\ell = 1}^6 U_{2, \ell}(k).
\end{align*}
By Lemma~\ref{lem:gauss:chaos}, for any $0 \le a < b \le n$ with $b - a > \eps_n$ and $\bm\alpha, \bm\beta \in \R^p$, we attain
\begin{align*}
\p\l( \l\vert \frac{1}{(b - a)^2} \bm\alpha^\top \mbf A^{a, b} \bm\beta - \l(1 + \frac{1}{b - a} \r) \bm\alpha^\top \bm\Sigma^2 \bm\beta \r\vert \ge C \Vert \bm\Sigma \Vert \vert \bm\Sigma^{1/2} \bm\alpha \vert_2 \vert \bm\Sigma^{1/2} \bm\beta \vert_2 \frac{\sqrt{pz}}{b - a} \r) \le 2 e^{-z}
\end{align*}
for $z \in (0, \eps_n^{1/3})$, with the constant $C$ defined in the lemma.
Also, by Lemma~\ref{lem:gauss:chaos:two}, for any $0 \le a < b \le n$ and $0 \le c < d \le n$ with $(a, b] \cap (c, d] = \emptyset$ and $\max(b - a, d - c) > \eps_n$, 
\begin{align*}
& \p\l( \l\vert \frac{1}{(b - a)(d - c)} \bm\alpha^\top \mbf B_{a, b}^{c, d} \bm\beta - \bm\alpha^\top \bm\Sigma^2 \bm\beta \r\vert \r. 
\\
& \qquad \l. \ge C \Vert \bm\Sigma \Vert \vert \bm\Sigma^{1/2} \bm\alpha \vert_2 \vert \bm\Sigma^{1/2} \bm\beta \vert_2 \sqrt{\frac{pz}{(b - a)(d - c)}} \r) \le 2 e^{-z}.
\end{align*}
As an illustration, consider $U_{2, 5}(k)$,
\begin{align*}
\p\l( \l\vert U_{2, 5}(k) - \frac{2 (n - \cp)}{n - k} \bm\beta_0^\top \bm\Sigma^2 \bm\beta_1 \r\vert \ge {2}C \Vert \bm\Sigma \Vert \vert \bm\Sigma^{1/2} \bm\beta_0 \vert_2 \vert \bm\Sigma^{1/2} \bm\beta_1 \vert_2 \sqrt{\frac{(n - \cp) pz}{k (n - k)^2}} \r) \le 2e^{-z},
\end{align*}
and we can similarly controlling each summand involved in $U_2(k)$ apart from $U_{2, 2}(k)$ when $\cp - \eps_n < k < \cp$. 
For this case, notice that by Lemma~\ref{lem:gauss:chaos},
\begin{align*}
\p\l( \l\vert U_{2, 2}(k) - \frac{(\cp - k)^2}{(n - k)^2} \l(1 + \frac{1}{\cp - k}\r) \bm\beta_0^\top \bm\Sigma^2 \bm\beta_0 \r\vert \ge  C \Vert \bm\Sigma \Vert \vert \bm\Sigma^{1/2} \bm\beta_0 \vert_2^2 \frac{(\cp - k) \sqrt{pt}}{(n - k)^2} \r) \\
\le 2\exp\l( -\min\l\{t, [(\cp - k)t]^{1/4} \r\} \r).
% \le 2\exp\l( - [(\cp - k)t]^{1/4} \r)
\end{align*}
If $t \le (\cp - k)^{1/3}$, this leads to
\begin{align}
\p\l( \l\vert U_{2, 2}(k) - \frac{(\cp - k)^2}{(n - k)^2} \l(1 + \frac{1}{\cp - k}\r) \bm\beta_0^\top \bm\Sigma^2 \bm\beta_0 \r\vert \ge  C \Vert \bm\Sigma \Vert \vert \bm\Sigma^{1/2} \bm\beta_0 \vert_2^2 \frac{ \sqrt{pz}}{n - k} \r) \le 2e^{-z},
\label{eq:u22}
\end{align}
with $z = t$.
On the other hand, for $t > (\cp - k)^{1/3}$, we set $z = [(\cp - k)t]^{1/4}$ to fulfil $z \in (0, \eps_n^{1/3})$, which leads to
\begin{align*}
\p\l( \l\vert U_{2, 2}(k) - \frac{(\cp - k)^2}{(n - k)^2} \l(1 + \frac{1}{\cp - k}\r) \bm\beta_0^\top \bm\Sigma^2 \bm\beta_0 \r\vert \ge  C \Vert \bm\Sigma \Vert \vert \bm\Sigma^{1/2} \bm\beta_0 \vert_2^2 \frac{\sqrt{pz}}{n-k}\frac{\sqrt{(\cp - k)z^3}}{n - k} \r) \le 2e^{-z}.
\end{align*}
Then, since  
\[
\frac{\sqrt{(\cp - k)z^3}}{n - k} \le \frac{\sqrt{(\cp - k)\eps_n}}{n - k} \le 1,
\]
this implies that~\eqref{eq:u22} holds.
In summary, we derive that for any $z \in (0, \eps_n^{1/3})$,
\begin{align*}
& \p\l[ \l\vert U_2(k) - f^\circ_k \r\vert \ge 
9 C \Vert \bm\Sigma \Vert (\vert \bm\Sigma^{1/2} \bm\beta_0 \vert_2 \vee \vert \bm\Sigma^{1/2} \bm\beta_0 \vert_2)^2 \frac{\sqrt{pz}}{\min(k, n - k)} \r] \le 12 e^{-z},
\end{align*}
where
\begin{align*}
f^\circ_k := \frac{(n - \cp)^2}{(n - k)^2} \bm\delta^\top \bm\Sigma^2 \bm\delta - \frac{\cp - k}{(n - k)^2} \bm\delta^\top \bm\Sigma^2 (\bm\beta_0 + \bm\beta_1) + \frac{1}{k} \bm\beta_0^\top \bm\Sigma^2 \bm\beta_0 + \frac{1}{n - k} \bm\beta_1^\top \bm\Sigma^2 \bm\beta_1
\end{align*}
for $k \le \cp$; the case where $k > \cp$ follows analogously with
\begin{align*}
f^\circ_k = \frac{\cp^2}{k^2} \bm\delta^\top \bm\Sigma^2 \bm\delta + \frac{k - \cp}{k^2} \bm\delta^\top \bm\Sigma^2 (\bm\beta_0 + \bm\beta_1) + \frac{1}{k} \bm\beta_0^\top \bm\Sigma^2 \bm\beta_0 + \frac{1}{n - k} \bm\beta_1^\top \bm\Sigma^2 \bm\beta_1.
\end{align*}

\paragraph{Bound on $U_3(k)$.}
There exist $\mbf z_t \sim_{\iid} \mc N_p(\mbf 0, \mbf I), \, t \in [n]$, such that $\mbf x_t = \bm\Sigma^{1/2} \mbf z_t$.
Then, we have, for $\mbf Z = [\mbf z_1, \ldots, \mbf z_n]^\top$, 
\begin{align*}
\bm\vep^\top \mbf C(k) \mbf X\mbf X^\top \mbf C(k) \bm\vep = \bm\vep^\top \mbf C(k) \mbf Z \bm\Sigma \mbf Z^\top \mbf C(k) \bm\vep.
\end{align*}
Thus, conditional on $\bm\vep$, we have $\mbf Z^\top \mbf C(k) \bm\vep/\vert\mbf C(k) \bm\vep\vert_2\sim\mc N_n(\mbf 0, \mbf I_n)$ and  further by Lemma~\ref{lem:chisq}, 
\begin{align*}
\p\l( \frac{\vert U_3 \vert}{\vert \mbf C(k) \bm\vep \vert_2^2} \ge 4 \Vert \bm\Sigma \Vert_F \sqrt{z} \r)
= \p\l[ \frac{\vert \bm\vep^\top \mbf C(k) (\mbf Z \bm\Sigma \mbf Z^\top - \tr(\bm\Sigma) \mbf I_n) \mbf C(k) \bm\vep \vert}{\vert \mbf C(k) \bm\vep \vert_2^2} \ge 4 \Vert \bm\Sigma \Vert_F \sqrt{z} \r]
\le 2e^{-z}
\end{align*}
for any $z \in (0, p)$.
By the same lemma, we have for any $z > 0$,
\begin{align*}
\p\l( \frac{\vert \mbf C(k) \bm\vep \vert_2^2}{\sigma^2} \ge \tr(\mbf C(k)^2) + 2 \Vert \mbf C(k)^2 \Vert_F \sqrt{z} + 2 \Vert \mbf C(k)^2 \Vert z \r) \le e^{-z}.
\end{align*}
Further limiting to $z \in (0, \eps_n)$, it holds
\begin{align*}
\tr(\mbf C(k)^2) &= \frac{n}{k(n - k)}, \quad \text{and}
\\
\Vert \mbf C(k)^2 \Vert_F \sqrt{z} &= \sqrt{\frac{z}{k^3} + \frac{z}{(n - k)^3}} \le \sqrt{\frac{1}{k^2} + \frac{1}{(n - k)^2}} \le \frac{n}{k(n - k)}.
\end{align*}
Combining the above observations, we derive that, for any $z \in (0, \eps_n)$,
\begin{align*}
\p\l[ \vert U_3 \vert \ge 20 \Vert \bm\Sigma \Vert \sigma^2 \frac{n \sqrt{pz}}{k(n - k)} \r] \le 3e^{-z}.
\end{align*}

\paragraph{Bound on $U_4(k)$.} Let us write
\begin{align*}
\frac{1}{2} U_4(k) = & \, 
\frac{1}{k^2} \bm\beta_0^\top \mbf X_{0, k}^\top (\mbf X_{0, k} \mbf X_{0, k}^\top - \tr(\bm\Sigma) \mbf I_k) \bm\vep_{0, k} 
+ \frac{1}{(n - k)^2} \bm\beta_0^\top \mbf X_{k, \cp}^\top (\mbf X_{k, \cp} \mbf X_{k, \cp}^\top - \tr(\bm\Sigma) \mbf I_{\cp - k}) \bm\vep_{k, \cp}
\\
& + \frac{1}{(n - k)^2} \bm\beta_1^\top \mbf X_{\cp, n}^\top (\mbf X_{\cp, n} \mbf X_{\cp, n}^\top - \tr(\bm\Sigma) \mbf I_{n - \cp}) \bm\vep_{\cp, n}
- \frac{1}{k(n - k)} \bm\beta_0^\top \mbf X_{0, k}^\top \mbf X_{0, k} \mbf X_{k, \cp}^\top \bm\vep_{k, \cp} 
\\
& - \frac{1}{k(n - k)} \bm\beta_0^\top \mbf X_{0, k}^\top \mbf X_{0, k} \mbf X_{\cp, n}^\top \bm\vep_{\cp, n} 
- \frac{1}{k(n - k)} \bm\beta_0^\top \mbf X_{k, \cp}^\top \mbf X_{k, \cp} \mbf X_{0, k}^\top \bm\vep_{0, k} 
\\
& + \frac{1}{(n - k)^2} \bm\beta_0^\top \mbf X_{k, \cp}^\top \mbf X_{k, \cp} \mbf X_{\cp, n}^\top \bm\vep_{\cp, n}
+ \frac{1}{k(n - k)} \bm\beta_1^\top \mbf X_{\cp, n}^\top \mbf X_{\cp, n} \mbf X_{0, k}^\top \bm\vep_{0, k}
\\
& 
+ \frac{1}{(n - k)^2} \bm\beta_1^\top \mbf X_{\cp, n}^\top \mbf X_{\cp, n} \mbf X_{k, \cp}^\top \bm\vep_{k, \cp}
=: \sum_{\ell = 1}^9 U_{4, \ell}(k).
\end{align*}
Consider $U_{4, 1}(k)$ which, conditional on $\mbf x_t, \, t \in [n]$, follows a centered normal distribution with
\begin{align*}
\Var\l( \l. \frac{k^2 U_{4, 1}(k)}{\sigma} \r\vert \mbf X \r) \le 
\vert \mbf X_{0, k} \bm\beta_0 \vert_2^2 \Vert \mbf X_{0, k}\mbf X_{0, k}^\top - \tr(\bm\Sigma) \mbf I_k \Vert^2.
\end{align*}
Notice that $\vert \mbf X_{0, k} \bm\beta_0 \vert_2^2$ follows a scaled $\chi_k^2$-distribution such that by Lemma~\ref{lem:chisq},
\begin{align*}
\p\l( \frac{\vert \mbf X_{0, k} \bm\beta_0 \vert_2^2}{\vert \bm\Sigma^{1/2} \bm\beta_0 \vert_2^2} \ge 5 k \r) \le e^{-\eps_n}.
\end{align*}
Combining this observation with Lemma~\ref{lem:op:xx} and Mill's ratio, we obtain
\begin{align*}
\p\l( \vert U_{4, 1}(k) \vert \ge 
50\sqrt{10} \Vert \bm\Sigma \Vert \sigma \vert \bm\Sigma^{1/2} \bm\beta_0 \vert_2 \frac{\sqrt{p z}}{k} \r) \le 5 e^{-z}
\end{align*}
for any $z \in (0, \eps_n)$.
Next, focusing on $U_{4, 6}(k)$, its conditional distribution is also a centered normal distribution with
\begin{align*}
\Var\l( \l. \frac{k(n - k) U_{4, 6}(k)}{\sigma} \r\vert \mbf X \r) \le 
\vert \mbf X_{k, \cp} \bm\beta_0 \vert_2^2 \Vert \mbf X_{0, k}\mbf X_{k, \cp}^\top \Vert^2.
\end{align*}
As in $U_{4, 1}(k)$, Lemma~\ref{lem:chisq} gives
\begin{align*}
\p\l( \frac{\vert \mbf X_{k, \cp} \bm\beta_0 \vert_2^2}{\vert \bm\Sigma^{1/2} \bm\beta_0 \vert_2^2} \ge \cp - k + 2\sqrt{(\cp - k)(n-k)} + 2 (n - k) \r) \le e^{-\eps_n}.
\end{align*}
Then, as $\mbf X_{0, k}\mbf X_{k, \cp}^\top$ is a sub-matrix of $\mbf X_{0, \cp}\mbf X_{0, \cp}^\top - \tr(\bm\Sigma) \mbf I_\cp$, Lemma~\ref{lem:op:xx} and Mill's ratio imply 
\begin{align*}
\p\l( \vert U_{4, 6}(k) \vert \ge 
50\sqrt{10} \Vert \bm\Sigma \Vert \sigma \vert \bm\Sigma^{1/2} \bm\beta_0 \vert_2 \frac{\sqrt{\cp pz}}{k\sqrt{n - k}} \r) \le 5 e^{-z}
\end{align*}
for any $z \in (0, \eps_n)$.
Analogous arguments apply to the rest of $U_{4, \ell}(k)$'s, from which we obtain
% \footnote{HC:
% \begin{align*}
% & \frac{1}{k} + \frac{\sqrt{(n - k)(\cp - k)}}{(n - k)^2} + \frac{n - \cp}{(n - k)^2} + \frac{\sqrt{k \cp}}{k(n - k)} + \frac{\sqrt{k(n - \cp + k)}}{k(n - k)} 
% \\
% & + \frac{\sqrt{(n - k) \cp}}{k(n - k)} + \frac{\sqrt{(n - k)(n - k)}}{(n - k)^2} + \frac{\sqrt{(k + n - \cp)(n - \cp)}}{k(n - k)} + \frac{\sqrt{(n - k)(n - \cp)}}{(n - k)^2}
% \\
% \le & \, \frac{1}{k} + \frac{4}{n - k} + \frac{4\sqrt{2}}{\min(k, n - k)},
% \end{align*}
% all from noticing that e.g.\ if $k < n - k$, then $\cp < k + (\cp - k) < 2(n - k)$ while otherwise, $\cp < 2k$
% }
\begin{align*}
\p\l[ \vert U_4(k) \vert \ge 450 \sqrt{10} \Vert \bm\Sigma \Vert \sigma (\vert \bm\Sigma^{1/2} \bm\beta_0 \vert_2 \vee \vert \bm\Sigma^{1/2} \bm\beta_1 \vert_2) \frac{n\sqrt{pz}}{k(n - k)} \r] \le 45e^{-z}.
\end{align*}

Collecting the bounds on $U_1(k), \ldots, U_4(k)$, we have for $k \le\cp$,
\begin{align*}
& \p\l\{ \l\vert T_k - f_k \r\vert \ge C_0 \Vert \bm\Sigma \Vert \Psi^2 \sqrt{pz} \r\} \le C_1 e^{-z}
\end{align*}
for $z \in (0, \eps_n^{1/3})$ and some constants $C_0, C_1 \in (0, \infty)$, which completes the proof.

\subsubsection{Proof of Theorem~\ref{thm:qcscan}}

The proof follows closely the steps in the proofs of Theorems~1 and~2 of \cite{kovacs2024optimistic}.
For notational convenience, we assume that there is no rounding in determining the dyadic search locations and the probe points.
Throughout, we omit the subscript `Q' where there is no confusion.

We verify the claims in the theorem for
\begin{align}
c = \sqrt{2} C_0 + 1, \, 
c_0 \ge \max(16c, 128 C_0), \,
c_1 \ge 8\sqrt{3} C_0 \text{ \ and \ }
c_2 \ge 96 C_1.
% \l( \frac{12}{\log(2)} + 65 \r) C_1.
\label{eq:const}
\end{align}
Note that under~\eqref{eq:dense:dlb}, it holds that
\begin{align}
\frac{\Psi^2}{\vert \bm\Sigma \bm\delta \vert_2^2 \min(\cp, n - \cp)} \le \frac{1}{c_0\sqrt{p \log\log(n)}} \le \frac{1}{2}
\label{eq:delta:beta}
\end{align}
provided that $c_0 \sqrt{p\log\log(n)} \ge 2$.

In Step~0, we address the case of no change point ($q = 0$) while in Steps~1--2, we consider the case where $\cp < n$, i.e.\ there exists a single change point ($q = 1$).
WLOG, we suppose that $\cp \le n / 2$. 
By design, any probe point at which $T_k$ is evaluated in Algorithm~\ref{alg:os} satisfies $\min(k, n - k) > \varpi$.
% the set of probe points on which $T_k$ is evaluated in Algorithm~\ref{alg:os}, belong to a grid $\mc G_*$ that is completely determined by the signal and is independent of the noise.
% Moreover, it can be shown that
% \begin{align}
% \vert \mc G_* \vert \le 2 (4\log(n))^2 + 2 \log_2(n) \le 40 \log^2(n),
% \label{eq:set:o}
% \end{align}
% provided that $3 \le \log(n)$.
% This follows from the observation that $\vert \mc D \vert \le 2 \log_2(n)$ and $\vert \mc P \vert \le 16 \log^2(n)$, where $\mc P$ is defined in Step~2 below.
% Also, by design, we have $\mc G \cap \{ k: \, 0 \le k \le \varpi \text{ or } n - \varpi + 1 \le k \le n\} = \emptyset$.

\paragraph{Step 0.} 
When $q = 0$, we have from~\eqref{eq:f},
\begin{align*}
f_k = \frac{1}{k} \bm\beta_0^\top \Sigma^2 \bm\beta_0 + \frac{1}{n - k} \bm\beta_1^\top \Sigma^2 \bm\beta_1 \le 2 \Vert \bm\Sigma \Vert \Psi^2.
\end{align*}
Then, noting that $\vert \mc D \vert \le 2 \log_2(n)$ and is deterministic, an application of Proposition~\ref{prop:dense} shows that with $c = \sqrt{2} C_0 + 1$, we have 
\begin{align*}
\p\l( \max_{k \in \mc D} T_k > \zeta \r) 
&\le \p\l( \max_{k \in \mc D} \vert T_k - f_k \vert \ge \zeta - 2 \Vert \bm\Sigma \Vert \Psi^2 \r)
\\
&\le \p\l( \max_{k \in \mc D} \vert T_k - f_k \vert \ge C_0 \Vert \bm\Sigma \Vert \Psi^2 \sqrt{2 p\log\log(n)} \r)
\\
&\le 2 \log_2(n) \cdot C_1 \exp(- 2 \log\log(n)) \le \frac{2 C_1}{\log(2)} \bigl(\log(n)\bigr)^{-1},
\end{align*}
provided that $\sqrt{p\log\log(n)} \ge 2$. 
That is, with probability greater than $1 - \tfrac{2 C_1}{\log(2)} \bigl(\log(n)\bigr)^{-1}$, Algorithm~\ref{alg:os} returns $(\wh\cp, \wh q) = (n, 0)$.

\paragraph{Step 1.}
We can always find $l^* \in \{1, \ldots, \ell - 1\}$ such that $2^{-l^* - 1}n < \cp \le 2^{-l^*} n$. In this step, we establish that the probe point $k^*$ selected from $\mc D$ is such that the event $\mc A_1 := \{ k^* \in \{ 2^{-l^*  - 1}n, 2^{-l^*}n \} \text{ \ and \ } T_{k^*} > \zeta \}$ holds with probability tending to one. 

Consider some $m \in \mc D \cap (0, 2^{-l^* - 1}n)$, for which it holds that $2^{-l^* - 1}n - m \ge 2^{-l^* - 2}n \ge \cp/4$. % > 4 \varpi$.
Then, by Proposition~\ref{prop:dense} and Lemma~\ref{lem:f}, we have for any $m \in \mc D \cap (0, 2^{-l^* - 1}n)$,
\begin{align}
& \p\l( T_{2^{-l^* - 1}n} < T_m \r)
\nn \\
\le & \, \p\l( \l\vert T_{2^{-l^* - 1}n} - T_m - (f_{2^{-l^* - 1}n} - f_m) \r\vert \ge f_{2^{-l^* - 1}n} - f_m \r)
\nn \\
\le & \, 2 C_1 \exp\l\{ - \l[ \frac{ f_{2^{-l^* - 1}n} - f_m }{2 C_0 \Vert \bm\Sigma \Vert \Psi^2 \sqrt{p}} \r]^2 \r\}
\nn \\
\le & \, 2 C_1 \exp\l\{ - \l[ 
\frac{\vert \bm\Sigma \bm\delta \vert_2^2 (2^{-l^*-1}n - m) (n - \cp)^2}{2 C_0 (n - 2^{-l^*-1}n)(n - m) \Vert \bm\Sigma \Vert \Psi^2 \sqrt{p}}
\r]^2 \r\} % since $n - k, n - m \le n$ and $n - \cp \ge n/2$
\nn \\
\le & \, 2 C_1 \exp\l\{ - \l[ 
\frac{\vert \bm\Sigma \bm\delta \vert_2^2 \cp}{64 C_0 \Vert \bm\Sigma \Vert \Psi^2 \sqrt{p}} \r]^2 \r\} 
\nn \\
\le & \, 2 C_1 \exp( - 2 \log\log(n)),
\label{eq:thm:qcscan:one}
\end{align}
where the last inequality follows from~\eqref{eq:dense:dlb} provided that $c_0 \ge 128 C_0$.

Next, denote by $m$ the smallest element of $\mc D \cap (2^{-l^*}n, n]$, for which we have either $m = 2^{-l^* + 1}n$ or $m = 3n/4 = 3 \cdot 2^{-l^* - 1}n$ such that $m < 3 \cp$.
Also, $m - 2^{-l^*}n \ge \min(\cp/2, (n - \cp)/4)$.
Then, by analogous arguments leading to~\eqref{eq:thm:qcscan:one}, for any $m \in \mc D \cap (2^{-l^*}n, n]$,
\begin{align}
& \p\l( T_{2^{-l^*}n} < T_m \r)
\nn \\
\le & \, \p\l( \l\vert T_{2^{-l^*}n} - T_m - (f_{2^{-l^*}n} - f_m) \r\vert \ge f_{2^{-l^*}n} - f_m \r)
\nn \\
% \le & \, 2 C_1 \exp\l\{ - \l[ \frac{ f_{2^{-l^* n + \varpi} - f_m}{2 C_0 \Vert \bm\Sigma \Vert \Psi^2 \sqrt{p}} \r]^2 \r\}
%\nn \\
\le & \, 2 C_1 \exp\l\{ - \l[ 
\frac{\vert \bm\Sigma \bm\delta \vert_2^2 (m - 2^{-l^*}n) \cp^2}{2 C_0 \cdot 2^{-l^*}n m \Vert \bm\Sigma \Vert \Psi^2 \sqrt{p}}
\r]^2 \r\}
\nn \\
\le & \, 2 C_1 \exp\l\{ - \l[ 
\frac{\vert \bm\Sigma \bm\delta \vert_2^2 \Delta}{48 C_0 \Vert \bm\Sigma \Vert \Psi^2 \sqrt{p}} \r]^2 \r\} 
\nn \\
\le & \, 2 C_1 \exp( - 2 \log\log(n)),
\label{eq:thm:qcscan:two}
\end{align}
provided that $c_0 \ge 96 C_0$.
Finally, from the definition of $l^*$ and~\eqref{eq:f}, we have
\begin{align*}
\min\l( f_{2^{-l^* - 1} n}, f_{2^{-l^*} n} \r) \ge \min\l\{ \frac{2^{-l^* - 1} n (n - \cp)^2}{n (n - 2^{-l^* - 1}n)}, \, \frac{(n - 2^{-l^*} n ) \cp^2}{n \cdot 2^{-l^*}n} \r\} \vert \bm\Sigma \bm\delta \vert_2^2 \ge \frac{\Delta}{8} \vert \bm\Sigma \bm\delta \vert_2^2.
\end{align*}
From this and Proposition~\ref{prop:dense}, it holds that
\begin{align*}
\p\l( \max_{k \in \{2^{-l^* - 1} n, 2^{-l^*} n\}} T_k \le \zeta \r)
&\le \p\l[ \max_{k \in \{2^{-l^* - 1} n, 2^{-l^*} n\}} \vert T_k - f_k \vert \ge \min\l( f_{2^{-l^* - 1} n}, f_{2^{-l^*} n} \r) - \zeta \r]
\\
&\le \p\l[ \max_{k \in \{2^{-l^* - 1} n, 2^{-l^*} n\}} \vert T_k - f_k \vert \ge \frac{\Delta}{16} \vert \bm\Sigma \bm\delta \vert_2^2 \r] 
\\
&\le 2C_1 \exp\l\{ -\l[ \frac{\vert \bm\Sigma \bm\delta \vert_2^2 \Delta}{16 C_0 \Vert \bm\Sigma \Vert \Psi^2 \sqrt{p}} \r]^2 \r\}
\\
&\le 2C_1 \exp( - \log\log(n)),
\end{align*}
provided that $c_0 \ge \max(16 c, 16C_0)$.
Noting that
\begin{align*}
\mc A_1 \;\supset\; \l\{ T_{2^{-l^* - 1}n} > \max_{m \in \mc D \cap (0, 2^{-l^* - 1}n)} \bar{T}_m \r\}
\bigcap \l\{ T_{2^{-l^*}n} > \max_{m \in \mc D \cap (2^{-l^*}n, n]} \bar{T}_m \r\} 
\\
\bigcap \l\{ \max_{k \in \{2^{-l^* - 1} n, 2^{-l^*} n\}} T_k \le \zeta \r\},
\end{align*}
and that $\vert \mc D \vert \le 2 \log_2(n)$, we have from~\eqref{eq:thm:qcscan:one}--\eqref{eq:thm:qcscan:two},
\begin{align*}
\p(\mc A_1) \ge 1 - (4 \log^{-1}(2) + 2) C_1 \bigl(\log(n)\bigr)^{-1}.
\end{align*}

\paragraph{Step 2.}
Based on Step~1, provided that $\mc A_1$ holds, it is sufficient to consider the case where we have $(s, e) = (k^*/2, 2k^*)$ at the start of the \texttt{OS} algorithm.
Besides, it is enough to consider the pairs of the probe points $(t, u)$ considered by the \texttt{OS} algorithm which lie on the same side of $\cp$, since otherwise no matter which side is dropped off, no mistake occurs.
Thus we define
\begin{align*}
\mc P :=& \, \l\{(t, u) \;:\; \text{ pairs of the probe points in the same step of a call \texttt{OS}$(a, k^*, b)$} \r. 
\\
& \quad \l. \text{such that $t, u$ lie on the same side of $\cp$ with $\vert t - \cp \vert < \vert u - \cp \vert$}, \r. 
\\
& \quad \l. \text{and \ } \vert t - u \vert > c_1 \vert \bm\Sigma \bm\delta \vert_2^{-2} \Vert \bm\Sigma \Vert \Psi^2 \sqrt{p\log\log(n)} \r\}.
\end{align*}
It can be shown that $\mc P$ is a subset of a deterministic set, say $\mc Q$, which depends only on $f_k$ defined in~\eqref{eq:f}, and further that $\vert \mc Q \vert \le 16\log^2(n)$ (see page~52 of \citealp{kovacs2024optimistic}).
The claim on the rate of estimation follows if $\p(\mc A_2^c) \lesssim \bigl(\log(n)\bigr)^{-1}$, where
\begin{align*}
\mc A_2 := \l\{ T_t > T_u \quad \forall \, (t, u) \in \mc P \r\}.
\end{align*}
First, consider $(t, u) \in \mc P$ where $u < t \le \cp$.
By Proposition~\ref{prop:dense} and Lemma~\ref{lem:f}, we have
\begin{align}
& \p\l( T_t <  T_u \r)
\le \p\l( \l\vert T_t - T_u - (f_t - f_u) \r\vert \ge f_t- f_u \r)
\nn \\
\le & \, 4 C_1 \exp\l\{ - \l[ \frac{ f_t - f_u}{2 C_0 \Vert \bm\Sigma \Vert \Psi^2 \sqrt{p}} \r]^2 \r\}
\nn \\
\le & \, 4 C_1 \exp\l\{ - \l[ 
\frac{\vert \bm\Sigma \bm\delta \vert_2^2 (t - u) (n - \cp)^2}{2 C_0 (n - t)(n - u) \Vert \bm\Sigma \Vert \Psi^2 \sqrt{p}}
\r]^2 \r\} % since $n - k, n - m \le n$ and $n - \cp \ge n/2$
\nn \\
\le & \, 4 C_1 \exp\l\{ - \l( \frac{c_1 \sqrt{\log\log(n)}}{8 C_0} \r)^2 \r\} 
\le 4 C_1 \exp( - 3 \log\log(n)),
\label{eq:thm:qcscan:three}
\end{align}
provided that $c_1 \ge 8 \sqrt{3} C_0$.
Similarly, consider $(t, u) \in \mc P$ satisfying $\cp \le t < u$.
By construction, we have $t < u \le 2 \cp$.
Then,
\begin{align}
& \p\l( T_t < T_u \r)
\nn \\
% \le & \, 4 C_1 \exp\l\{ - \l[ \frac{ f_t - f_u}{2 C_0 \Vert \bm\Sigma \Vert \Psi^2 \sqrt{p}} \r]^2 \r\}
% \nn \\
\le  \,& 4 C_1 \exp\l\{ - \l[ 
\frac{\vert \bm\Sigma \bm\delta \vert_2^2 (u - t) \cp^2}{2 C_0 t u \Vert \bm\Sigma \Vert \Psi^2 \sqrt{p}}
\r]^2 \r\} % since $n - k, n - m \le n$ and $n - \cp \ge n/2$
\nn \\
\le  \,& 4 C_1 \exp\l\{ - \l( \frac{c_1 \sqrt{\log\log(n)}}{8 C_0} \r)^2 \r\} 
\le 4 C_1 \exp( - 3 \log\log(n)),
\label{eq:thm:qcscan:four}
\end{align}
as long as $c_1 \ge 8 \sqrt{3} C_0$.
Finally, noting that $\vert \mc P \vert \le 16\log^2(n)$, we have
\begin{align*}
\p(\mc A_2) \ge 1 - \p(\mc A_1 \cap \mc A_2^c) - \p(\mc A_1^c) \ge 1 - 80 C_1 \bigl(\log(n)\bigr)^{-1}.
\end{align*}
Further, provided that $c_0 \ge 2 c_1$, we have $\vert \wh\cp  - \cp \vert \le \Delta/4$ under~\eqref{eq:dense:dlb}.

\subsubsection{Proof of Proposition~\ref{prop:lope}}

Let us define, for some $\mbf a, \mbf b \in \R^p$,
\begin{align*}
\mc B_1(\mbf a) &:= \l\{ \max_{\substack{ 0 \le s < e \le n}} \frac{1}{\sqrt{\max(e - s, \bar{C}_1\log(p \vee n))}} \l\vert \sum_{t = s + 1}^e \mbf a^\top \mbf x_t \vep_t \r\vert \le \bar{C}_2 \vert \bm\Sigma^{1/2} \mbf a \vert_2 \sigma \sqrt{\log(p \vee n)} \r\},
\\
\mc B_2(\mbf a, \mbf b) &:= \Biggl\{ \max_{\substack{ 0 \le s < e \le n}} \frac{1}{\sqrt{\max(e - s, \bar{C}_1\log(p \vee n))}} \l\vert \sum_{t = s + 1}^e \mbf a^\top (\mbf x_t \mbf x_t^\top - \bm\Sigma) \mbf b \r\vert 
\\
& \qquad \qquad \qquad \qquad \qquad  \le \bar{C}_2 \vert \bm\Sigma^{1/2} \mbf a \vert_2 \vert \bm\Sigma^{1/2} \mbf b \vert_2 \sqrt{\log(p \vee n)} \Biggr\},
\\
\mc B_3(\mbf a) &:= \l\{ \frac{\underline{\sigma}}{2} \vert \mbf a \vert_2^2 - \frac{C_{\text{RE}} \log(p)}{\underline{\sigma} n} \vert \mbf a \vert_1^2 \le \frac{1}{n} \sum_{t \in [n]} \mbf a^\top \mbf x_t \mbf x_t^\top \mbf a \le  \frac{3\Vert \bm\Sigma \Vert}{2} \vert \mbf a \vert_2^2 + \frac{C_{\text{RE}} \log(p)}{\underline{\sigma} n} \vert \mbf a \vert_1^2 \r\},
\end{align*}
with some $C_{\text{RE}} \ge 2 (54 C)^2$, where $C$ is as in Lemma~\ref{lem:hw}.
By the same lemma, we have 
\begin{align*}
& \p\l[ \l( \bigcap_{\mbf a \in \{\bm\beta_0, \bm\beta_1, \mbf e_\ell, \, \ell \in [p] \}} \mc B_1(\mbf a) \r)^c \; \r]
\le (p + 2) \cdot n^2 \cdot 2 \exp\l( - \frac{\bar{C}_2^2 \log(p \vee n)}{C^2} \r),
\\
& \p\l[ \l( \bigcap_{\substack{\mbf a \in \{\bm\beta_0, \bm\beta_1 \} \\ \mbf b \in \{\mbf e_\ell, \, \ell \in [p]\}}} \mc B_2(\mbf a, \mbf b) \r)^c \; \r]
\le 2(p + 2) \cdot n^2 \cdot 2 \exp\l( - \frac{\bar{C}_2^2 \log(p \vee n)}{C^2} \r).
\end{align*}
Also, defining $\mathbb{K}(b) = \mathbb{B}_0(b) \cap \mathbb{B}_2(1)$ with some $b \ge 1$,
\begin{align*}
\p\l\{ \sup_{\mbf a \in \mathbb{K}(2b)} \frac{1}{n} \l\vert \sum_{t \in [n]} \mbf a^\top(\mbf x_t\mbf x_t^\top - \bm\Sigma) \mbf a \r\vert > \frac{\underline{\sigma}}{54} \r\} &\le 2 \exp\l[ - \l( \frac{\underline{\sigma}\sqrt{n}}{54C} \r)^2 + 2 b \log(p) \r] 
\\
& \le 2 \exp\l[ - \frac{1}{2} \l( \frac{\underline{\sigma}\sqrt{n}}{54C} \r)^2 \r],
\end{align*}
where the first inequality follows from the $\varepsilon$-net argument (see e.g.\ Lemma~F.2 of \citeauthor{basu2015regularized}, \citeyear{basu2015regularized}), and the second inequality holds with the choice
\begin{align*}
b = \l\lfloor \frac{\underline{\sigma}^2 n}{4 (54 C)^2 \log(p)} \r\rfloor,
\end{align*}
which is an integer under the condition on $n$.
This, in combination of Lemma~12 of \cite{loh2012high}, gives
\begin{align*}
\p\l[ \l( \bigcap_{\mbf a \in \R^p} \mc B_3(\mbf a) \r)^c \, \r] \le 2 \exp\l[ - \frac{1}{2} \l( \frac{\underline{\sigma}\sqrt{n}}{54C} \r)^2 \r],
\end{align*}
with the choice of $C_{\text{RE}}$.
In summary, we establish that $\p(\mc B) \ge 1 - c' (p \vee n)^{-1}$ for large enough $c', \bar{C}_2 \in (0, \infty)$, where
\begin{align}
\label{eq:set:b}
\mc B := \l[ \bigcap_{\mbf a \in \{\bm\beta_0, \bm\beta_1, \mbf e_\ell, \, \ell \in [p] \}} \mc B_1(\mbf a) \r] \bigcap \l[ \bigcap_{\substack{\mbf a, \mbf b \in \\ \{\bm\beta_0, \bm\beta_1, \mbf e_\ell, \, \ell \in [p]\}}} \mc B_2(\mbf a, \mbf b) \r]  \bigcap \l[ \bigcap_{\mbf a \in \R^p} \mc B_3(\mbf a) \r] .
\end{align}
When either (i) or (ii) holds, by the arguments employed in the proofs of Theorems~\ref{thm:mcscan} and~\ref{thm:qcscan}, we have $\vert \wh\cp - \cp \vert \le \eta \Delta$ for some $\eta \in (0, 1)$ that depends on $\bar{c}_0$ or $c_0$, with large probability.
Then, evoking Proposition~S3.6 in the supplement of \cite{cho2024detection}, we have
\begin{align*}
& \l\vert \wh{\bm\delta} - \bm\delta(\wh\cp) \r\vert_2 \lesssim \frac{\sigma_X \Psi \sqrt{\mathfrak{s}_\delta \log(p \vee n)}}{\underline{\sigma} \sqrt{\Delta}}
\text{ \ and \ }
\\
& \l\vert \wh{\bm\delta} - \bm\delta(\wh\cp) \r\vert_1 \lesssim \frac{\sigma_X \Psi \mathfrak{s}_\delta \sqrt{\log(p \vee n)}}{\underline{\sigma} \sqrt{\Delta}},
\end{align*}
where $\bm\delta(\wh\cp) = \bm\Sigma^{-1} (\bm\gamma_{\wh\cp, n} - \bm\gamma_{0, \wh\cp})$, deterministically on the event $\mc B$.

\paragraph{When $\wh\cp = \wh\cp_{\Mc}$.} By Theorem~\ref{thm:mcscan}, we have with probability greater than $1 - \bar{c}_2 (p\log(n))^{-1}$,
\begin{align*}
\l\vert \bm\delta(\wh\cp_{\Mc}) - \bm\delta \r\vert_1 &= 
\frac{\vert \wh\cp_{\Mc} - \cp \vert \vert \bm\delta \vert_1}{\cp \cdot \mathbb{I}_{\{\wh\cp_{\Mc} \le \cp\}} +  (n - \cp) \cdot \mathbb{I}_{\{\wh\cp_{\Mc} > \cp\}}}
\\
&\le \frac{\bar{c}_1 \vert \bm\Sigma \bm\delta \vert_\infty^{-2} \sigma_X^2 \Psi^2 \log(p \log(n)) \cdot \vert \bm\delta \vert_1 }{ \Delta}
\\
&\le \frac{\bar{c}_1 \sigma_X^2 \Psi^2 \log(p \log(n)) \cdot \mathfrak{s}_\delta}{\underline{\sigma} \sqrt{\Delta} \cdot \sqrt{\vert \bm\Sigma \bm\delta \vert_\infty^2 \Delta}} 
\\
&\lesssim \frac{\sigma_X \Psi \mathfrak{s}_\delta \sqrt{\log(p \log(n))}}{\underline{\sigma} \sqrt{\Delta}},
\end{align*}
where in the second inequality we use that $\vert \bm\delta \vert_1 \le \underline{\sigma}^{-1} \mathfrak{s}_\delta \vert \bm\Sigma \bm\delta \vert_\infty$, which follows from that
\begin{align*}
\underline{\sigma} \mathfrak{s}_\delta^{-1} \vert \bm\delta \vert_1^2
\le \underline{\sigma} \vert \bm\delta \vert_2^2
\le \bm\delta^\top \bm\Sigma \bm\delta
\le \vert \bm\delta \vert_1 \vert \bm\Sigma \bm\delta \vert_\infty.
\end{align*}
In the last inequality, we evoke the condition~\eqref{eq:sparse:dlb}.
Similarly, we have $\vert \bm\delta \vert_2 \le \underline{\sigma}^{-1} \sqrt{\mathfrak{s}_\delta} \vert \bm\Sigma \bm\delta \vert_\infty$ so that
\begin{align*}
\l\vert \bm\delta(\wh\cp_{\Mc}) - \bm\delta \r\vert_2 &\le \frac{ \bar{c}_1 \vert \bm\Sigma \bm\delta \vert_\infty^{-2} \sigma_X^2 \Psi^2 \log(p \log(n)) \cdot \vert \bm\delta \vert_2 }{\Delta}
\\
&\lesssim \frac{\sigma_X \Psi \sqrt{\mathfrak{s}_\delta \log(p \log(n))}}{\underline{\sigma} \sqrt{\Delta}}.
\end{align*}

\paragraph{When $\wh\cp = \wh\cp_{\Qc}$.} By Theorem~\ref{thm:qcscan}, we have with probability greater than $1 - c_2 (p\log(n))^{-1}$,
\begin{align*}
\l\vert \bm\delta(\wh\cp_{\Qc}) - \bm\delta \r\vert_1 &\le \frac{c_1 \vert \bm\Sigma \bm\delta \vert_2^{-2} \Vert \bm\Sigma \Vert \Psi^2 \sqrt{p\log\log(n)} \cdot \vert \bm\delta \vert_1 }{ \Delta}
\\
&\le \frac{c_1 \Vert \bm\Sigma \Vert \Psi^2 \sqrt{\mathfrak{s}_{\delta} p \log\log(n)}}{\underline{\sigma} \sqrt{\Delta} \cdot \sqrt{\vert \bm\Sigma \bm\delta \vert_2^2 \Delta}} 
\\
&\lesssim \frac{\Vert \bm\Sigma \Vert^{1/2} \Psi \sqrt{\mathfrak{s}_\delta} (p\log\log(n))^{1/4}}{\underline{\sigma} \sqrt{\Delta}},
\end{align*}
where in the second inequality we use that $\vert \bm\delta \vert_1 \le \underline{\sigma}^{-1} \sqrt{\mathfrak{s}_\delta} \vert \bm\Sigma \bm\delta \vert_2$, which follows from that
\begin{align*}
\underline{\sigma} \mathfrak{s}_\delta^{-1} \vert \bm\delta \vert_1^2
\le \underline{\sigma} \vert \bm\delta \vert_2^2
\le \bm\delta^\top \bm\Sigma \bm\delta
\le \underline{\sigma}^{-1} \vert \bm\Sigma \bm\delta \vert_2^2.
\end{align*}
In the last inequality, we evoke the condition~\eqref{eq:dense:dlb}.
Similarly, we have $\vert \bm\delta \vert_2 \le \underline{\sigma}^{-1} \vert \bm\Sigma \bm\delta \vert_2$ so that
\begin{align*}
\l\vert \bm\delta(\wh\cp_{\Qc}) - \bm\delta \r\vert_2 &\le \frac{ c_1 \vert \bm\Sigma \bm\delta \vert_2^{-2} \Vert \bm\Sigma \Vert \Psi^2 \sqrt{p\log\log(n)} \cdot \vert \bm\delta \vert_2 }{\Delta}
\\
&\lesssim \frac{\Vert \bm\Sigma \Vert^{1/2} \Psi (p\log\log(n))^{1/4}}{\underline{\sigma} \sqrt{\Delta}}.
\end{align*}

Finally, the conclusions follow by the triangle inequality.

\subsubsection{Proof of Proposition~\ref{prop:v}}

In what follows, we condition on the consistent detection of the change point guaranteed in Theorem~\ref{thm:qcscan}.
Throughout, we omit the subscript `Q' where there is no confusion.

Let us write
\begin{align*}
\l\vert \frac{n}{\wh\cp(n - \wh\cp)} T_{\wh\cp} - \bm\delta^\top \bm\Sigma^2 \bm\delta \r\vert \le \frac{n}{\wh\cp(n - \wh\cp)} \vert T_{\wh\cp} - f_{\wh\cp} \vert + \l\vert \frac{n}{\wh\cp(n - \wh\cp)} f_{\wh\cp} - \bm\delta^\top \bm\Sigma^2 \bm\delta \r\vert =: U_1 + U_2.
\end{align*}
Recall the definitions of $\mc A_l, \, l = 1, 2$, from the proof of Theorem~\ref{thm:qcscan}.
For $c_1$ large enough, we have $\vert \wh\cp - \cp \vert \le \Delta/4$ on $\mc A_1 \cap \mc A_2$, which implies that
\begin{align}
\frac{\wh\cp (n - \wh\cp)}{n} \le \frac{\Delta}{4}, \text{ \ i.e. \ }
\frac{n}{\wh\cp(n - \wh\cp)} \ge \frac{4}{\Delta}.
\nn % \label{eq:k:ratio}
\end{align}
Also, Steps~1--2 of the proof of Theorem~\ref{thm:qcscan} show that, on $\mc A_1 \cap \mc A_2$, the grid of probe points on which $T_k$ is evaluated in Algorithm~\ref{alg:os}, denoted by $\mc G$, has its cardinality bounded as $\vert \mc G \vert \lesssim \log^2(n)$.
This, in combination with Proposition~\ref{prop:dense}, gives
\begin{align*}
& \p\l[ U_1 \ge \frac{C \Vert \bm\Sigma \Vert \Psi^2 \sqrt{p \log\log(n)}}{\Delta}; \mc A_1 \cap \mc A_2 \r] 
\\
\le & \; \p\l[ \vert T_{\wh\cp} - f_{\wh\cp} \vert \ge \frac{C \Vert \bm\Sigma \Vert \Psi^2 \sqrt{p \log\log(n)}}{4} ; \mc A_1 \cap \mc A_2 \r] 
\\
\le & \; \p\l[ \max_{k \in \mc G} \vert T_k - f_k \vert \ge \frac{C \Vert \bm\Sigma \Vert \Psi^2 \sqrt{p \log\log(n)}}{4} ; \mc A_1 \cap \mc A_2 \r] 
\\
\lesssim & \; \log^2(n) \cdot C_1 \exp(-3\log\log(n)) \lesssim C_1 \bigl(\log(n)\bigr)^{-1},
\end{align*}
provided that $C \ge 12 C_0$.
Since $\p(\mc A_1 \cap \mc A_2) \ge 1 - c_2 \bigl(\log(n)\bigr)^{-1}$ with some $c_2$ that depends on~$C_1$, we have
\begin{align*}
\p\l[ U_1 \ge \frac{C \Vert \bm\Sigma \Vert \Psi^2 \sqrt{p \log\log(n)}}{\Delta} \r] \lesssim \bigl(\log(n)\bigr)^{-1}
\end{align*}
with the constant omitted in $\lesssim$ depending only on $C_1$.

In studying the `bias' $U_2$, WLOG, we suppose that $\wh\cp \le \cp$; analogous arguments apply to the case when $\wh\cp \ge \cp$.
Then, from~\eqref{eq:f}, it follows that
\begin{align*}
U_2 \le & \, \l\vert 1 - \frac{(n - \cp)^2}{(n - \wh\cp)^2} \r\vert \bm\delta^\top \bm\Sigma^2 \bm\delta  + \l\vert \frac{1}{\wh\cp} + \frac{\cp - \wh\cp}{(n - \wh\cp)^2} \r\vert \bm\beta_0^\top \bm\Sigma^2 \bm\beta_0
+ \l\vert \frac{1}{n - \wh\cp} - \frac{\cp - \wh\cp}{(n - \wh\cp)^2} \r\vert \bm\beta_1^\top \bm\Sigma^2 \bm\beta_1 
\\
=: & \, U_{2, 1} + U_{2, 2} + U_{2, 3}.
\end{align*}
On $\mc A_1 \cap \mc A_2$, it holds that
\begin{align*}
U_{2, 1} &= \frac{\vert \wh\cp - \cp \vert (n - \wh\cp + n - \cp)}{(n - \wh\cp)^2} \bm\delta^\top \bm\Sigma^2 \bm\delta
\le \frac{2 \vert \wh\cp - \cp \vert}{n - \cp} \bm\delta^\top \bm\Sigma^2 \bm\delta
\\
&\le \frac{2 c_1 \Vert \bm\Sigma \Vert \Psi^2 \sqrt{p\log\log(n)} }{\Delta}.
\end{align*}
As for $U_{2, 2}$, thanks to~\eqref{eq:dense:dlb},
\begin{align*}
U_{2, 2} \le \frac{7}{3\Delta} \Vert \bm\Sigma \Vert \vert \bm\Sigma^{1/2} \bm\beta_0 \vert_2^2 \le \frac{7}{3c_0 \sqrt{p \log\log(n)}} \bm\delta^\top \bm\Sigma^2 \bm\delta
\end{align*}
and $U_{2, 3}$ is similarly bounded. 
Collecting the bounds on $U_1$ and $U_2$, the proof is complete.

\subsubsection{Supporting lemmas}

\begin{lem}[Lemma~B.5 of \citealp{cho2022high}]
\label{lem:hw}
Suppose that $\mbf z_t \sim_{\iid} \mc N_m(\mbf 0, \mbf I_m)$ for all~$t\in[n]$.
Then there exists some universal constant $C \in (0, \infty)$ % and any sequence $\eps_n \to \infty$ as $n \to \infty$ while $n^{-1} \eps_n \to 0$, such that for all $0 \le s < e \le n$ with $e - s \ge \eps_n$ and $z \in (0, \eps_n^{1/2})$,
such that for all $0 \le s < e \le n$ and $z > 0$,
\begin{align*}
\sup_{\mbf a, \mbf b \in \mathbb{B}_2(1)} \p\l\{ \frac{1}{\sqrt{e - s}} \l\vert \sum_{t = s + 1}^e \l( \mbf a^\top \mbf z_t \mbf z_t^\top \mbf b - \mbf a^\top \mbf b \r) \r\vert \ge C z \r\} \le 2 \exp\l[ - \min(z^2,  z \sqrt{e - s}) \r].
\end{align*}
\end{lem}

\begin{lem}
\label{lem:f:bar}
Recall $\bar{\mbf f}_k$ defined in~\eqref{eq:f:bar}.
Then, when $1 \le m \le k \le \cp$. we have
\begin{align*}
\vert \bar{\mbf f}_k \vert_\infty - \vert \bar{\mbf f}_m \vert_\infty  = \vert \bar{\mbf f}_k - \bar{\mbf f}_m \vert_\infty \ge \frac{3(k - m)(n - \cp)(k + n - m)}{8\sqrt{n k (n - k)} (n - m)} \vert \bm\Sigma \bm\delta \vert_\infty
\end{align*}
and for $\cp \le k \le m \le n$,
\begin{align*}
\vert \bar{\mbf f}_k \vert_\infty - \vert \bar{\mbf f}_m \vert_\infty = \vert \bar{\mbf f}_k - \bar{\mbf f}_m \vert_\infty \ge \frac{3(m - k) \cp (n - k + m)}{8\sqrt{n k (n - k)} m} \vert \bm\Sigma \bm\delta \vert_\infty.
\end{align*}
\end{lem}

\begin{proof}
For $m \le k \le \cp$, we have
\begin{align*}
& \sqrt{\frac{k (n - k)}{n}} \frac{n - \cp}{n - k} - \sqrt{\frac{m (n - m)}{n}} \frac{n - \cp}{n - m} 
\\
=& \, \sqrt{\frac{k (n - k)}{n}} \frac{n - \cp}{n - k} \l( 1 - \sqrt{\frac{m (n - k)}{k (n - m)}} \r)
\\
\ge& \, \sqrt{\frac{k (n - k)}{n}} \frac{n - \cp}{n - k} \l[ 1 - \l(1 - \frac{k - m}{2k} \r) \l(1 - \frac{k - m}{2(n - m)} \r) \r] 
\\
=& \, \sqrt{\frac{k (n - k)}{n}} \frac{n - \cp}{n - k} \cdot \frac{k - m}{2k (n - m)} \l( k + n - m - \frac{k - m}{2} \r)
\\
\ge& \, \frac{3(k - m)(n - \cp)(k + n - m)}{8\sqrt{n k (n - k)} (n - m)}.
\end{align*}
Similar arguments apply to the case of $\cp \le k \le m$.
\end{proof}

%%%%%%%%%%%%%%%%%%%%%%%%%%%%%%%%%%%%%%%
\begin{lem}[Lemma~B.1 of \citealp{VeGa18}] 
\label{lem:op:xx}
Let $\mbf x_t \sim_{\iid} \mc N_p(\mbf 0, \bm \Sigma)$ for $t \in[n]$ and $ \mbf X = [\mbf x_1, \ldots, \mbf x_n]^\top$. Then, for all $z > 0$, we have
\begin{align*}
\p\l( \l\Vert \mbf X\mbf X^\top - \tr(\bm\Sigma)\mbf I_n \r\Vert \le 2 \sqrt{\tr(\bm\Sigma)} \Vert \bm\Sigma \Vert^{1/2} (\sqrt n + 10 + \sqrt{2z}) + 3 \Vert \bm\Sigma \Vert(n + 100 + 2z) \r) \ge 1 - 2 e^{-z}.
\end{align*}
As a consequence, for all $n \ge 20$, we have
\begin{align*}
\p\l( \Vert \mbf X\mbf X^\top - \tr(\bm\Sigma)\mbf I_n \Vert \le 25 \Vert \bm\Sigma \Vert (\sqrt{np} + n) \r) \ge 1 - 2e^{-n}.
\end{align*}
\end{lem}
See also \citet{KoLo17} for a systematic treatment of this type of concentration inequalities. 

\begin{lem}[Lemma~C.1 of \citealp{VeGa18}] 
\label{lem:chisq}
Let $\mbf x \sim \mc N_m(\mbf 0, \mbf I_m)$ and $\mbf A \in \R^{m \times m}$ be a symmetric matrix.
For any $z > 0$,
\begin{align*}
\p\l( \mbf x^\top \mbf A \mbf x \ge \tr(\mbf A) + 2 \Vert \mbf A \Vert_F \sqrt{z} + 2 \Vert \mbf A \Vert z \r) \le e^{-z}.
\end{align*}
When $\mbf A = \mbf I_m$, the above bound simplifies as
\begin{align*}
\p( \mbf x^\top \mbf x \ge m + 2 \sqrt{mz} + 2z ) \le e^{-z},
\end{align*}
and we also have
\begin{align*}
\p( \mbf x^\top \mbf x \le m - 2 \sqrt{mz} ) \le e^{-z}.
\end{align*}
\end{lem}

\begin{lem}[Equation (27) of \citealp{VeGa18}] 
\label{lem:tr:xx}
Let $\mbf x_t \sim_{\iid} \mc N_p(\mbf 0, \bm \Sigma)$ for $t \in[n]$ and $ \mbf X = [\mbf x_1, \ldots, \mbf x_n]^\top$. Then, for all $z \in (0, np)$, we have
\begin{align*}
\p\l( \l\vert \tr\l(\mbf X\mbf X^\top - \tr(\bm\Sigma) \mbf I_n\r) \r\vert  \ge 8 \Vert \bm\Sigma \Vert \sqrt{npz} \r) \le 2 e^{-z}.
\end{align*}
\end{lem}
Note that Lemma~\ref{lem:tr:xx} also follows from Lemma~\ref{lem:chisq}.

% \begin{lem}[Theorem~6.2.1 of \cite{vershynin2018high}]
% \label{lem:hw:ineq}
% Let $\mbf x = (X_1, \ldots, X_m)^\top \in \R^m$ be a random vector with independent, mean zero, sub-Gaussian coordinates. 
% Let $\mbf A$ be an $m \times m$ matrix. 
% Then for any $z \ge 0$, we have
% \begin{align*}
% \p\l( \vert \mbf x^\top \mbf A \mbf x - \E(\mbf x^\top \mbf A \mbf x) \vert \ge z \r) \le 2 \exp\l[ - c \min\l( \frac{z^2}{K^4 \Vert \mbf A \Vert_F^2}, \frac{z}{K^2 \Vert \mbf A \Vert}\r) \r],
% \end{align*}
% where $K = \max_{t \in [m]} \Vert X_t \Vert_{\psi_2}$.
% \end{lem}

\begin{lem}[Lemma~B.2 of \citealp{VeGa18}]
\label{lem:gauss:chaos}
Let $\mbf x_t \sim_{\iid} \mc N_p(\mbf 0, \bm \Sigma)$ for $t \in[n]$, $p \ge n$ and $ \mbf X = [\mbf x_1, \ldots, \mbf x_n]^\top$. Then, for all $z \ge 0$, any $\bm\alpha, \bm\beta \in \R^p$ and some constant $C > 0$, 
\begin{align*}
& \p\l[\l\vert  \frac{1}{n^2} \bm\alpha^\top \mbf X^\top \l( \mbf X \mbf X^\top - \tr(\bm\Sigma) \mbf I_n \r) \mbf X \bm\beta - \l(1 + \frac{1}{n}\r) \bm\alpha^\top \bm\Sigma^2 \bm\beta \r\vert \r.
\\
& \qquad \l. \ge C \Vert \bm\Sigma \Vert \vert \bm\Sigma^{1/2} \bm\alpha \vert_2 \vert \bm\Sigma^{1/2} \bm\beta \vert_2 \frac{\sqrt{pz}}{n} \r] \le 2 \exp\l(-\min\bigl\{z, (nz)^{1/4}\bigr\}\r).
\end{align*}
% \begin{align*}
% \p\l[ \l\vert 
% \frac{1}{n^2} \bm\beta^\top \mbf X^\top (\mbf X \mbf X^\top - \tr(\bm\Sigma) \mbf I_n) \mbf X \bm\beta - \l(1 + \frac{1}{n} \r) \vert \bm\Sigma \bm\beta \vert_2^2 \r\vert \ge C \vert \bm\Sigma^{1/2} \bm\beta \vert_2^2 \Vert \bm\Sigma \Vert \frac{\sqrt{pz}}{n} \r] \le 2 e^{-z}.
% \end{align*}
\end{lem}

\begin{proof}
While the lemma is stated slightly more generally than Lemma~B.2 of \cite{VeGa18}, the proof takes analogous arguments and thus is omitted; see also the proof of Lemma~\ref{lem:gauss:chaos:two}.
\end{proof}

\begin{lem}
\label{lem:gauss:chaos:two}
Let $\mbf x_t \sim_{\iid} \mc N_p(\mbf 0, \bm \Sigma)$ for $t \in[n]$, $p \ge n$ and $ \mbf X = [\mbf x_1, \ldots, \mbf x_n]^\top$. 
Consider some $0 < k < n$. 
For all $z \ge 0$, any $\bm\alpha, \bm\beta \in \R^p$ and some constant $C > 0$, we have
\begin{align*}
& \p\l[\l\vert  \frac{1}{k (n - k)} \bm\alpha^\top \mbf X_{0, k}^\top \mbf X_{0, k} \mbf X_{k, n}^\top \mbf X_{k, n} \bm\beta - \bm\alpha^\top \bm\Sigma^2 \bm\beta \r\vert \r. 
\\ 
& \qquad \l. \ge C \Vert \bm\Sigma \Vert \vert \bm\Sigma^{1/2} \bm\alpha \vert_2 \vert \bm\Sigma^{1/2} \bm\beta \vert_2 \sqrt{\frac{pz}{k (n - k)}} \r] \le 2 \exp\l(-\min\bigl\{z, (nz)^{1/4}\bigr\}\r).
\end{align*}
\end{lem}

\begin{proof}
The proof largely follows the proof of Lemma~B.2 of \cite{VeGa18}, which in turn relies on the deviation inequality derived for Gaussian chaos in Theorem~1.3 of \cite{adamczak2015concentration}. 

Let us denote the eigendecomposition of $\bm\Sigma$ by $\bm\Sigma = \mbf E\mbf M \mbf E^\top$, where $\mbf M$ is the diagonal matrix with eigenvalues, $\mu_i, \, i \in [p]$, on the diagonal, and $\mbf E$ the orthogonal matrix containing eigenvectors.
Then, we have $\mbf Z = [Z_{ti}, \, t \in [n], \, i \in [p]] := \mbf X \mbf E \mbf M^{-1/2}$ satisfy $Z_{ti} \sim_{\iid} \mc N(0, 1)$.
Next, we may write
\begin{align*}
W &:= \frac{1}{k (n - k)} \bm\alpha^\top \mbf X_{0, k}^\top \mbf X_{0, k} \mbf X_{k, n}^\top \mbf X_{k, n} \bm\beta
\\
&= \frac{1}{k (n - k)} (\underbrace{\mbf M^{1/2} \mbf E^\top \bm\alpha}_{\mbf a})^\top \mbf Z_{0, k}^\top \mbf Z_{0, k} \mbf M \mbf Z_{k, n}^\top \mbf Z_{k, n} \underbrace{\mbf M^{1/2} \mbf E^\top \bm\beta}_{\mbf b}.
\end{align*}
We can further write $f(\mbf Z) := k (n - k) W$, as
\begin{align*}
f(\mbf Z) &= \sum_{t = 1}^k \sum_{u = k + 1}^n \sum_{i, j, \ell = 1}^p a_i \mu_j b_\ell Z_{t i} Z_{t j} Z_{u j} Z_{u \ell}
\\
&= \sum_{t_1, t_2, t_3, t_4 = 1}^n \sum_{i_1, i_2, i_3, i_4 = 1}^p \underbrace{a_{i_1} \mu_{i_2} b_{i_4} \mathbb{I}_{\{t_1 = t_2 \le k\}} \mathbb{I}_{\{t_3 = t_4 > k \}} \mathbb{I}_{\{i_2 = i_3\}}}_{c_{(t_1, i_1), (t_2, i_2), (t_3, i_3), (t_4, i_4)}} \cdot Z_{t_1 i_1} Z_{t_2 i_2} Z_{t_3 i_3} Z_{t_4 i_4}.
\end{align*}
To proceed further, we need some notation. 
For positive integers $d$ and $q$, we denote a $d$-indexed matrix $\mbf C = (c_{i_1, \ldots, i_d})_{i_1, \ldots,  i_d = 1}^q$; also we write $\mbf i = ( i_1, \ldots,  i_d) \in [q]^d$ and $\mbf i_J = ( i_l, \, l \in J)$ for $J \subseteq [d]$.
Also, let $\mc P_d$ the set of partitions of $[d]$ into non-empty disjoint subsets.
Given a partition $\mc J = \{J_1, \ldots, J_\ell\}$, we define the tensor product norm:%(Equation (32) of \cite{VeGa18})
\begin{align*}
\Vert \mbf C \Vert_{\mc J} := \sup\l\{ \sum_{\mbf i \in [q]^d} c_{\mbf i} \prod_{l = 1}^\ell \mbf z^{(l)}_{\mbf i_{J_l}} : \, \Vert \mbf z^{(l)} \Vert_F \le 1, \, 1 \le l \le \ell \r\},
\end{align*}
where $\mbf z^{(l)}$ is a $\vert J_l \vert$-indexed matrix.
Then, $f(\mbf Z)$ is a polynomial of degree $D = 4$ of the $q = np$ independent standard Gaussian random variables $Z_{ti}$'s, to which Theorem~1.3 of \cite{adamczak2015concentration} is applicable, such that with some constant $C_0 > 0$,
\begin{align}
\label{eq:dev:fz}
\p\l( \vert f(\mbf Z) - \E(f(\mbf Z)) \vert \ge w \r) \le 2 \exp\l[ - \min_{1 \le d \le D} \min_{\mc J \in \mc P_d} \l( \frac{C_0^d w}{\Vert \E(\partial^d f(\mbf Z) \Vert_{\mc J}} \r)^{2/\vert \mc J \vert} \r],
\end{align}
where $\partial^d f$ denotes the $d$th derivative of $f$; we later set
\begin{align}
\label{eq:wz}
w = C \sqrt{pk(n - k) z} \Vert \bm\Sigma \Vert \vert \vert \mbf a \vert_2 \vert \mbf b \vert_2
\end{align}
with some constant $C > 0$ that depends on $C_0$. 
When studying the tail bound in~\eqref{eq:dev:fz}, we only have to consider the derivatives of order $d \in \{2, 4\}$ since all the other terms are null.

\paragraph{When $d = 4$.} Denoting by $\mbf C = [c_{(t_1, i_1), (t_2, i_2), (t_3, i_3), (t_4, i_4)}]$ with $(t_l, i_l) \in [n] \times [p]$, we derive a bound on $\Vert \E(\partial^4 f(\mbf Z)) \Vert_{\mc J}$ from that on $\Vert \mbf C \Vert_{\mc J}$ by triangle inequality, as
$[\E(\partial^4 f(\mbf Z))]_{\mbf i_1, \mbf i_2, \mbf i_3, \mbf i_4} = \sum_\pi C_{\mbf i_{\pi(1)}, \mbf i_{\pi(2)}, \mbf i_{\pi (3)}, \mbf i_{\pi(4)}}$, where the sum is taken over all permutations $\pi$ of $[4]$.
When $\mc J = \{[4]\}$ (hence $\vert \mc J \vert = 1$), it holds
\begin{align*}
\Vert \mbf C \Vert_{\mc J}^2 = \Vert \mbf C \Vert_F^2 &= \sum_{t = 1}^k \sum_{u = k + 1}^n \sum_{i, j, \ell = 1}^p a_i^2 \mu_j^2 b_\ell^2
\\
&= k (n - k) \tr(\bm\Sigma^2) \vert \mbf a \vert_2^2 \vert \mbf b \vert_2^2 
\le k (n - k) p \Vert \bm\Sigma \Vert^2 \vert \mbf a \vert_2^2 \vert \mbf b \vert_2^2.
\end{align*}
When $\vert \mc J \vert = 2$, by the construction of the elements of $\mbf C$ that involve indicator functions, there is a loss by a multiplicative factor of $n$ or $p$ in $\Vert \mbf C \Vert_{\mc J}^2$ compared to $\Vert \mbf C \Vert_F^2$, i.e.\
\begin{align*}
\Vert \mbf C \Vert_{\mc J}^2 \le \frac{\Vert \mbf C \Vert_F^2}{\min(n, p)} \le \frac{k(n - k)}{n} p \Vert \bm\Sigma \Vert^2 \vert \mbf a \vert_2^2 \vert \mbf b \vert_2^2. 
\end{align*}
For $\mc J$ with $\vert \mc J \vert \ge 3$, by the definition of $\Vert \cdot \Vert_{\mc J}$, there exists some $\mc J'$ with $\vert \mc J' \vert < \vert \mc J \vert$ (obtained by merging two elements of $\mc J$) such that $\Vert \mbf C \Vert_{\mc J} \le \Vert \mbf C \Vert_{\mc J'}$.
In summary,
\begin{align}
& \min_{\mc J \in \mc P_4} \l( \frac{C_0^8 w^2}{\Vert \E(\partial^4 f(\mbf Z)) \Vert_{\mc J}^2} \r)^{1/\vert \mc J\vert} 
\nn \\
\ge & \, \min \l[ \frac{C_0^8 w^2}{k (n - k) p \Vert \bm\Sigma \Vert^2 \vert \mbf a \vert_2^2 \vert \mbf b \vert_2^2}, \, \min_{l \in \{2, 3, 4\}} \l( \frac{C_0^8 n w^2}{k(n - k) p \Vert \bm\Sigma \Vert^2 \vert \mbf a \vert_2^2 \vert \mbf b \vert_2^2} \r)^{1/l} \r]\nn \\
\ge & \, \min \l[ z, (nz)^{1/4}\r]
\label{eq:d:four}
\end{align}
where in the last inequality we use $w$ as in~\eqref{eq:wz} with sufficiently large $C$. 
% \footnote{
% \begin{align*}
% & \frac{w^2}{k (m - k) p \Vert \bm\Sigma \Vert^2 \vert \mbf a \vert_2^2 \vert \mbf b \vert_2^2} < \l( \frac{m w^2}{k(m - k) p \Vert \bm\Sigma \Vert^2 \vert \mbf a \vert_2^2 \vert \mbf b \vert_2^2} \r)^{1/4}
% \\
% \Longleftrightarrow & \, w^6 < m k^3 (m - k)^3 p^3 \Vert \bm\Sigma \Vert^6 \vert \mbf a \vert_2^6 \vert \mbf b \vert_2^6
% \\
% \Longleftrightarrow & \, z^3  < m
% \end{align*}
% }
 
\paragraph{When $d = 2$.} 
Let $\mbf C' = [c'_{(t_1, i_1), (t_2, i_2)}] = \E(\partial^2 f(\mbf Z))$ with $(t_l, i_l) \in [n] \times [p]$.
Recall that each summand of $f(\mbf Z)$ involve exactly two terms with the index $t_1 \le k$ and the other two with $t_2 > k$, from which we have 
\begin{align*}
c'_{(t_1, i_1), (t_2, i_2)} = \l\{ \begin{array}{ll}
(n - k)(a_{i_1} b_{i_2} \mu_{i_2}  + a_{i_2} b_{i_1} \mu_{i_1}) & \text{when \ } 1 \le t_1 = t_2 \le k, \\
k( a_{i_1} b_{i_2} \mu_{i_1}  + a_{i_2} b_{i_1} \mu_{i_2} ) & \text{when \ } k + 1 \le t_1 = t_2 \le n, \\
0 & \text{when \ } t_1 \ne t_2. \\
\end{array}\r.
\end{align*}
Therefore, when $\mc J = \{[2]\}$,
\begin{align*}
\Vert \mbf C' \Vert_{\mc J}^2 = \Vert \mbf C' \Vert_F^2 \le 4 \Vert \bm\Sigma \Vert^2 \vert \mbf a \vert_2^2 \vert \mbf b \vert_2^2 n k (n - k).
\end{align*}
When $\mc J = \{\{1\}, \{2\}\}$, we have $\Vert \mbf C' \Vert_{\mc J} = \Vert \mbf C' \Vert \le 2 \max(k, n-k)\Vert \bm\Sigma \Vert \vert \mbf a \vert_2 \vert \mbf b \vert_2$.
In summary,
\begin{align*}
& \min_{\mc J \in \mc P_2} \l(\frac{C_0^4w^2}{\Vert \E(\partial^2 f(\mbf Z)) \Vert_{\mc J}^2} \r)^{1/\vert \mc J\vert} 
 \\
\ge & \, \min \l[ \frac{C_0^4 w^2}{4 n k (n - k) \Vert \bm\Sigma \Vert^2 \vert \mbf a \vert_2^2 \vert \mbf b \vert_2^2}, \, \l( \frac{C_0^4w^2}{4 \max(k, n-k)^2 \Vert \bm\Sigma \Vert^2 \vert \mbf a \vert_2^2 \vert \mbf b \vert_2^2} \r)^{1/2} \r]
\end{align*}
which dominates the RHS of~\eqref{eq:d:four}, when we set $w$ as in~\eqref{eq:wz} with sufficiently large $C$.
\medskip
% \footnote{HC: while \cite{VeGa18} considers $n \le p \le n^2$, I think a slightly stronger upper bound is required because
% \begin{align*}
% & \frac{w^2}{m k (m - k) \Vert \bm\Sigma \Vert^2 \vert \mbf a \vert_2^2 \vert \mbf b \vert_2^2} < \l( \frac{w^2}{k (m - k) \Vert \bm\Sigma \Vert^2 \vert \mbf a \vert_2^2 \vert \mbf b \vert_2^2} \r)^{1/2}
% \\
% \Longleftrightarrow & \, w < m \sqrt{k (m - k)} \Vert \bm\Sigma \Vert \vert \mbf a \vert_2 \vert \mbf b \vert_2
% \\
% \Longleftrightarrow & \, z  < (m/\sqrt{p})^2
% \end{align*}
% where we require $z \asymp \log(n)$ or $\log\log(n)$ with $m \ge \eps_n$ whenever this lemma is evoked.
% }\footnote{\cred We only need to require that the first term in the last display for $d = 4$ is smaller than the second term here, namely,
% \begin{align*}
%     &\frac{w^2}{k (m - k) p \Vert \bm\Sigma \Vert^2 \vert \mbf a \vert_2^2 \vert \mbf b \vert_2^2} < \l( \frac{w^2}{k (m - k) \Vert \bm\Sigma \Vert^2 \vert \mbf a \vert_2^2 \vert \mbf b \vert_2^2} \r)^{1/2} \\
%     \iff & z <p^2,
% \end{align*}
% which is the case when $z < m^{1/3} \le n^{1/3} \le p^{1/3}$.}

Collecting the bounds and plugging them into~\eqref{eq:dev:fz}, the conclusion follows.
\end{proof}

% \begin{lem}
% \label{lem:gauss:chaos:three}
% For all $z \in (0, n^{1/3})$, any $\bm\alpha, \bm\beta \in \R^p$ and some constant $C > 0$, we have
% \begin{align*}
% & \p\l[\l\vert  \frac{1}{n^2} \bm\alpha^\top \mbf X^\top \l( \mbf X \mbf X^\top - \tr(\bm\Sigma) \mbf I_n \r) \mbf X \bm\beta - \l(1 + \frac{1}{n}\r) \bm\alpha^\top \bm\Sigma^2 \bm\beta \r\vert \r.
% \\
% & \qquad \l. \ge C \Vert \bm\Sigma \Vert \vert \bm\Sigma^{1/2} \bm\alpha \vert_2 \vert \bm\Sigma^{1/2} \bm\beta \vert_2 \frac{\sqrt{pz}}{n} \r] \le 2 e^{-z}.
% \end{align*}
% \end{lem}

% \begin{proof}
% The proof follows the same steps as the proof of Lemmas~\ref{lem:gauss:chaos} and~\ref{lem:gauss:chaos:two} verbatim and thus is omitted.
% \end{proof}

\begin{lem}
\label{lem:f} Recall $f_k$ defined in~\eqref{eq:f}. Then, when $1 \le m \le k \le \cp$, we have
\begin{align*}
f_k - f_m = \frac{(k - m)(n - \cp)}{(n - k)(n - m)} \l( (n - \cp) \vert \bm\Sigma \bm\delta \vert_2^2 - \vert \bm\Sigma \bm\beta_0 \vert_2^2 + \vert \bm\Sigma \bm\beta_1 \vert_2^2 \r)
\end{align*}
and for $\cp \le k \le m \le n$,
\begin{align*}
f_k - f_m = \frac{(m - k)\cp}{km} \l( \cp \vert \bm\Sigma \bm\delta \vert_2^2 + \vert \bm\Sigma \bm\beta_0 \vert_2^2 - \vert \bm\Sigma \bm\beta_1 \vert_2^2 \r).
\end{align*}
Further, under the condition~\eqref{eq:delta:beta}, we have
\begin{align*}
f_k - f_m \ge \frac{\vert \bm\Sigma \bm\delta \vert_2^2 \vert k - m \vert}{2} \l[ \frac{(n - \cp)^2}{(n - k)(n - m)} \mathbb{I}_{\{1 \le m \le k \le \cp \}} + \frac{\cp^2}{km} \mathbb{I}_{\{\cp \le k \le m \le n\}} \r].
\end{align*}
\end{lem}

\begin{proof}
For $m \le k \le \cp$, we have
\begin{align*}
f_k - f_m &= \frac{(n - \cp)^2}{n} \l(\frac{k}{n - k} - \frac{m}{n - m}\r) \vert \bm\Sigma \bm\delta \vert_2^2 
\\
& \qquad - \frac{1}{n}\l( \frac{(\cp - k)k}{n - k} - \frac{(\cp - m)m}{n - m} - k + m \r) \l( \vert \bm\Sigma \bm\beta_1 \vert_2^2 - \vert \bm\Sigma \bm\beta_0 \vert_2^2 \r) 
\\
&= \frac{(k - m)(n - \cp)}{(n - k)(n - m)} \l[ (n - \cp) \vert \bm\Sigma \bm\delta \vert_2^2 + \vert \bm\Sigma \bm\beta_1 \vert_2^2 - \vert \bm\Sigma \bm\beta_0 \vert_2^2 \r]
\\
&\ge \frac{(k - m)(n - \cp)^2}{2(n - k)(n - m)} \vert \bm\Sigma \bm\delta \vert_2^2,
\end{align*}
where the last inequality follows provided that $(n - \cp) \vert \bm\Sigma \bm\delta \vert_2^2 \ge 2 \vert \bm\Sigma \bm\beta_0 \vert_2^2$ under~\eqref{eq:delta:beta}.
Similar calculations yield the claim when $m \ge k \ge \cp$.
\end{proof}

\subsection{Proofs for the results in Section~\ref{sec:four}}

\subsubsection{Proof of Proposition~\ref{prop:adapt}}

Upon the close inspection of the proof of Theorem~\ref{thm:qcscan}, it is clear that the conclusions of Theorems~\ref{thm:mcscan} and~\ref{thm:qcscan} hold on the following events
\begin{subequations}
\begin{align}
\bar{\mc E} &:= \bigcap_{k \in \mc D \cap \bar{\mc Q}} \l\{ \vert \mbf X^\top \wt{\mbf Y}(k) - \bar{\mbf f}_k \vert_\infty \le \bar{C}_0' \sigma_X \Psi \sqrt{\log(p\log(n))}  \r\}, \\
\mc E &:= \bigcap_{k \in \mc D \cap \mc Q} \l\{ \vert T_k - f_k \vert \le C_0' \Vert \bm\Sigma \Vert \Psi^2 \sqrt{p\log\log(n)}  \r\},
\end{align}
\label{eq:set:e} 
\end{subequations}
where $\mc D$ denotes the dyadic grid in Algorithm~\ref{alg:os} and $\bar{\mc Q}$ and $\mc Q$ deterministic grids that depend only on $\bar{\mbf f}_k$ or $f_k$, respectively (see Step~2 in the proof of Theorem~\ref{thm:qcscan}), and $\bar{C}_0'$ and $C_0'$ some universal constants that depend on $\bar{C}_0$ and $C_0$, respectively.
By Propositions~\ref{prop:sparse} and~\ref{prop:dense}, and from that $\vert \mc D \vert \le 2\log_2(n)$ and $\vert \bar{\mc Q} \vert, \vert \mc Q \vert \le 16\log^2(n)$, it follows that $\p(\bar{\mc E} \cap \mc E) \ge 1 - \bar{c}_2 (p\log(n))^{-1} - c_2 \bigl(\log(n)\bigr)^{-1}$ with $\bar{c}_2$ and $c_2$ depending on $\bar{C}_1$ and $C_1$, respectively.
Below we condition our arguments on $\bar{\mc E} \cap \mc E$. 

Notice that when $q = 0$, we have $\bar{\mbf f}_k = \mbf 0$ and $f_k \le 2 \Vert \bm\Sigma \Vert \Psi^2$. 
Hence from the construction of $\bar{\mc E}$ and $\mc E$, it is apparent that for large enough $\bar{c}$ and $c$, we have $\wh q = 0$ when $q = 0$, and thus $\wh\cp_{\Oc} = n$.
When $q = 1$, by the arguments adopted in Step~1 of the proof of Theorem~\ref{thm:qcscan}, we obtain $\wh q = 1$ provided that either~\eqref{eq:sparse:dlb} or~\eqref{eq:dense:dlb} holds.

Suppose that $\mathbb{I}_{\{\bar{T}_{\wh\cp_{\Mc}} \le \zeta_{\Mc} \}} \cdot \mathbb{I}_{\{T_{\wh\cp_{\Qc}} > \zeta_{\Qc} \}} = 1$.
Then it follows that~\eqref{eq:sparse:dlb} is not met, and thus by assumption, the condition~\eqref{eq:dense:dlb} holds, which guarantees that $\wh\cp_{\Qc}$ attains the estimation rate derived in Theorem~\ref{thm:qcscan}.
Then,
\begin{align*}
\frac{\vert \bm\Sigma \bm\delta \vert_\infty^2 \Delta}{\sigma_X^2 \log(p\log(n))} \lesssim \Psi^2 \lesssim \frac{\vert \bm\Sigma \bm\delta \vert_2^2 \Delta}{{ \Vert \bm\Sigma \Vert} \sqrt{p\log\log(n)}},
\end{align*}
so the conclusion of the proposition follows. 
Analogous arguments as above apply when $\mathbb{I}_{\{\bar{T}_{\wh\cp_{\Mc}} > \zeta_{\Mc} \}} \cdot \mathbb{I}_{\{T_{\wh\cp_{\Qc}} \le \zeta_{\Qc} \}} = 1$.

Next, suppose that $\mathbb{I}_{\{\bar{T}_{\wh\cp_{\Mc}} > \zeta_{\Mc} \}} \cdot \mathbb{I}_{\{T_{\wh\cp_{\Qc}} > \zeta_{\Qc} \}} = 1$.
By the definition of $\bar{\mbf f}_k$ in~\eqref{eq:f:bar} and Lemma~\ref{lem:f:bar}, it follows that
\begin{align}
\zeta_{\Mc} < \bar{T}_{\wh\cp_{\Mc}} \le \vert \bm\Sigma \bm\delta \vert_\infty\sqrt{\Delta} + \bar{C}_0' \sigma_X \Psi \sqrt{\log(p\log(n))},
\label{eq:tbar:upper}
\end{align}
indicating that~\eqref{eq:sparse:dlb} is satisfied with some $\bar{c}_0$ that depends on $\bar{c}$.
Thus, it follows that $\vert \wh\cp_{\Mc} - \cp \vert \le \eta \Delta$ for some $\eta \in (0, 1)$ that depends on $\bar{c}_0$ (and consequently $\bar{c}$), which leads to
\begin{align}
\sqrt{\frac{1 - \eta}{2(1 + \eta)}} \vert \bm\Sigma \bm\delta \vert_\infty \sqrt{\Delta} - \bar{C}_0' \sigma_X \Psi \sqrt{\log(p\log(n))} \le \bar{T}_{\wh\cp_{\Mc}}.
\label{eq:tbar:lower}
\end{align}
That is, $\bar{T}_{\wh\cp_{\Mc}}^2 \asymp \vert \bm\Sigma \bm\delta \vert_\infty^2 \Delta$ from~\eqref{eq:tbar:upper} and~\eqref{eq:tbar:lower}.
Similarly, by the definition of $f_k$ in~\eqref{eq:f} and Lemma~\ref{lem:f} and \eqref{eq:delta:beta}, for large enough $c_0$ in~\eqref{eq:dense:dlb} that depends on $c$, we have
\begin{multline}
\frac{1 - \eta}{1 + \eta} \vert \bm\Sigma \bm\delta \vert_2^2 \Delta
- C_0' \Vert \bm\Sigma \Vert \Psi^2 \sqrt{p\log\log(n)}
\\
\le T_{\wh\cp_{\Qc}} \le \frac{3}{2} \vert \bm\Sigma \bm\delta \vert_2^2 \Delta + C_0' \Vert \bm\Sigma \Vert \Psi^2 \sqrt{p\log\log(n)},
\label{eq:t:sim}
\end{multline}
leading to $T_{\wh\cp_{\Qc}}^2 \asymp \vert \bm\Sigma \bm\delta \vert_2^2 \Delta$.
Hence, the estimator $\wh\cp_{\Oc} = \wh\cp_{\Mc} \cdot \mathbb{I}_{\{C_{\Mc/\Qc} > 1 \}} + \wh\cp_{\Qc} \cdot \mathbb{I}_{\{C_{\Mc/\Qc} \le 1 \}}$ attains the desired rate of estimation.

{
\subsubsection{Proof of Corollary~\ref{cor:adapt:short}}

We refer to Lemma~2 of \cite{moen2024minimax} for the proof of the probabilistic statements below. 

Denoting the diagonal entries of $\bm\Sigma$ by $\Sigma_{ii}$, there exists $C \in (0, \infty)$ such that
\begin{align}
& \p\l( \max_{i \in [p]} \l\lvert\frac{\wh\Sigma_{ii}}{\Sigma_{ii}} - 1\r\rvert \le C\sqrt{\frac{\log(p)}{n}} \r) \ge 1 - p^{-1},  
\nn \\
& \p\l( \bigcap_{\ell \in \{1, 2\}} \bigcap_{t \in \mc T} \l\{ \l\vert \wh\Psi^2_{\ell, t} - \E(\wh\Psi^2_{\ell, t}) \r\vert \le \frac{\Psi^2}{2} \sqrt{\frac{\log\log(n)}{t}} \r\} \r) \ge 1 - C \log^{-1}(n).
\label{eq:psi:prob}
\end{align}
Next, let us set $\wh\Psi^2 = \max_{t \in \mc T} (\wh\Psi^2_{1, t} \vee \wh\Psi^2_{2, t})$.
We observe that
\begin{align*}
\E(\wh\Psi^2_{1, t}) &= \frac{t \wedge \cp}{t} \vert \bm\Sigma \bm\beta_0 \vert_2^2 + \frac{(t - \cp) \vee 0}{t} \vert \bm\Sigma \bm\beta_1 \vert_2^2 + \sigma^2 \le 3\Psi^2, \text{ \ and}
\\
\E(\wh\Psi^2_{2, t}) &= \frac{t \wedge (n - \cp)}{t} \vert \bm\Sigma \bm\beta_1 \vert_2^2 + \frac{(t - n + \cp) \vee 0}{t} \vert \bm\Sigma \bm\beta_0 \vert_2^2 + \sigma^2 \le 3\Psi^2,
\end{align*}
which leads to $\wh\Psi^2 \lesssim \Psi^2$ with probability tending to one thanks to~\eqref{eq:psi:prob} and the construction of $\mc T$.
WLOG, let us suppose that $\vert \bm\Sigma^{1/2}\bm\beta_1 \vert_2 \ge \vert \bm\Sigma^{1/2}\bm\beta_0 \vert_2$.
Then, for any $t \le \min(n/2, n - \cp)$, we have
$\wh\Psi^2_{2, t} \ge \Psi^2/2$ uniformly with probability tending to one.
Summarising these observations, we have
\begin{align*}
\p\l( M^{-1} \le \frac{\wh\Psi^2}{\Psi^2} \le M \r) \ge 1 - C\log^{-1}(n),
\end{align*}
where $M \in (1, \infty)$ is an absolute constant.

Finally, supposing that $p \le c_* n$ for a constant $c_* \in (0, \infty)$, there exists $C \in (0, \infty)$ that depends only on $c_*$, which gives
\begin{align*}
& \p\l( \Vert \wh{\bm\Sigma} - \bm\Sigma \Vert \le C \Vert \bm\Sigma \Vert \sqrt{\frac{p}{n}} \r) \ge 1 - e^{-p}.
\end{align*}
Summarising the above concludes the proof.

}

\subsubsection{Proof of Proposition~\ref{prop:refine}}

In what follows, the arguments hold deterministically, conditional on $\mc B$ in~\eqref{eq:set:b} and $\bar{\mc E} \cap \mc E$ in~\eqref{eq:set:e}, the latter of which implies the consistency of $\wh\cp_{\Oc}$ as derived in \Cref{prop:adapt}.
On $\mc B$, it holds with $\varpi_{\Rf} \ge \bar{C}_1\log(p \vee n)$ that
\begin{multline}
\label{eq:gamma:a}
\max_{\mbf a \in \{\bm\delta, \mbf e_\ell \in [p]\}} \max_{\varpi_{\Rf} < k < n - \varpi_{\Rf}} \sqrt{\frac{k (n - k)}{n}} \vert \bm\Sigma^{1/2} \mbf a \vert_2^{-1} \l\vert \mbf a^\top \l( \wh{\bm\gamma}_{k, n} - \wh{\bm\gamma}_{0, k} - \bm\gamma_{k, n} + \bm\gamma_{0, k} \r) \r\vert 
\\
\lesssim \Psi\sqrt{\log(p \vee n)}.
\end{multline}
To see this, consider $\varpi_{\Rf} < k \le \cp$. Then,
\begin{align*}
& \sqrt{\frac{k (n - k)}{n}} \l\vert \mbf a^\top \l( \wh{\bm\gamma}_{k, n} - \wh{\bm\gamma}_{0, k} - \bm\gamma_{k, n} + \bm\gamma_{0, k} \r) \r\vert
\\
= & \, \sqrt{\frac{k(n - k)}{n}} \Biggl\vert \frac{1}{n - k} \sum_{t = k + 1}^n \mbf a^\top \mbf x_t \mbf x_t^\top \bm\beta_1  - \frac{1}{n - k} \sum_{t = k + 1}^{\cp} \mbf a^\top \mbf x_t \mbf x_t^\top \bm\delta - \frac{1}{k} \sum_{t = 1}^k \mbf a^\top \mbf x_t \mbf x_t^\top \bm\beta_0 
\\
& \qquad \qquad \qquad \qquad - \frac{n - \cp}{n - k} \mbf a^\top \bm\Sigma \bm\delta  - \frac{1}{k} \sum_{t = 1}^k \mbf a^\top \mbf x_t \vep_t + \frac{1}{n - k} \sum_{t = k + 1}^n \mbf a^\top \mbf x_t \vep_t \Biggr\vert
\\
\le & \, \frac{1}{\sqrt{n - k}} \l\vert \sum_{t = k + 1}^n \mbf a^\top (\mbf x_t \mbf x_t^\top - \bm\Sigma) \bm\beta_1 \r\vert
+ \frac{1}{\sqrt{n - k}} \l\vert \sum_{t = k + 1}^\cp \mbf a^\top (\mbf x_t \mbf x_t^\top - \bm\Sigma) \bm\delta \r\vert
\\
& \, + \frac{1}{\sqrt{k}} \l\vert \sum_{t = 1}^k \mbf a^\top (\mbf x_t \mbf x_t^\top - \bm\Sigma) \bm\beta_0 \r\vert
+ \frac{1}{\sqrt{k}} \l\vert \sum_{t = 1}^k \mbf a^\top \mbf x_t \vep_t^\top \r\vert + \frac{1}{\sqrt{n - k}} \l\vert \sum_{t = k + 1}^n \mbf a^\top \mbf x_t \vep_t^\top \r\vert
\\
\le & \, 6 \bar{C}_2 \vert \bm\Sigma^{1/2} \mbf a \vert_2 \Psi \sqrt{\log(p \vee n)}.
\end{align*}
The case of $\cp < k < n - \varpi_{\Rf}$ is handled analogously.
Then by~\eqref{eq:gamma:a}, setting $\mbf a = \bm\delta$ gives
\begin{align}
\max_{\varpi_{\Rf} < k < n - \varpi_{\Rf}} \sqrt{\frac{k (n - k)}{n}} \l\vert \bm\delta^\top \l( \wh{\bm\gamma}_{k, n} - \wh{\bm\gamma}_{0, k} - \bm\gamma_{k, n} + \bm\gamma_{0, k} \r) \r\vert \lesssim \vert \bm\Sigma^{1/2} \bm\delta \vert_2 \Psi\sqrt{\log(p \vee n)}.
\nn
\end{align}
By the arguments adopted in the proof of \Cref{prop:lope}, on $\mc B$, we have
\begin{align*}
\l\vert \wh{\bm\delta} - \bm\delta(\wh\cp_{\Oc}) \r\vert_1 \lesssim \frac{\sigma_X \Psi \mathfrak{s}_\delta \sqrt{\log(p \vee n)}}{\underline{\sigma} \sqrt{\Delta}},
\end{align*}
and on $\bar{\mc E} \cap \mc E$, 
\begin{align*}
\l\vert \bm\delta(\wh\cp_{\Oc}) - \bm\delta \r\vert_1 \lesssim \frac{ \Psi}{\underline{\sigma}\sqrt{\Delta}} \min\l\{ \mathfrak{s}_\delta \sigma_X \sqrt{\log(p\log(n))}, \, { \Vert \bm\Sigma \Vert^{1/2}} \sqrt{\mathfrak{s}_\delta} \bigl(p\log\log(n)\bigr)^{1/4} \r\}
\end{align*}
which, { combined with that $\Vert \bm\Sigma \Vert^{1/2} \ge \sigma_X$,} leads to
% \footnote{
% since it is assumed that $p \ge n$ in Theorem~\ref{thm:qcscan} and hence here also, and because
% \begin{align*}
% \l\vert \bm\delta(\wh\cp) - \bm\delta \r\vert_1 \lesssim \frac{\Vert \bm\Sigma \Vert^{1/2} \Psi}{\underline{\sigma}\sqrt{\Delta}}
% \times
% \begin{cases}
% \mathfrak{s}_\delta \sqrt{\log(p \vee n)} & \text{if \ } \mathfrak{s}_\delta \lesssim \tfrac{\sqrt{p\log\log(n)}}{\log(p\log(n))},
% \\
% \sqrt{\mathfrak{s}_\delta} (p\log\log(n))^{1/4} & \text{otherwise}
% \end{cases} 
% \end{align*}
% regardless of $\mathfrak{s}_\delta$, the same conclusion follows
% }
\begin{align}
\label{eq:lope:lone}
\l\vert \wh{\bm\delta} - \bm\delta \r\vert_1 \lesssim \frac{\Vert \bm\Sigma \Vert^{1/2} \Psi \mathfrak{s}_\delta \sqrt{\log(p \vee n)}}{\underline{\sigma}\sqrt{\Delta}}.
% \times
% \begin{cases}
% \mathfrak{s}_\delta \sqrt{\log(p \vee n)} & \text{if \ } \mathfrak{s}_\delta \lesssim \tfrac{\sqrt{p\log\log(n)}}{\log(p\log(n))},
% \\
% \sqrt{\mathfrak{s}_\delta} (p\log\log(n))^{1/4} & \text{otherwise.}
% \end{cases}
\end{align}
Then, combining~\eqref{eq:gamma:a} (with $\mbf a \in \{\mbf e_\ell, \, \ell \in [p]\}$) and~\eqref{eq:lope:lone},
\begin{align}
& \max_{\varpi_{\Rf} < k < n - \varpi_{\Rf}} \sqrt{\frac{k (n - k)}{n}} \l\vert (\wh{\bm\delta} - \bm\delta)^\top \l( \wh{\bm\gamma}_{k, n} - \wh{\bm\gamma}_{0, k} - \bm\gamma_{k, n} + \bm\gamma_{0, k} \r) \r\vert 
\nn \\
\le & \, \vert \wh{\bm\delta} - \bm\delta \vert_1 \cdot \max_{\varpi_{\Rf} < k < n - \varpi_{\Rf}} \sqrt{\frac{k (n - k)}{n}} \l\vert \wh{\bm\gamma}_{k, n} - \wh{\bm\gamma}_{0, k} - \bm\gamma_{k, n} + \bm\gamma_{0, k} \r\vert_\infty 
\nn \\
\lesssim & \, \frac{\mathfrak{s}_\delta { \Vert \bm\Sigma \Vert^{1/2} \sigma_X} \Psi^2 \log(p \vee n)}{\underline{\sigma} \sqrt{\Delta}}.
\nn
\end{align}
Altogether, we have
\begin{align}
& \max_{\varpi_{\Rf} < k < n - \varpi_{\Rf}} \sqrt{\frac{k (n - k)}{n}} \l\vert \wh{\bm\delta}^\top \l( \wh{\bm\gamma}_{k, n} - \wh{\bm\gamma}_{0, k} - \bm\gamma_{k, n} + \bm\gamma_{0, k} \r) \r\vert 
\nn \\
\lesssim & \, \vert \bm\Sigma^{1/2} \bm\delta \vert_2 \Psi\sqrt{\log(p \vee n)} + \frac{\mathfrak{s}_\delta { \Vert \bm\Sigma \Vert^{1/2} \sigma_X} \Psi^2 \log(p \vee n)}{\underline{\sigma} \sqrt{\Delta}}.
\label{eq:lope:bound}
\end{align}
Also, by Proposition~\ref{prop:lope} with $C$ being the constant hidden in notation `$\lesssim$' therein, we have
\begin{align}
\bm\delta^\top \bm\Sigma \bm\delta - (\wh{\bm\delta} - \bm\delta)^\top \bm\Sigma \bm\delta 
& \ge \bm\delta^\top \bm\Sigma \bm\delta - \Vert \bm\Sigma \Vert^{1/2} \vert \bm\Sigma^{1/2} \bm\delta \vert_2 \vert \wh{\bm\delta} - \bm\delta \vert_2 
\nn \\
& \ge \bm\delta^\top \bm\Sigma \bm\delta - \frac{C\vert \bm\Sigma^{1/2} \bm\delta \vert_2 \Vert \bm\Sigma \Vert \Psi \sqrt{\mathfrak{s}_\delta \log(p \vee n)}}{\underline{\sigma} \sqrt{\Delta}}
\ge \frac{1}{2} \bm\delta^\top \bm\Sigma \bm\delta,
\label{eq:dsd}
\end{align}
as long as
$\vert \bm\Sigma^{1/2} \bm\delta \vert_2^2 \Delta \ge  4C^2 \underline{\sigma}^{-2} \Vert \bm\Sigma \Vert^2 \Psi^2 \mathfrak{s}_\delta \log(p \vee n)$, 
a condition readily met under~\eqref{eq:prop:refine} for large enough $\wt{c}_0$.
Given these observations, we first show that $\vert \wh{\cp}_{\Rf} - \cp \vert \le \Delta/2$. 
This follows from the fact that for any $\varpi_{\Rf} < k < \cp - \Delta/2$, we have from~\eqref{eq:lope:bound} and~\eqref{eq:dsd},
\begin{align*}
& \, \wh{\bm\delta}^\top \l[ \mbf X^\top \wt{\mbf Y}(k) - \mbf X^\top \wt{\mbf Y}(\cp) \r]
\\
\le & \, \l( \sqrt{\frac{k (n - k)}{n}} \frac{n - \cp}{n - k} - \sqrt{\frac{\cp(n - \cp)}{n}} \r) \wh{\bm\delta}^\top \bm\Sigma \bm\delta 
\\
& \, + 2 \max_{\varpi_{\Rf} < m < n - \varpi_{\Rf}}\sqrt{\frac{m (n - m)}{n}} \l\vert \wh{\bm\delta}^\top \l( \wh{\bm\gamma}_{m, n} - \wh{\bm\gamma}_{0, m} - \bm\gamma_{m, n} + \bm\gamma_{0, m} \r) \r\vert
\\
\le & \, - \frac{1}{2\sqrt{2}} \l( 1 - \sqrt{\frac{1}{3}} \r) \vert \bm\Sigma^{1/2} \bm\delta \vert_2^2 \sqrt{\Delta} 
\\
& \, + C' \l[ \vert \bm\Sigma^{1/2} \bm\delta \vert_2 \Psi\sqrt{\log(p \vee n)} + \frac{\mathfrak{s}_\delta \Vert \bm\Sigma \Vert \Psi^2 \log(p \vee n)}{\underline{\sigma} \sqrt{\Delta}} \r]
< 0,
\end{align*}
where $C'$ denotes the hidden constant factor in~\eqref{eq:lope:bound}.
The last inequality is met provided that
\begin{align}
\vert \bm\Sigma^{1/2} \bm\delta \vert_2^2 \Delta &\gtrsim \Psi^2 \log(p \vee n) \max\l\{1, \underline{\sigma}^{-1} \mathfrak{s}_\delta \Vert \bm\Sigma \Vert \r\},
\nn %\label{eq:refine:cond:two}
\end{align}
which in turn is satisfied under~\eqref{eq:prop:refine}.
The case where $\cp + \Delta/2 < k < n - \varpi_{\Rf}$, is handled analogously.

From here on, WLOG, suppose that $\wh{\cp}_{\Rf} \le \cp$ and also, focus on the case of $\vert \wh{\cp}_{\Rf} - \cp \vert \gtrsim \log(p \vee n)$ as otherwise, the claim holds from that $\vert \bm\Sigma^{1/2} \bm\delta \vert_2^{-2} \Psi^2 \gtrsim 1$.
By the construction of $\wh\cp_{\Rf}$, we have
\begin{align*}
& 0 \le \sqrt{\frac{\wh{\cp}_{\Rf}(n - \wh{\cp}_{\Rf})}{n}} \wh{\bm\delta}^\top ( \wh{\bm\gamma}_{\wh{\cp}_{\Rf}, n} - \wh{\bm\gamma}_{0, \wh{\cp}_{\Rf}} ) - \sqrt{\frac{\cp(n - \cp)}{n}} \wh{\bm\delta}^\top ( \wh{\bm\gamma}_{\cp, n} - \wh{\bm\gamma}_{0, \cp} )
\\
= &\, \wh{\bm\delta}^\top \Biggl( \sqrt{\frac{\wh{\cp}_{\Rf}}{n(n - \wh{\cp}_{\Rf})}} \sum_{t = \wh{\cp}_{\Rf} + 1}^n \mbf x_t \vep_t - \sqrt{\frac{n - \wh{\cp}_{\Rf}}{n\wh{\cp}_{\Rf}}} \sum_{t = 1}^{\wh{\cp}_{\Rf}} \mbf x_t \vep_t 
\\
& \qquad \qquad \qquad \qquad \qquad  - \sqrt{\frac{\cp}{n(n - \cp)}} \sum_{t = \cp + 1}^n \mbf x_t \vep_t + \sqrt{\frac{n - \cp}{n\cp}} \sum_{t = 1}^{\cp} \mbf x_t \vep_t \Biggr) 
\\
& + \wh{\bm\delta}^\top \l[ \sqrt{\frac{\wh{\cp}_{\Rf}}{n(n - \wh{\cp}_{\Rf})}} \sum_{t = \wh{\cp}_{\Rf} + 1}^{\cp} (\mbf x_t \mbf x_t^\top - \bm\Sigma) \bm\beta_0 + \sqrt{\frac{\wh{\cp}_{\Rf}}{n(n - \wh{\cp}_{\Rf})}} \sum_{t = \cp + 1}^{n} (\mbf x_t \mbf x_t^\top - \bm\Sigma) \bm\beta_1 \r.
\\
& \quad \quad \quad - \sqrt{\frac{n - \wh{\cp}_{\Rf}}{n\wh{\cp}_{\Rf}}} \sum_{t = 1}^{\wh{\cp}_{\Rf}} (\mbf x_t \mbf x_t^\top - \bm\Sigma) \bm\beta_0 - \sqrt{\frac{\cp}{n(n - \cp)}} \sum_{t = \cp + 1}^n (\mbf x_t \mbf x_t^\top - \bm\Sigma) \bm\beta_1 
\\
& \qquad \qquad \l. + \sqrt{\frac{n - \cp}{n\cp}} \sum_{t = 1}^{\cp} (\mbf x_t \mbf x_t^\top - \bm\Sigma) \bm\beta_0 \r]
- \l( \sqrt{\frac{\cp (n - \cp)}{n}} - \sqrt{\frac{\wh{\cp}_{\Rf} (n - \wh{\cp}_{\Rf})}{n}} \frac{n - \cp}{n - \wh{\cp}_{\Rf}} \r) \wh{\bm\delta}^\top \bm\Sigma \bm\delta 
\\
=: & \, U_1 + U_2 - U_3.
\end{align*}
Using the arguments leading to~\eqref{eq:lope:bound} in combination with those adopted in the proof of Theorem~2 of \cite{cho2024detection}, we have
\begin{align*}
\max(\vert U_1 \vert, \vert U_2 \vert) \lesssim 
\l( \vert \bm\Sigma^{1/2} \bm\delta \vert_2 \Psi \sqrt{\log(p \vee n)} + \frac{\Vert \bm\Sigma \Vert \Psi^2 \mathfrak{s}_\delta \log(p \vee n)}{\underline{\sigma}\sqrt{\Delta}} \r)
\sqrt{\frac{\vert \wh{\cp}_{\Rf} - \cp \vert}{\Delta}}
\end{align*}
while from that $\vert \wh\cp_{\Rf} - \cp \vert \le \Delta/2$, \eqref{eq:dsd} and Lemma 7 of \cite{wang2018high},
\begin{align*}
U_3 \ge \frac{1}{3\sqrt{6}} \cdot \frac{\vert \bm\Sigma^{1/2} \bm\delta \vert_2^2 \vert \wh{\cp}_{\Rf} - \cp \vert}{\sqrt{\Delta}}.
\end{align*}
Altogether, under~\eqref{eq:prop:refine}, we derive that
\begin{align*}
\vert \wh{\cp}_{\Rf} - \cp \vert &\lesssim \vert \bm\Sigma^{1/2} \bm\delta \vert_2^{-2} \Psi^2 \log(p \vee n) \max\l\{ 1, \, \frac{ \Vert \bm\Sigma \Vert^2 \Psi^2 \mathfrak{s}_\delta^2 \log(p \vee n)}{ \underline{\sigma}^2 \vert \bm\Sigma^{1/2} \bm\delta \vert_2^2 \Delta } \r\}
\\
&\lesssim \vert \bm\Sigma^{1/2} \bm\delta \vert_2^{-2} \Psi^2 \log(p \vee n).
\end{align*}

{

\subsubsection{Minimax localisation rate under strong signal strength}\label{ss:opt:loc}

The localisation rate obtained in Proposition~\ref{prop:refine} is minimax optimal up to a logarithmic factor. This can be seen by projecting the model onto the oracle direction $\bm\delta$ and considering the resulting univariate mean change problem for $\{\bm \delta^\top \mbf x_t Y_t\}_{t\in [n]}$, for which the optimal localisation rate for this problem is well understood, see e.g.~\citet{VeFrLeRe23}.
{Regarding the problem under~\eqref{eq:amoc} as that of detecting a change point in the mean of the $p$-variate sequence $\{\mbf x_t Y_t \}_{t \in [n]}$, such an oracle approach can be shown to attain the maximum high-dimensional efficiency as defined by \cite{aston2018high}.}

We now provide a formal argument using standard lower bound techniques based on
two hypotheses. Consider the single change point model under~\eqref{eq:amoc}, and define the model space
\begin{multline*}
\mc P_{\bm\Sigma, \sigma}^{\mathfrak{s}_{\delta}, n, p}(\tau) = \biggl\{ \p_{\cp, \bm\beta_0, \bm\beta_1} \;:\; \{(Y_t, \mbf x_t)\}_{t \in [n]} \sim \p_{\cp, \bm\beta_0, \bm\beta_1} \text{ \ such that \ } \mbf x_t \sim_{\iid} \mc N_p(\mbf 0, \bm\Sigma) \text{ and}\\
\text{independently \ }  \vep_t \sim_{\iid}\mc N (0,\sigma^2), \; 
Y_t = \mbf x_t^\top \bigl( \bm\beta_0 \cdot \mathbb{I}_{\{t \le \cp\}} + \bm\beta_1 \cdot \mathbb{I}_{\{t > \cp\}} \bigr) + \vep_t,
\\
\bm\delta = \bm\beta_1 - \bm\beta_0, \,
\Delta = \min(\cp, n - \cp), \,
 \, \Psi = \max(\vert \bm\Sigma^{1/2} \bm\beta_0 \vert_2, \vert \bm\Sigma^{1/2} \bm\beta_1 \vert_2, \sigma), 
\\
\lvert \bm\delta\rvert_0 = \mathfrak{s}_{\delta}
\text{ \ and \ } \bm\delta^\top \bm\Sigma \bm\delta \cdot \Delta \ge \tau \Psi^2 \mathfrak{s}_\delta^2 \log(p \vee n)\biggr\}
\end{multline*} 
for some constant $\tau > 0$, with $\bm\Sigma$ positive definite and $\sigma > 0$.

\begin{lem}
Let $\tau$ be either a sufficiently large constant  (in particular, $\tau \ge \wt{c}_0  \underline{\sigma}^{-2} \Vert \bm\Sigma \Vert^2$) or satisfy $\tau \to \infty$ (at an arbitrary rate) as $n\to\infty$. If
\begin{equation}\label{s:sz:n}
n \ge 16 \tau \mathfrak{s}_{\delta}^2\log(p \vee n),
\end{equation}
then 
\[
\inf_{\wh\cp} \sup_{\p \in \mc P_{\bm\Sigma, \sigma}^{\mathfrak{s}_{\delta}, n, p}(\tau)}
\p\!\left(\left|\wh\cp\bigl((Y_t, \mbf x_t)_{t\in[n]}\bigr) - \cp\right| 
\ge \frac{\Psi^2}{8\,\bm\delta^\top \bm\Sigma \bm\delta} \right)
\ge \frac{1}{4},
\]
where the infimum is taken over all measurable functions $\wh\cp:\mathbb{R}^{n \times (p + 1)}\to [n]$.
\end{lem}

\begin{proof}
Consider the two hypotheses
\[
\p^{(0)} = \p_{n/4, \bm \delta, \bm 0}\quad\text{and}\quad
\p^{(1)} =\p_{n/4+2\gamma,  \bm \delta, \bm 0},
\]
where $\bm\delta$ satisfies $\vert\bm\delta\vert_0 = \mathfrak{s}_{\delta}$ and $\bm\delta^\top \bm\Sigma \bm\delta  = \sigma^2/4$, with $\gamma \le n/4$. In this setting, $\Psi^2 = \sigma^2$. Condition~\eqref{s:sz:n} ensures that $\p^{(0)}, \p^{(1)} \in \mc P_{\bm\Sigma, \sigma}^{\mathfrak{s}_{\delta}, n, p}(\tau)$, and therefore
\[
\inf_{\wh\cp}\; \sup_{\p \in \mc P_{\bm\Sigma, \sigma}^{\mathfrak{s}_{\delta}, n, p}(\tau)}
\p\!\left(\bigl|\wh\theta -\theta(\p)\bigr| \ge \gamma\right)
\ge \frac{1-\sqrt{\alpha/2}}{2},
\]
where $\alpha = \chi^2(\p^{(1)},\p^{(0)})$. By Lemma~\ref{l:cross},
\[
\alpha + 1
= \left(\E_{\p_{\bm 0}}\!\left[\left(\frac{\mathrm{d}\p_{\bm \delta}}{\mathrm{d}\p_{\bm 0}}\right)^2\right]\right)^{2\gamma}
= 2^\gamma.
\]
Choosing $\gamma := { \Psi^2}/({8\,\bm\delta^\top \bm\Sigma \bm\delta}) = 1/2$ yields $\alpha = \sqrt{2} - 1 \le 1/2$, and hence
\[
\frac{1-\sqrt{\alpha/2}}{2} \ge \frac{1}{4},
\]
which completes the proof.
\end{proof}
}

\clearpage

\section{Additional simulation results}
\label[appendix]{app:comp}

In all figures, the vertical dashed lines mark where the sparsity level in consideration corresponds to $\sqrt p$. 
Also, the $x$-axis is shown on a log scale, and the $y$-axis on a squared root scale (to accommodate zero errors). 

\subsection{Advanced optimistic search vs.\ full grid searches}
\label[appendix]{app:grid}

\begin{figure}[h!t!b!p!]
    \centering
    \includegraphics[width=\linewidth]{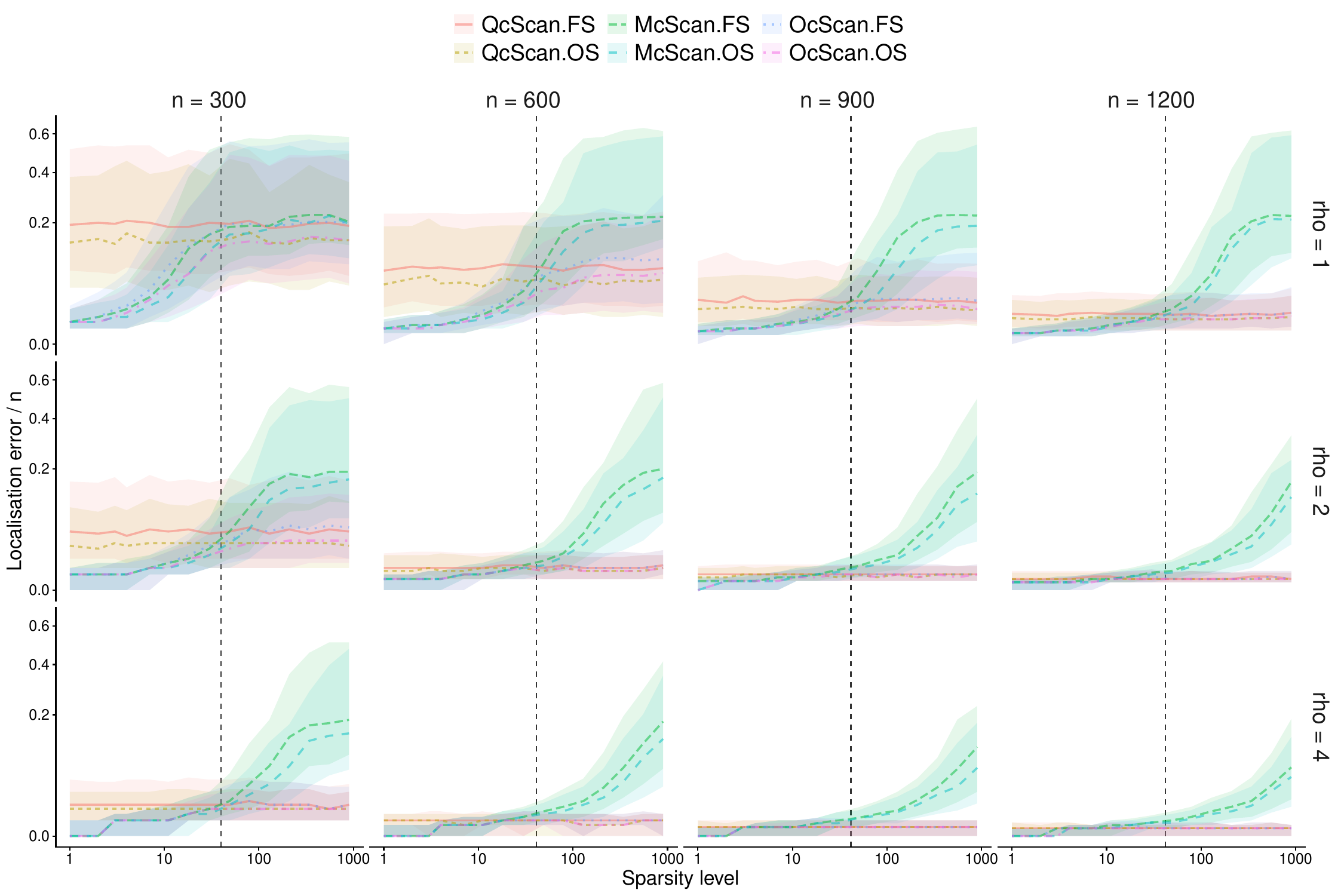}
    \caption{Estimation performance of McScan, QcScan and OcScan, using full grid search (FS) and advanced optimistic search (OS) in (M1).}
    \label{f:qos_quad_m1}
\end{figure}

\begin{figure}[h!t!b!p!]
    \centering
    \includegraphics[width=\linewidth]{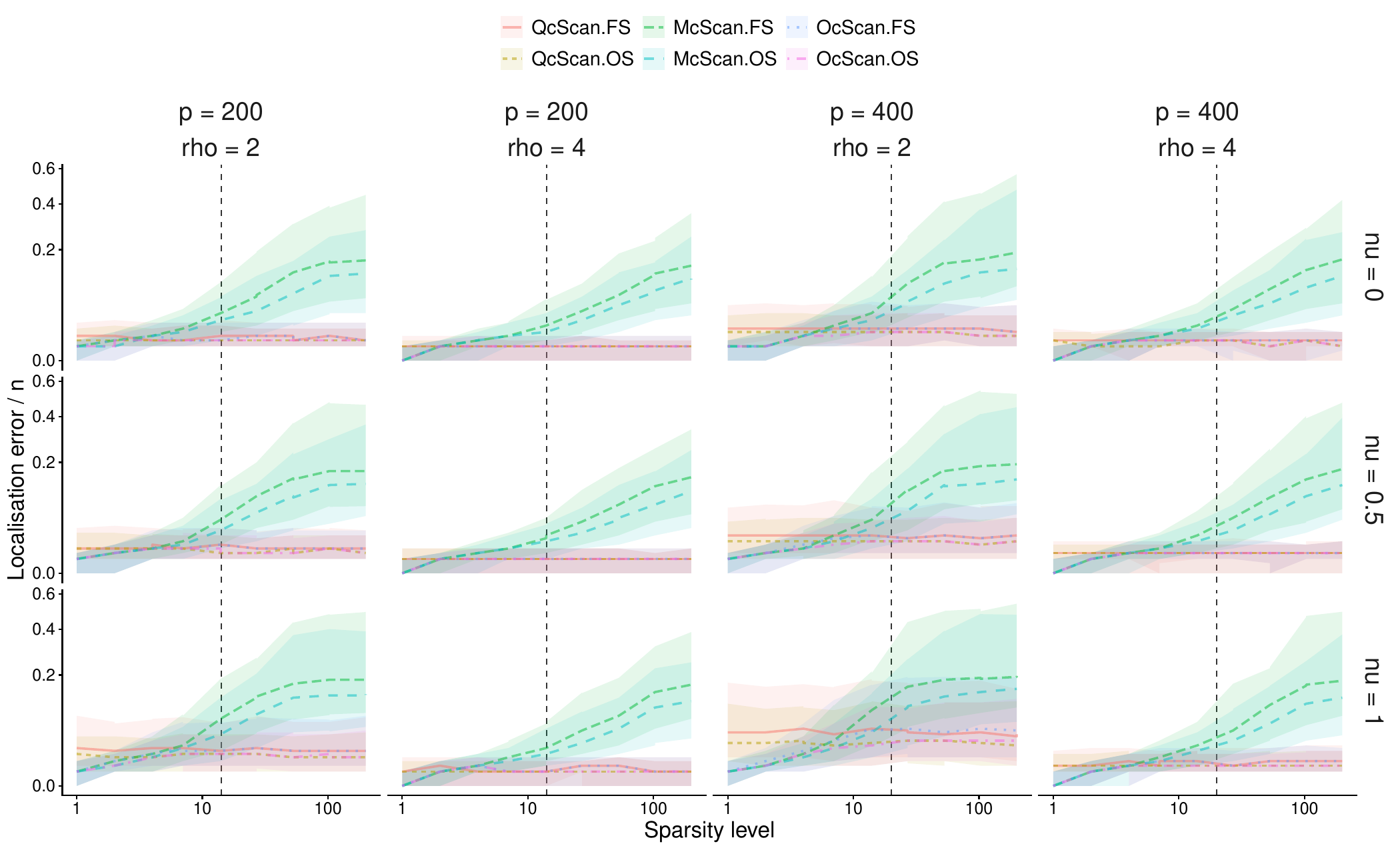}
    \caption{Estimation performance of McScan, QcScan and OcScan using full grid search (FS) and advanced optimistic search (OS) in (M2) with $\gamma = 0$.}
    \label{f:qos_quad_m2_ga0}
\end{figure}

\begin{figure}[h!t!b!p!]
    \centering
    \includegraphics[width=\linewidth]{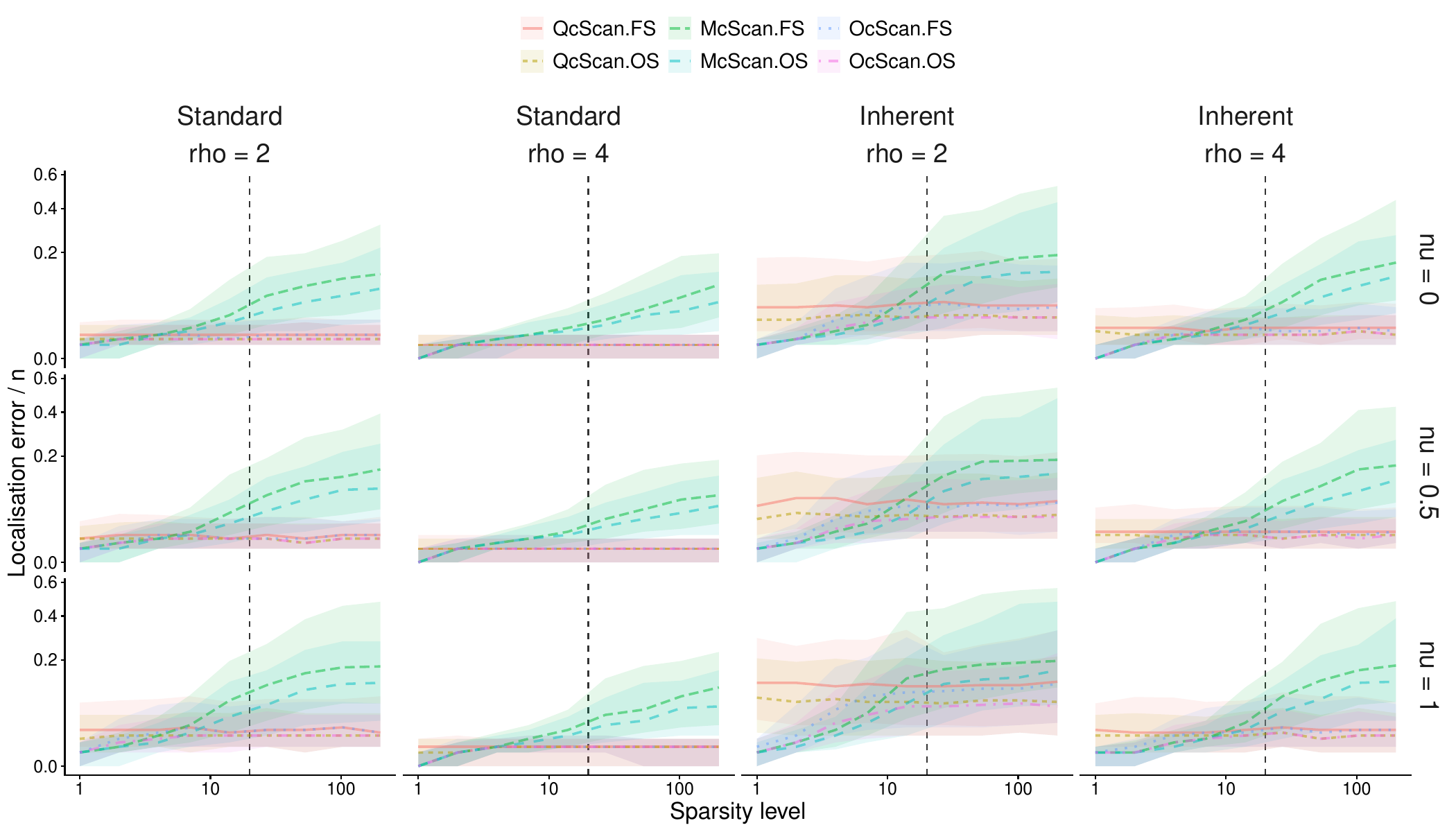}
    \caption{Estimation performance of McScan, QcScan and OcScan using full grid search (FS) and advanced optimistic search (OS) in (M2) with $\gamma = 0.6$ and $p = 400$.}
    \label{f:qos_quad_m2_ga0p6_p400}
\end{figure}

See \Cref{f:qos_quad_m1,f:qos_quad_m2_ga0,f:qos_quad_m2_ga0p6_p400} for the comparative performance of McSan, QcScan and OcScan applied with the advanced optimistic search in Algorithm~\ref{alg:os}, against with the full grid search under scenarios of (M1) and (M2) in Section~\ref{sec:sim}. The results indicate only minor differences in performance, with the variants using advanced optimistic searches showing a marginal but persistent advantage overall, possibly because the change point lies exactly on the dyadic grid explored by the advanced optimistic search.

\clearpage

\subsection{McScan vs.\ its refinement}
\label[appendix]{app:refine}

\begin{figure}[h!t!b!p!]
    \centering
    \includegraphics[width=\linewidth]{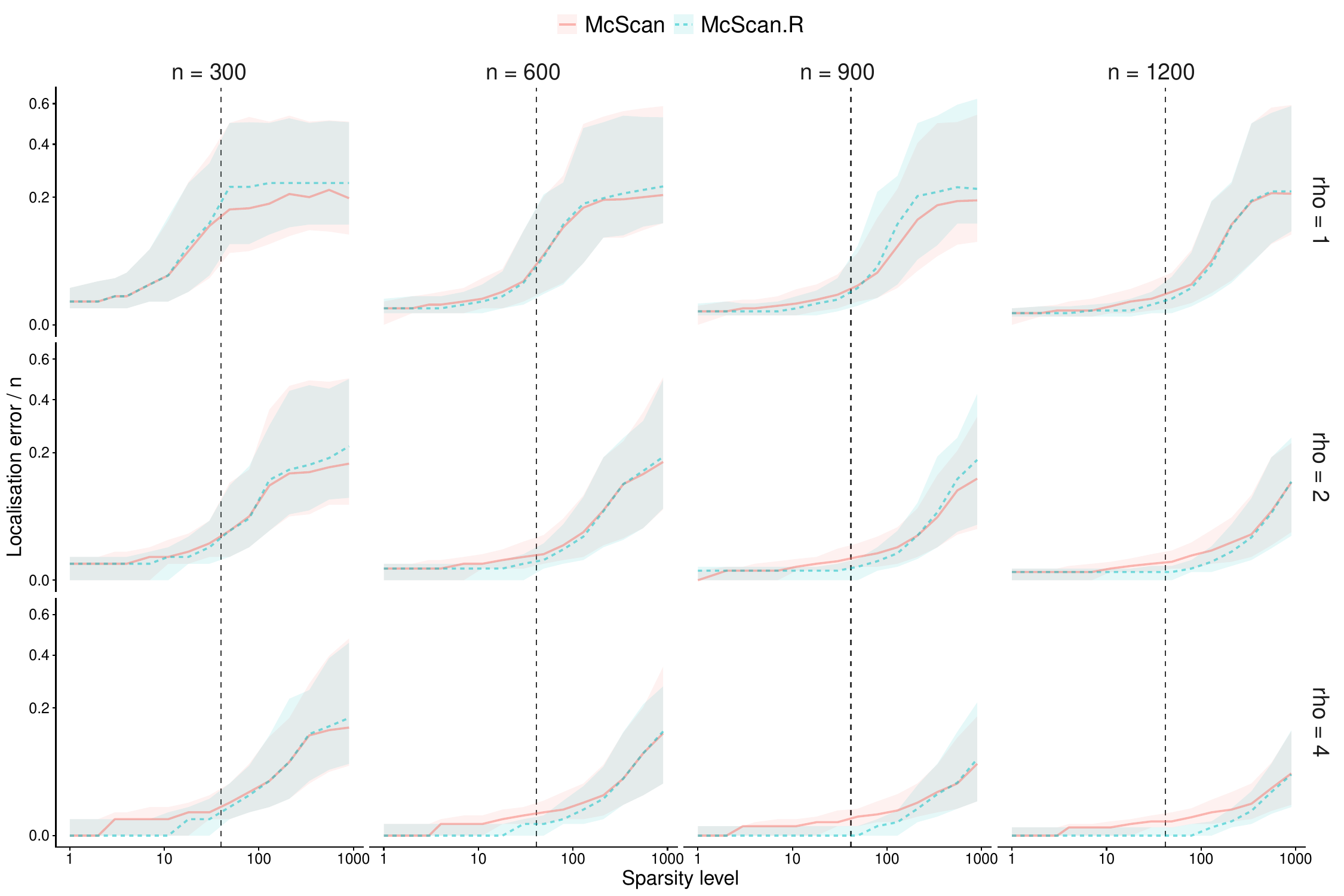}
    \caption{Estimation performance of McScan and McScan.R in (M1).}
    \label{f:mcs_rms_m1}
\end{figure}

\begin{figure}[h!t!b!p!]
    \centering
    \includegraphics[width=\linewidth]{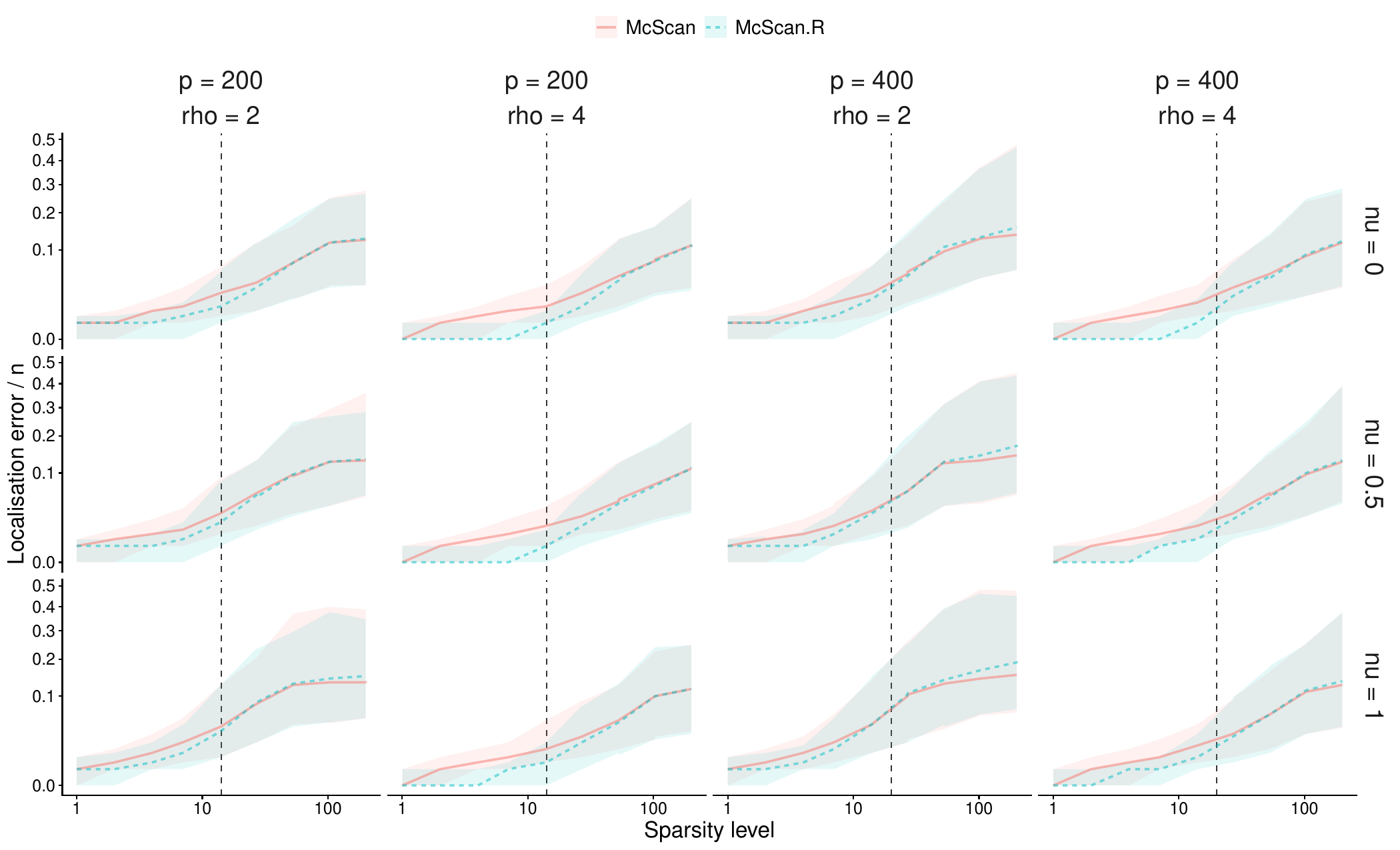}
    \caption{Estimation performance of McScan and McScan.R in (M2) with $\gamma = 0$.}
    \label{f:mss_rms_m2_ga0}
\end{figure}

\begin{figure}[h!t!b!p!]
    \centering
    \includegraphics[width=\linewidth]{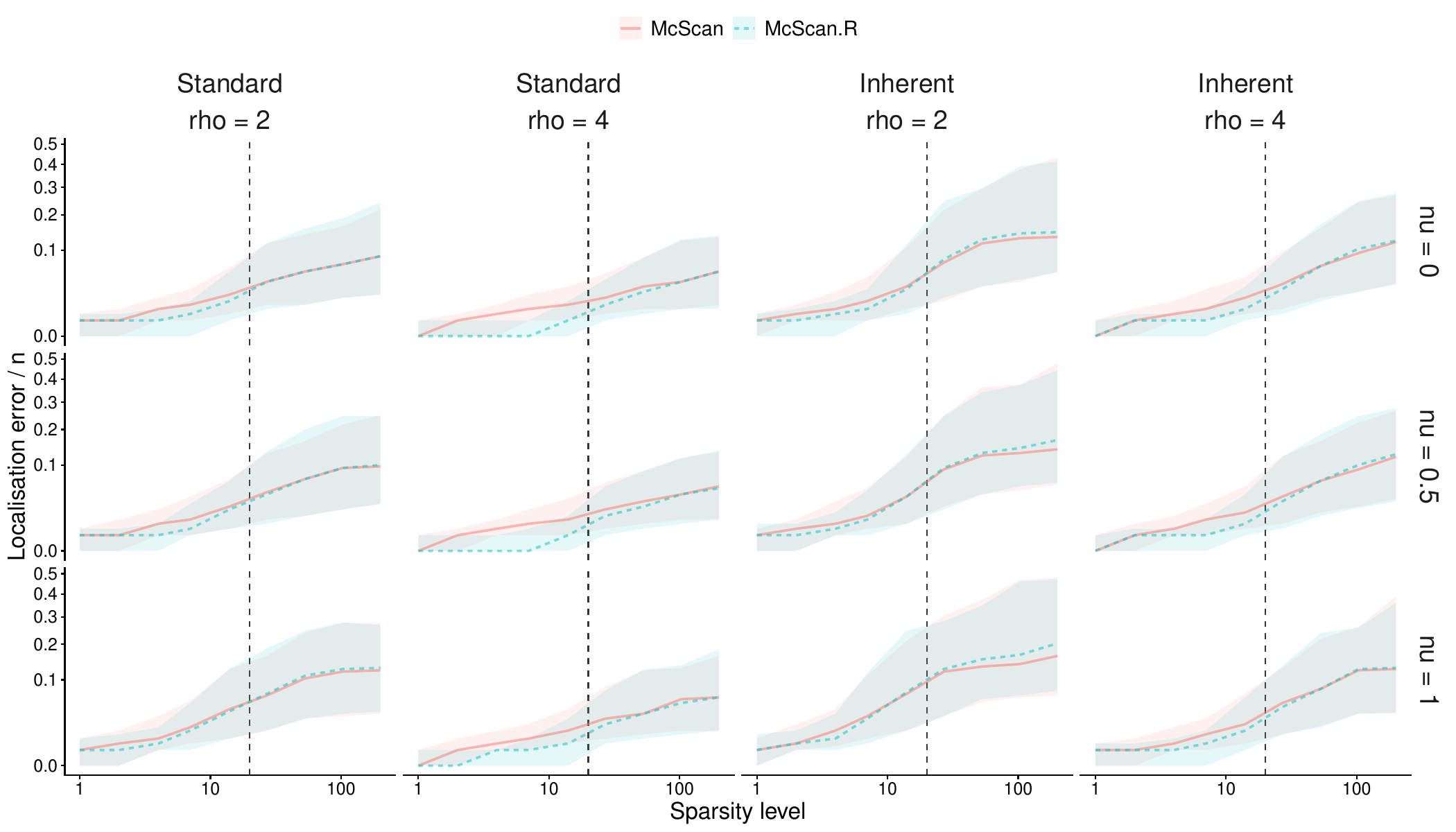}
    \caption{Estimation performance of McScan and McScan.R in (M2) with $\gamma = 0.6$ and $p = 400$.}
    \label{f:mss_rms_m2_ga0p6_p400}
\end{figure}

We compare McScan with its refinement referred to as McScan.R.
We first estimate the change point by McScan with which the differential parameter is estimated as in~\eqref{eq:lope}. Then we re-estimate the change point by scanning the covariance projected on the direction of estimated differential parameter.
See \Cref{f:mcs_rms_m1,f:mss_rms_m2_ga0,f:mss_rms_m2_ga0p6_p400} for the comparative performance of McScan and its refinement McScan.R in scenarios of (M1) and (M2) in Section~\ref{sec:sim}. 
We see the enhanced performance of McScan.R in the sparse regime (both in terms of $\mathfrak{s}$ and $\mathfrak{s}_\delta$) and when the size of change is large (with larger~$\rho$), which is in agreement with the stronger condition made in Proposition~\ref{prop:refine}.

\clearpage

\subsection{OcScan vs.\ its variants}
\begin{figure}[h!t!b!p!]
    \centering
    \includegraphics[width=\linewidth]{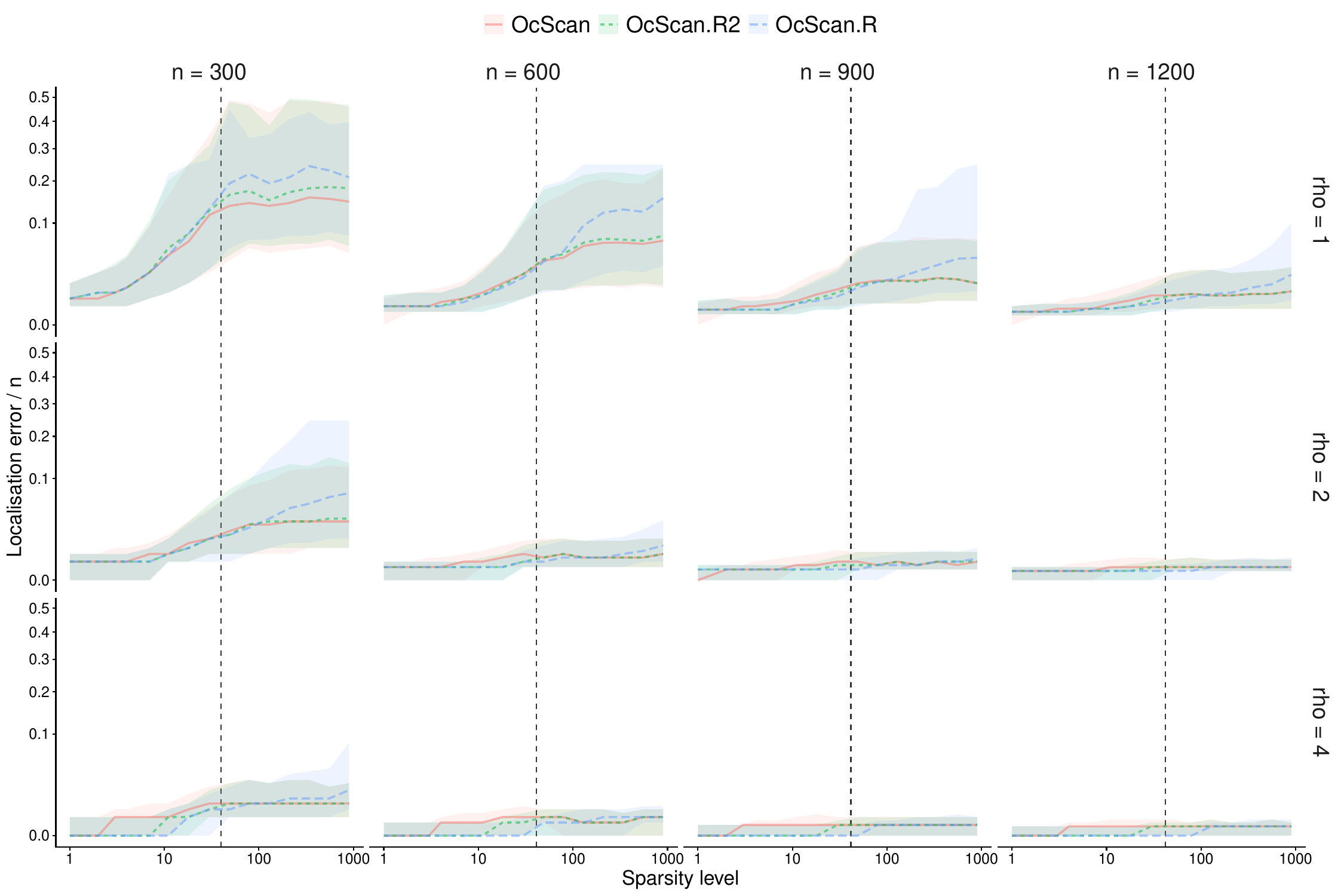}
    \caption{Estimation performance of OcScan, OcScan.R and OcScan.R2 in (M1).}
    \label{f:refine_adapt_m1}
\end{figure}

\begin{figure}[h!t!b!p!]
    \centering
    \includegraphics[width=\linewidth]{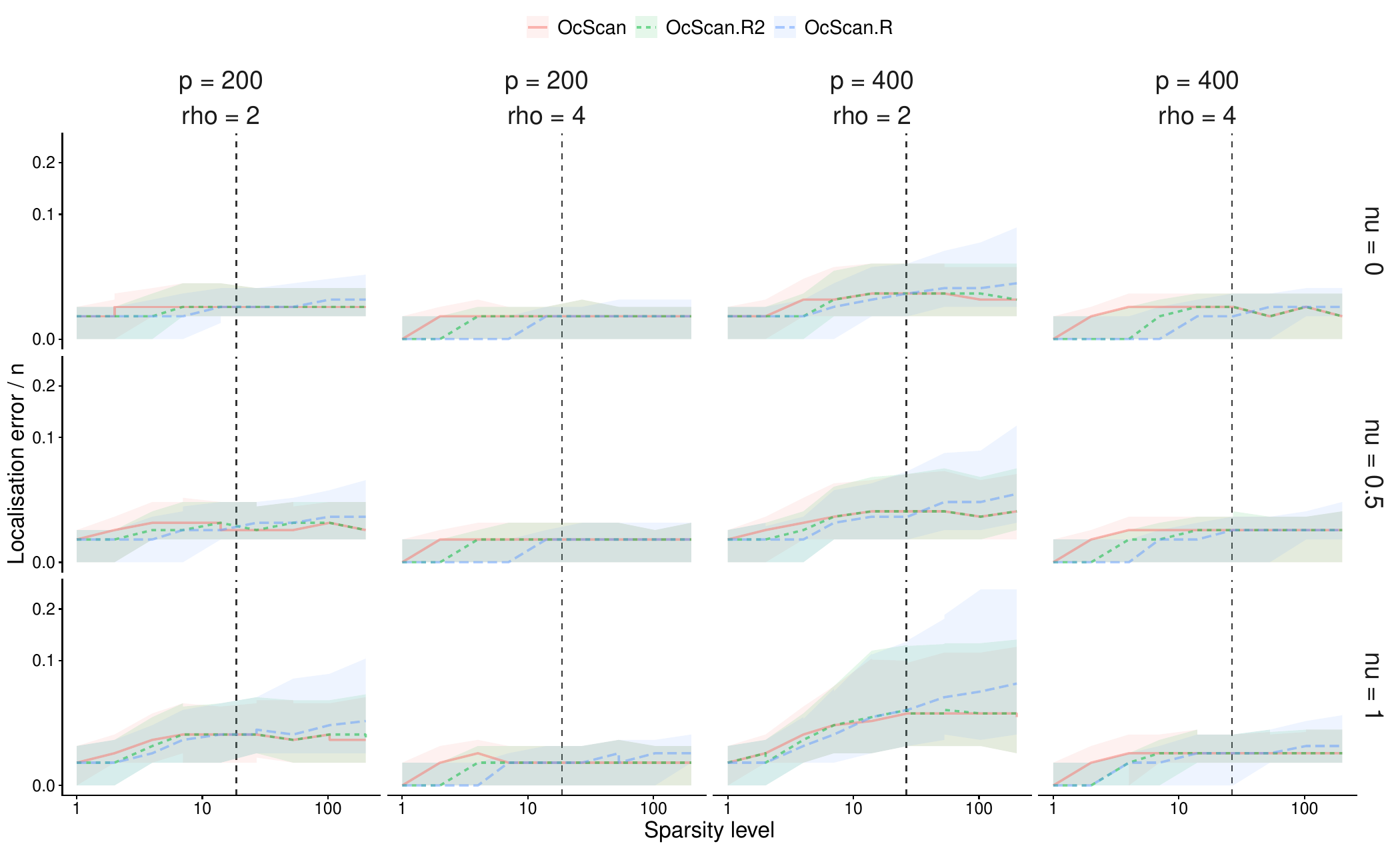}
    \caption{Estimation performance of OcScan, OcScan.R and OcScan.R2 in (M2) with $\gamma = 0$.}
    \label{f:refine_adapt_m2_ga0}
\end{figure}

\begin{figure}[h!t!b!p!]
    \centering
    \includegraphics[width=\linewidth]{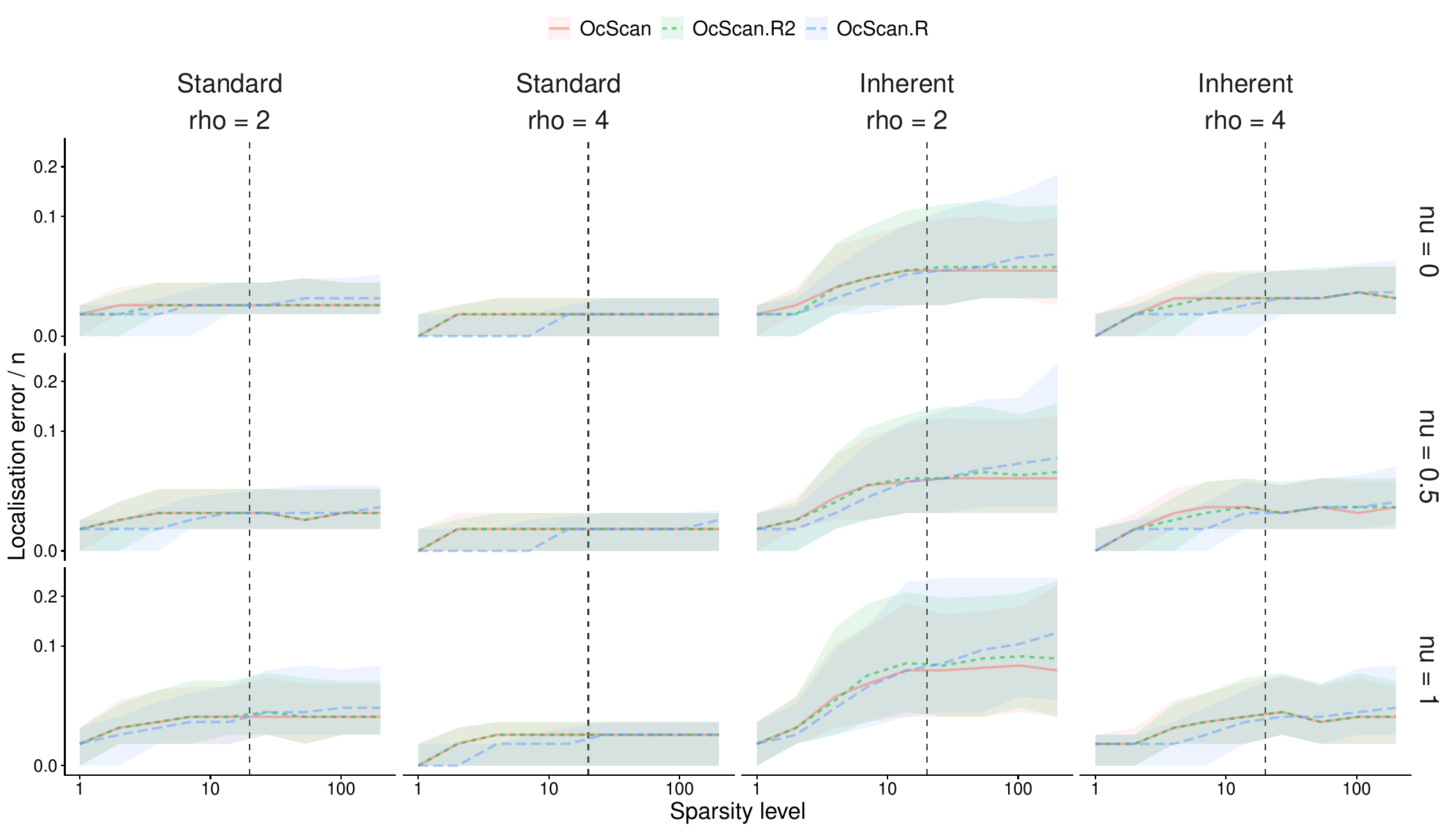}
    \caption{Estimation performance of OcScan, OcScan.R and OcScan.R2 in (M2) with $\gamma=0.6$ and $p = 400$.}
    \label{f:refine_adapt_m2_ga0p6_p400}
\end{figure}

\label[appendix]{app:ocscan}
We consider two variants of OcScan that incorporate a further refinement step. In addition to OcScan.R introduced in Section~\ref{sec:refine}, we investigate OcScan.R2 where, if the estimated change point of OcScan is closer to the one by McScan, we report the change point estimated by McScan.R (\Cref{app:refine}), and otherwise we report the one by OcScan. 
% \begin{description}[wide, leftmargin=0pt]
% \item[{OcScan.R}.] We proceed the same way as McScan.R, but with $\wh\cp$ in~\eqref{eq:lope} obtained by OcScan instead of McScan. 
% \item[{OcScan.R2}.] If the estimated change point of OcScan is closer to the one by McScan, we report the change point estimated by McScan.R, and otherwise we report the  one by OcScan. 
% \end{description}
See \Cref{f:refine_adapt_m1,f:refine_adapt_m2_ga0,f:refine_adapt_m2_ga0p6_p400} for the comparative performance of OcScan and its two variants in scenarios of (M1) and (M2) in Section~\ref{sec:sim}. 
We observe that the refinement step improves upon the performance of OcScan, particularly in the sparse regime. Among the three methods, OcScan.R is preferred, although it performs slightly worse in the dense regime.

\clearpage

\subsection{Additional simulations results under (M1)}
\label[appendix]{app:m1}

\begin{figure}[ht]
\centering
\includegraphics[width=\linewidth]{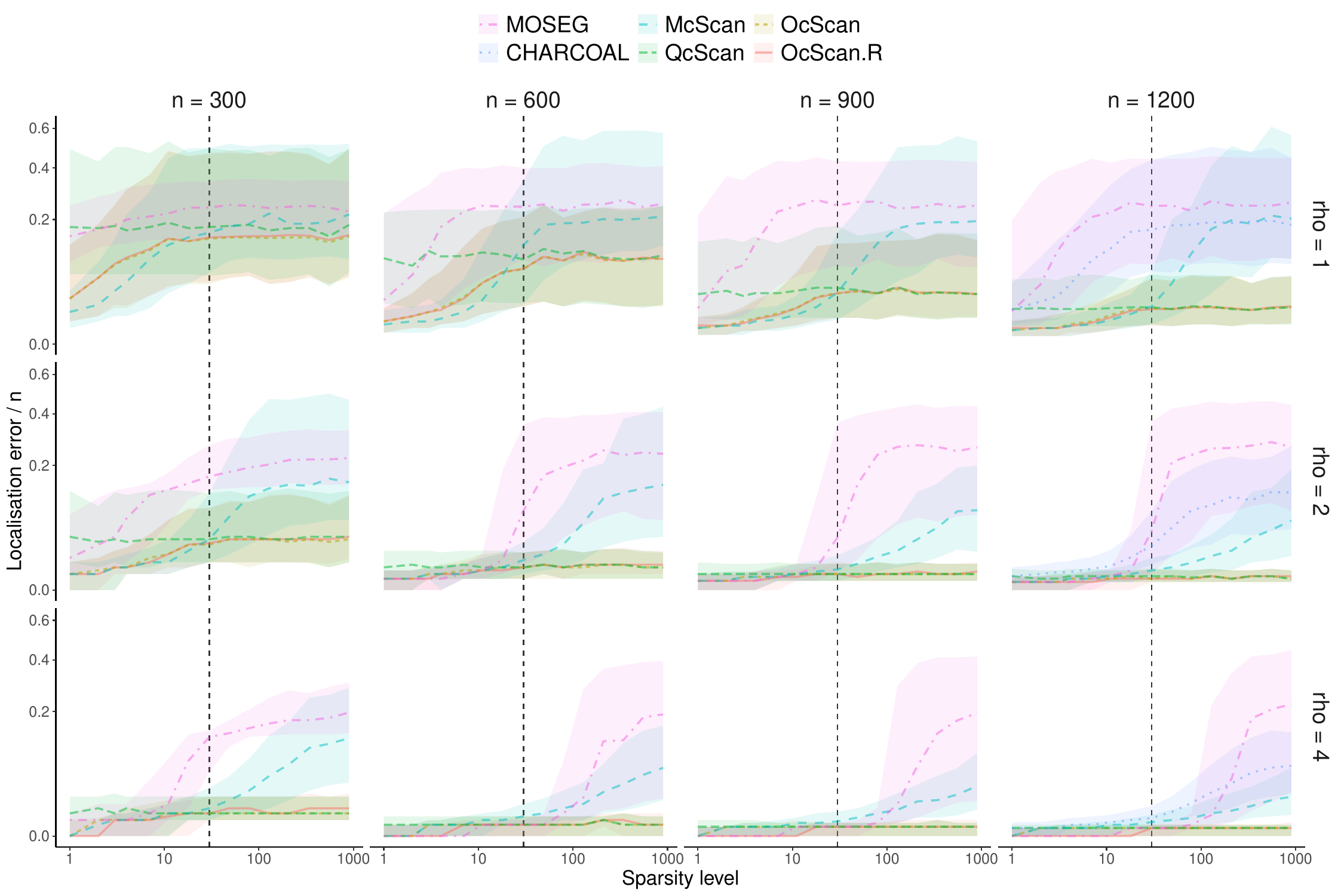}
\caption{Estimation performance of MOSEG, CHARCOAL, McScan, QcScan, OcScan and OcScan.R in modified (M1) with rank deficiency $r = p/2$; CHARCOAL is only applicable in the last column where $n = 1200 > p = 900$. 
The $x$-axis denotes $\mathfrak{s}_\delta = \vert \bm\delta \vert_0$. 
In each scenario, the results are based on 1000 repetitions, with median error curves shown alongside shaded regions representing the interquartile range. 
The vertical dashed lines mark where $\mathfrak{s}_\delta = \sqrt{p}$. 
The $x$-axis is shown on a log scale, and the $y$-axis on a squared root scale.} \label{f:rank450}
\end{figure}

\begin{figure}[ht]
\centering
\includegraphics[width=\linewidth]{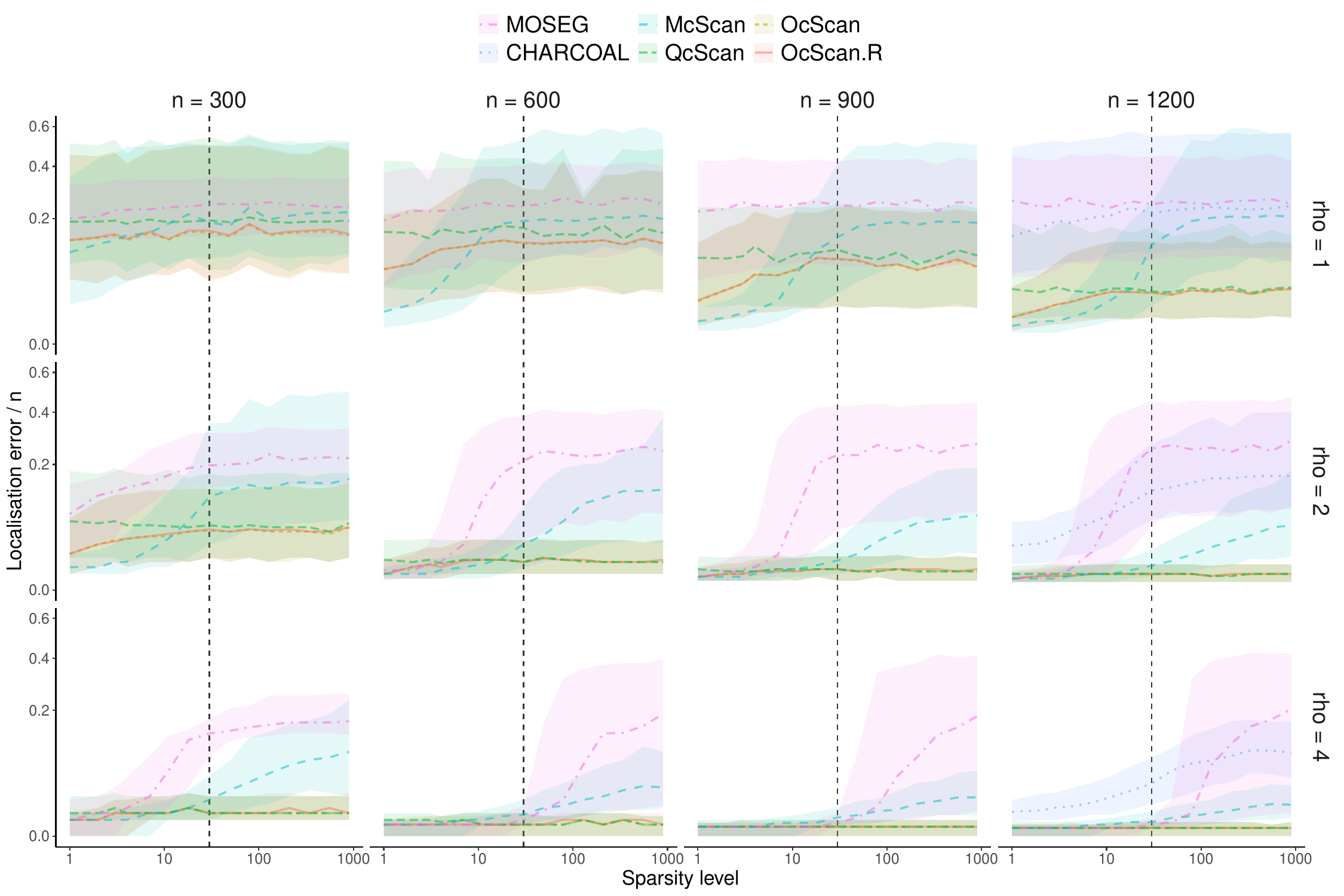}
\caption{Estimation performance of MOSEG, CHARCOAL, McScan, QcScan, OcScan and OcScan.R in modified (M1) with rank deficiency $r = p/4$; CHARCOAL is only applicable in the last column where $n = 1200 > p = 900$. 
The $x$-axis denotes $\mathfrak{s}_\delta = \vert \bm\delta \vert_0$. 
In each scenario, the results are based on 1000 repetitions, with median error curves shown alongside shaded regions representing the interquartile range. 
The vertical dashed lines mark where $\mathfrak{s}_\delta = \sqrt{p}$. 
The $x$-axis is shown on a log scale, and the $y$-axis on a squared root scale.} \label{f:rank225}
\end{figure}

Additional simulation results under the rank deficient variant of scenario (M1) in \Cref{ss:loc:err} are given in \Cref{f:rank450,f:rank225}. Together with \Cref{f:iderr}, these figures demonstrate that the performance advantage of the proposed covariance scanning methods over the competitors (MOSEG and CHARCOAL) becomes more pronounced as the rank $r$ decreases.

\clearpage

\subsection{Additional simulations results under (M2)}
\label[appendix]{app:m2}

We present further simulation results obtained under the scenario (M2) in Section~\ref{ss:loc:err}. 

\begin{figure}[h!t!b!p!]
\centering
\includegraphics[width=\textwidth]{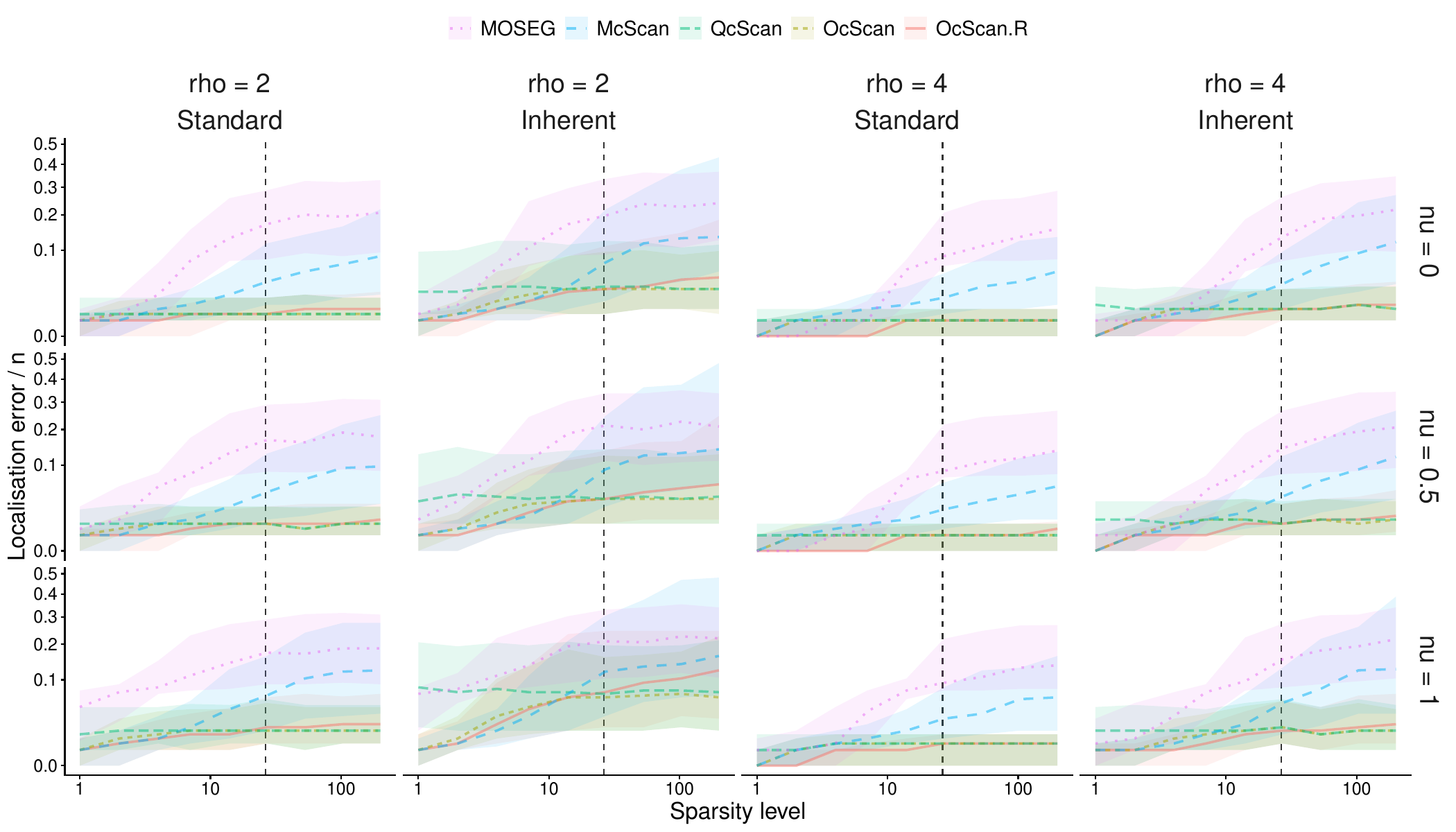}
\caption{Estimation performance of MOSEG, CHARCOAL, McScan, QcScan, OcScan and OcScan.R in (M2) with $\gamma = 0.6$ and $p = 400$. CHARCOAL is not applicable as $p = 400 > n = 300$.
In the left two columns, the $x$-axis denotes the standard sparsity $\mathfrak{s}_\delta = \vert \bm\delta \vert_0$ while in the right two columns, it is the inherent sparsity $\mathfrak{s} = \vert \bm\Sigma^{1/2}\bm\delta \vert_0$. 
In each scenario, the results are based on 1000 repetitions, with median error curves shown alongside shaded regions representing the interquartile range. 
The vertical dashed lines mark where $\mathfrak{s}_\delta = \sqrt{p}$ or $\mathfrak{s} = \sqrt{p}$. 
The $x$-axis is shown on a log scale, and the $y$-axis on a squared root scale. 
\label{f:toep_0p6_p400}}
\end{figure}

\begin{figure}[h!t!b!p!]
\centering
\includegraphics[width=\textwidth]{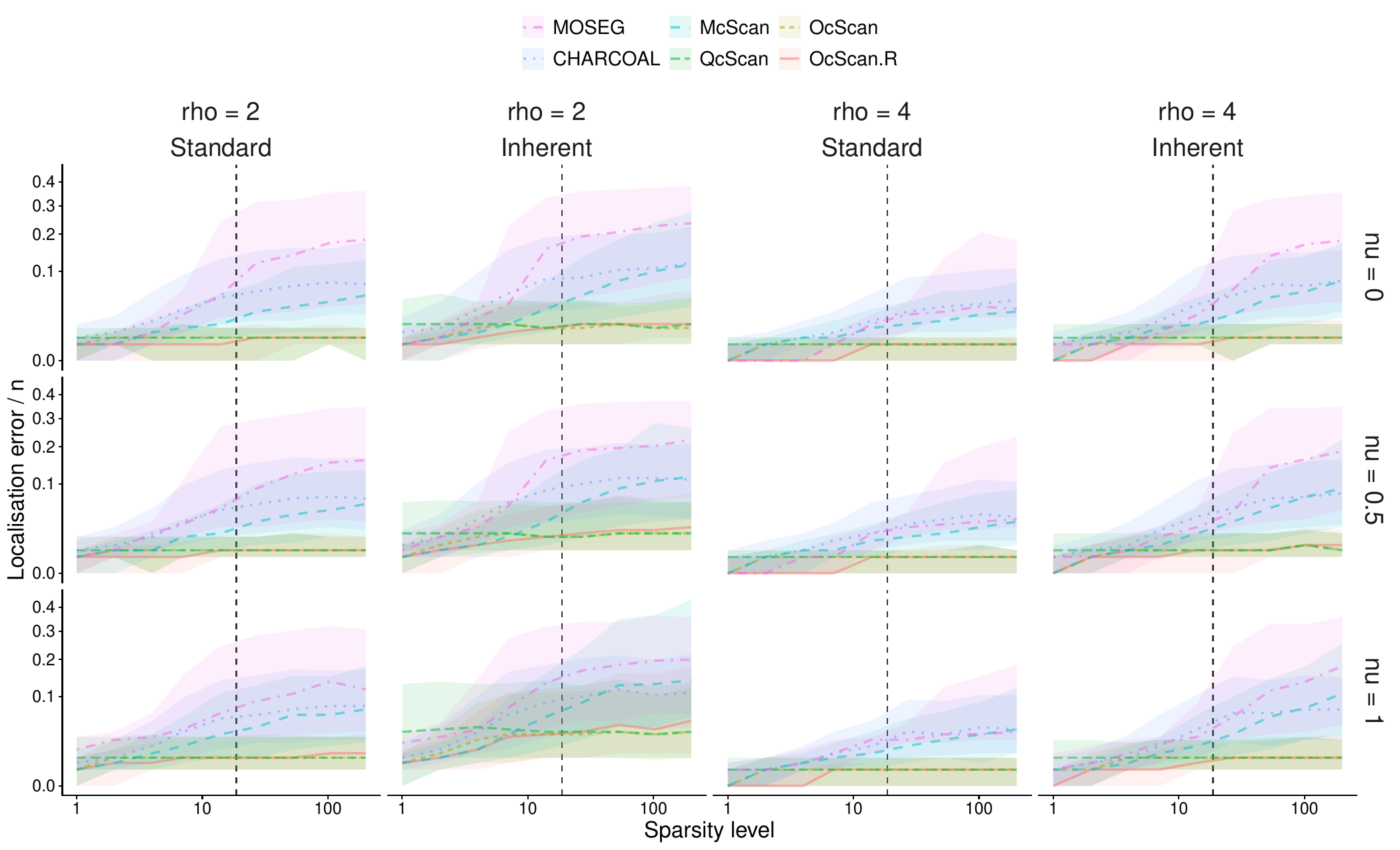}
\caption{
Estimation performance of MOSEG, CHARCOAL, McScan, QcScan, OcScan and OcScan.R in (M2) with $\gamma = -0.6$ and $p = 200$. 
In the left two columns, the $x$-axis denotes the standard sparsity $\mathfrak{s}_\delta = \vert \bm\delta \vert_0$ while in the right two columns, it is the inherent sparsity $\mathfrak{s} = \vert \bm\Sigma^{1/2}\bm\delta \vert_0$. 
In each scenario, the results are based on 1000 repetitions, with median error curves shown alongside shaded regions representing the interquartile range. 
The vertical dashed lines mark where $\mathfrak{s}_\delta = \sqrt{p}$ or $\mathfrak{s} = \sqrt{p}$. 
The $x$-axis is shown on a log scale, and the $y$-axis on a squared root scale. 
\label{f:toep_n0p6_p200}}
\end{figure}

\begin{figure}[h!t!b!p!]
\centering
\includegraphics[width=\textwidth]{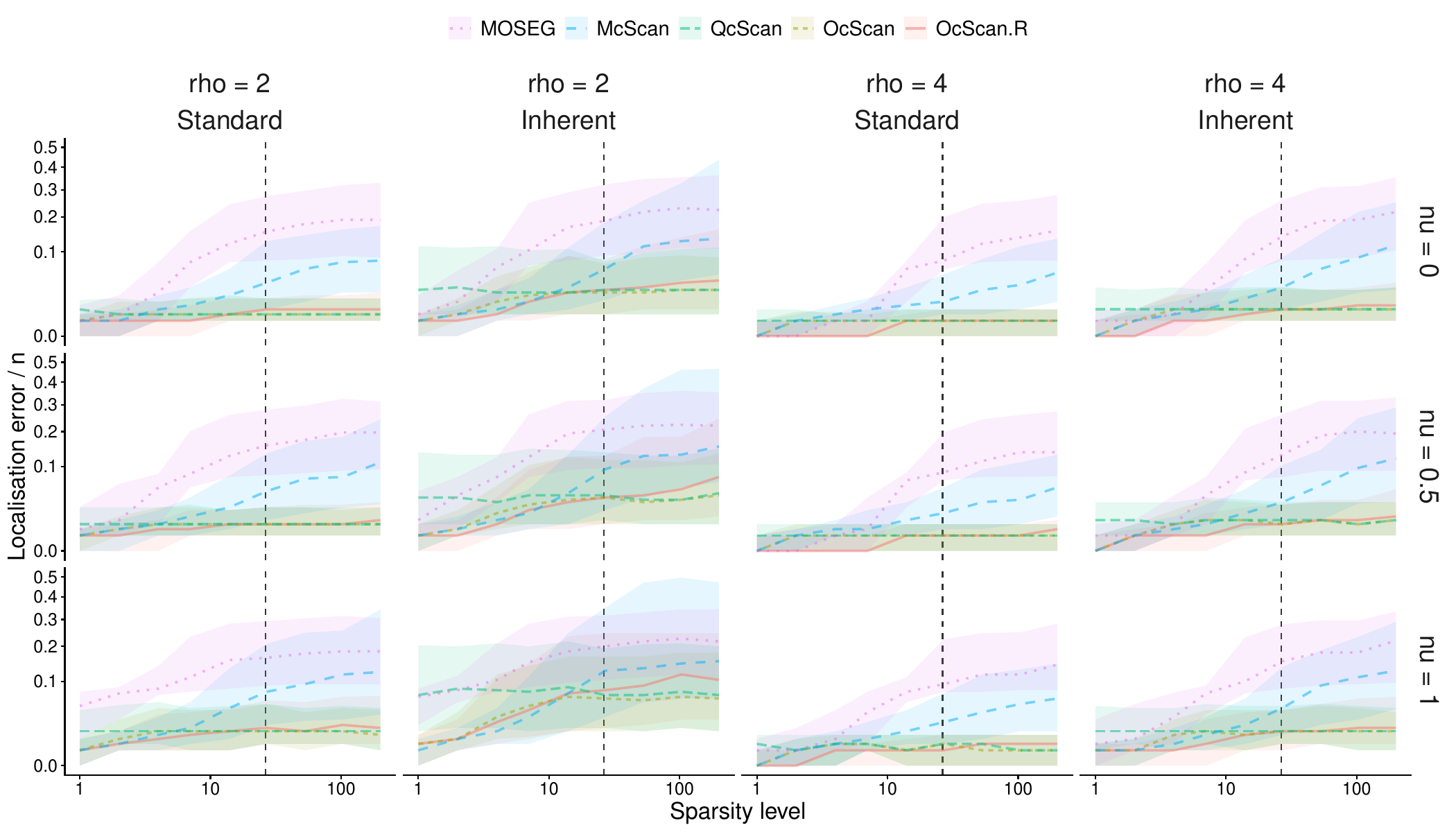}
\caption{Estimation performance of MOSEG, CHARCOAL, McScan, QcScan, OcScan and OcScan.R in (M2) with $\gamma = -0.6$ and $p = 400$. CHARCOAL is not applicable as $p = 400 > n = 300$.
In the left two columns, the $x$-axis denotes the standard sparsity $\mathfrak{s}_\delta = \vert \bm\delta \vert_0$ while in the right two columns, it is the inherent sparsity $\mathfrak{s} = \vert \bm\Sigma^{1/2}\bm\delta \vert_0$. 
In each scenario, the results are based on 1000 repetitions, with median error curves shown alongside shaded regions representing the interquartile range. 
The vertical dashed lines mark where $\mathfrak{s}_\delta = \sqrt{p}$ or $\mathfrak{s} = \sqrt{p}$. 
The $x$-axis is shown on a log scale, and the $y$-axis on a squared root scale. 
 \label{f:toep_n0p6_p400}}
\end{figure}

\clearpage

{
\subsection{Choice of thresholds}
\label[appendix]{app:thd}

We investigate the choice of constants in the thresholds 
$\zeta_{\Mc} = \bar{c} \wh\sigma_X \wh\Psi \sqrt{\log(p \log(n))}$ and 
$\zeta_{\Qc} = c \Vert \wh{\bm\Sigma} \Vert \wh\Psi^2 \sqrt{p\log\log(n)}$, 
defined in \Cref{cor:adapt:short}, for McScan and QcScan, respectively. 
Specifically, we consider model (M2) in \Cref{ss:loc:err} with $\nu = 0$, 
$n = 300$ and $p = 200$, under varying levels of intrinsic sparsity 
$\mathfrak{s} = \vert \bm\Sigma^{1/2} \bm\delta \vert_0 = \vert \bm\Sigma^{1/2} \bm\beta_0 \vert_0 \in \{1, 4, 14, 53, 200\}$. 
The simulation results are summarised in 
\Cref{f:thd:mcscan,f:thd:qcscan}. 
Based on these findings, we recommend the choices 
$\bar{c} = 1.3$ and $c = 0.7$.

\begin{figure}[h!t!b!p!]
\centering
\includegraphics[width=\textwidth]{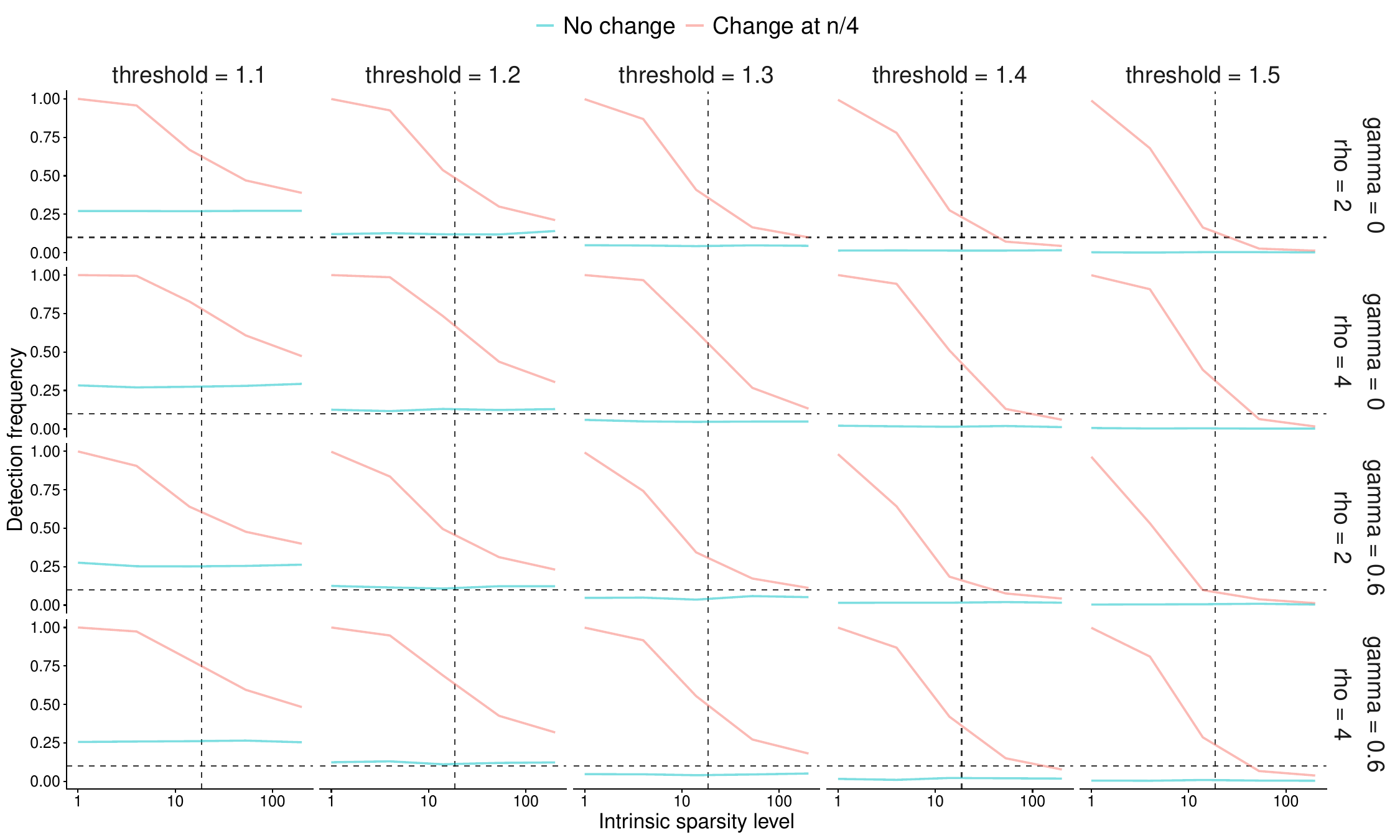}
\caption{Detection frequencies of McScan in (M2) with $\nu = 0$, $n = 300$ and $p = 200$, obtained with various threshold choices $\zeta_{\Mc} = \bar{c} \wh\sigma_X \wh\Psi \sqrt{\log(p \log(n))}$ where 
$\bar{c}\in \{1.1, 1.2, 1.3, 1.4, 1.5\}$ (left to right), when a single change point is present (red, $(q, \cp) = (1, n/4)$) and no change point exists (blue, $(q, \cp) = (0, n)$). Detection frequencies are computed 
from 1000 repetitions. The $x$-axis, displayed on a log scale, represents 
the intrinsic sparsity $\mathfrak{s} = \vert \bm\Sigma^{1/2}\bm\delta \vert_0$. 
The vertical dashed lines indicate where $\mathfrak{s} = \sqrt{p\log\log(n)}$, 
while the horizontal dashed lines mark the $10\%$ level.}
\label{f:thd:mcscan}
\end{figure}

\begin{figure}[h!t!b!p!]
\centering
\includegraphics[width=\textwidth]{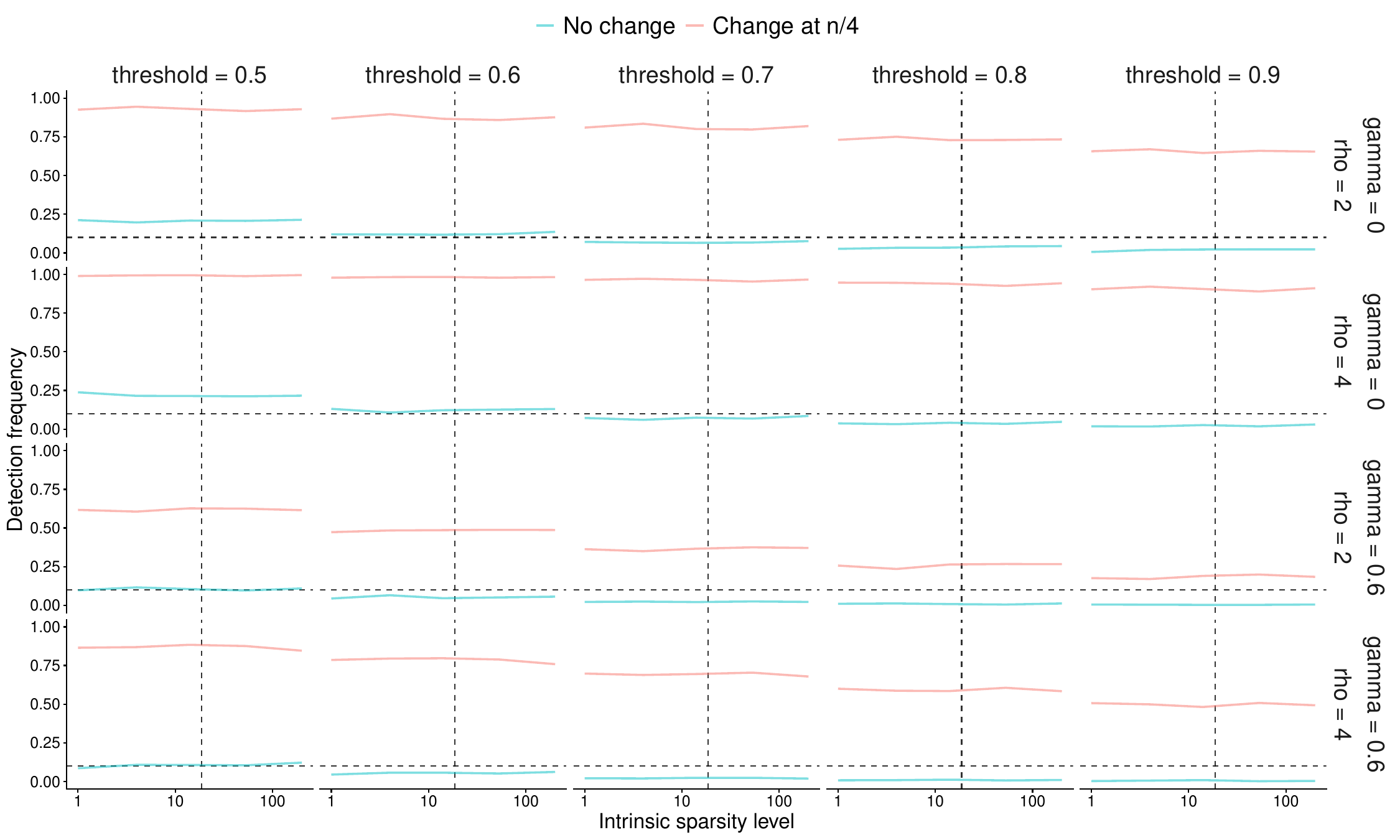}
\caption{Detection frequencies of QcScan in (M2) with $\nu = 0$, $n = 300$ and $p = 200$, obtained with various threshold choices $\zeta_{\Qc} = c \Vert \wh{\bm\Sigma} \Vert \wh\Psi^2 \sqrt{p\log\log(n)}$ where 
$c \in \{0.5, 0.6, 0.7, 0.8, 0.9\}$ (left to right), when a single change point is present (red, $(q, \cp) = (1, n/4)$) and no change point exists (blue, $(q, \cp) = (0, n)$). Detection frequencies are computed 
from 1000 repetitions. The $x$-axis, displayed on a log scale, represents the intrinsic sparsity $\mathfrak{s} = \vert \bm\Sigma^{1/2}\bm\delta \vert_0$. 
The vertical dashed lines indicate where $\mathfrak{s} = \sqrt{p\log\log(n)}$, 
and the horizontal dashed lines mark the $10\%$ level.}
\label{f:thd:qcscan}
\end{figure}

}

\end{document}